\documentclass[english,11pt]{article}
\usepackage[german,french,english]{babel}
\usepackage[cp850]{inputenc}
\usepackage{latexsym,graphicx, fancybox}
\usepackage{graphicx}
\usepackage{caption}

\usepackage{amssymb, amsmath, amsfonts}
\usepackage{amsmath, amssymb}
\usepackage{color}
\usepackage{tikz}
\usepackage{latexsym}
\usepackage{tocloft}        

\usepackage{mathrsfs}

\usepackage{natbib}

\usepackage{etoolbox}
\apptocmd{\thebibliography}{\setlength{\itemsep}{0.01pt}}{}{}

\usepackage{microtype}
	\usepackage{amsthm}
	
	\newcommand{\dd}{\,\mathrm{d}}
	\newcommand{\R}{\mathbb{R}}
	
	\newcommand{\E}{\mathbb{E}}
	\renewcommand{\P}{\mathbb{P}}

	\usepackage{authblk}

	\usepackage{authblk}
	
	\newtheorem{Theorem}{Theorem}[section]
	\newtheorem{Definition}[Theorem]{Definition}
	\newtheorem{Proposition}[Theorem]{Proposition}
	\newtheorem{Assumption}[Theorem]{Assumption}
	
	\newtheorem{Lemma}[Theorem]{Lemma}
	
	\newtheorem{Remark}[Theorem]{Remark}
	\newtheorem{Example}[Theorem]{Example}
	
	\newtheorem{assumption}{Assumption}[section]

	\usepackage{thm-restate}
	\usepackage{float}
	\usepackage{hyperref}
	
\begin{document}
	
\selectlanguage{english}
\title{\bf On a Stationarity Theory for Stochastic Volterra Integral Equations with Affine Drift. %
}

\author{
	Emmanuel Gnabeyeu\\
	Laboratoire de Probabilit\'es, Statistique et Mod\'elisation (LPSM), UMR~8001,\\
	\textit{Sorbonne Universit\'e and Universit\'e Paris Cit\'e}, Paris, France. \\
	E-mail: \texttt{emmanuel.gnabeyeu\_mbiada@sorbonne-universite.fr}
	\and
	Gilles Pag\`es\\
	Laboratoire de Probabilit\'es, Statistique et Mod\'elisation (LPSM), UMR~8001,\\
	Sorbonne Universit\'e, case 158, 4, pl. Jussieu, F-75252 Paris Cedex 5, France. \\
	E-mail: \texttt{gilles.pages@sorbonne-universite.fr}
}
\date{September 1, 2025}
  \maketitle
  \renewcommand{\abstractname}{Abstract}
\begin{abstract}
	This paper investigate the properties of solutions to forward Stochastic Volterra Integral Equations (SVIEs for short) with affine drift, specifically their stationarity, both over a finite horizon and in the long run. 
	We demonstrate that it is possible to induce a \emph{fake stationary regime}, in the sense that all marginal distributions share the same expectation and variance. This phenomenon can be achieved either through explicit closed-form specifications of the deterministic initial condition \(\phi\) and the mean-reversion function \(\mu\) appearing in the drift, or by introducing a deterministic stabilizing factor \(\varsigma\) in the diffusion coefficient, associated with the kernel, while keeping the function \(\mu\) otherwise fully flexible.
	We further look at the $L^p$-confluence properties \(p>0\) of such processes as time goes to infinity, namely we investigate whether the marginals of solutions associated with different initial values become asymptotically confluent in $L^p$. We finally study the functional weak long-run asymptotics  for some classes of diffusion coefficients. More precisely, we establish that, in both settings, the time-shifted solutions of such SVIEs converge weakly, in the functional sense, toward a family of  $L^2$-stationary  processes sharing the same covariance function.
	These results are then applied to a class of Exponential-Fractional Stochastic Volterra Integral Equations driven by an $\alpha$-gamma fractional integration kernel,
	 in the particular case \(\alpha \in (0,1]\), which corresponds to the rough-path regime.
	  Building on these fake stationary Volterra processes, we finally introduce a family of stabilized Rough volatility models.
	
\end{abstract}


\section{Introduction} 
\subsection{Literature review}
The theory of stochastic Volterra integral equations (SVIEs) has its origins in the 1980s and has been widely developed since then. These equations which have recently attracted much attention in the mathematical finance community have
been introduced mostly with non-singular kernel for modelling in population dynamics, biology
and physics \cite{mohammed1998}, in order to generalize modelling to non-Markovian stochastic systems with
some memory effect. They were also motivated particularly by the physics of heat transfer \cite{gripenberg1990} and have undergone extensive mathematical study. Early investigations can be traced back to the seminal work of Berger et al. (see \cite{Berger1980a},\cite{Berger1980b})
who derived existence and uniqueness results for SVIEs driven by Brownian motion with Lipschitz continuous coefficients. These initial results were subsequently extended in various directions. 
For instance, \cite{Protter1985} generalized the existence and uniqueness results to SVIEs driven by right-continuous semimartingales and smooth kernels. An example of such a kernel is \(K(t,s) = (t-s)^{H-\frac{1}{2}}\), where \(H\) is known as the Hurst coefficient. Other studies focused on extensions that incorporated anticipative integrands, utilizing Skorokhod integration and Malliavin calculus (this was explored in \cite{Pardoux1990} and \cite{Alos1997}). 
\cite{Cochran1995} and \cite{CoutinD2001} focused on SVIEs with singular kernels. In a more recent contribution, \cite{Wang2008} proved the existence and uniqueness of solutions to SVIEs  with singular kernels and non-Lipschitz coefficients, utilizing a condition analogous to that of \cite{Yamada1971} for stochastic differential equations. Additionally, \cite{ZhangXi2010} examined SVIEs in Banach spaces with locally Lipschitz coefficients and singular kernels.

In the late 1990s, attempts were made within the financial community to incorporate long-memory effects into continuous-time stochastic volatility models. This shift was largely motivated by the need to capture persistent dependencies observed in financial markets, particularly through fractional Brownian motion (see \cite{CoutinD2001, ComteR1998}). Earlier studies, such as those by Comte and Renault \cite{ComteR1998}, found that \(H > 1/2\) was a key parameter in capturing long memory in volatility dynamics.
In the early 2000s, research shifted to Volterra equations with singular kernels that blow up as \(s \to t\) (i.e., \(K(t,s) \to +\infty\) as \(s \to t\) or \(H < 1/2\)). This development was further motivated by the empirical observation in~\cite{GatheralJR2018}, that volatility paths exhibit low H\"older regularity (\(H \approx 0.1\)).
As a result, there has been a resurgence of interest in SVIEs within mathematical finance, particularly with the rise of rough volatility models, as highlighted in the work of \cite{ElEuchR2018}. These models, which use the above kernel, naturally capture this feature, as their paths have a H\"older continuity exponent \(H\). Singular kernel Volterra equations also arise as limiting dynamics in models of order books via nearly unstable Hawkes processes (see \cite{JaissonR2016, EGnabeyeuR2025}). 
In both context, such processes are used to mimic Fractional Brownian motion-driven stochastic differential equations (SDEs). More specifically, within these frameworks, Volterra equations with fractional kernels \(K\) provide a more tractable alternative than SDEs involving stochastic integrals with respect to true \(H\)-fractional Brownian motions, which would otherwise require the use of ``high-order" rough path theory or regularity structures.  

The existence of stationary solutions and the study of the long-time behaviour of stochastic differential equations constitute fundamental problems in stochastic analysis.  In recent years, these questions have attracted considerable attention in the context of stochastic Volterra integral equations with time-independent coefficients. 
The ergodicity of SVIEs has been rigorously established for the particular class of affine Volterra processes in~\cite{friesen2022volterra,Jacquieretal2022}. The former relies on the analysis of Volterra-type Riccati equations induced by the specific affine structure, whereas the latter exploits the associated Markovian lifts. Both works establish the existence of stationary solutions for Volterra CIR-type processes, without providing an explicit characterization of their dynamics.
~\cite{Pages2024} recently introduced the notion of \emph{fake stationary regime} for SVIEs with general diffusion coefficients. This notion refers to the invariance, under time shifts, of the marginal distributions or of moments up to a prescribed order. In this paper, we investigate these properties for one-dimensional SVIEs with affine drift and convolution kernel on \([0,T]\) in $\mathbb{R}$, for arbitrary \(T>0\), of the following form
\begin{equation}\label{eq:Volterrameanrevert}
	\begin{cases}
		X_t = X_0\phi(t) +\int_0^t K(t-s)(\mu(s)-\lambda X_s)ds + \int_0^t K(t-s)\sigma(s,X_{ s})dW_s, \quad X_0\perp\!\!\!\perp W.\\
		X_0 : (\Omega, \mathcal{F}, \mathbb{P}) \to (\mathbb{R}, {\cal B}or(\mathbb{R})) \text{ is a given initial random variable}
	\end{cases}
\end{equation}
where \( \phi \) is a deterministic continuous function, typically normalized such that \( \phi(0) = 1 \), 
$\lambda>0$, $\mu :[0, T] \to \R$ is a Borel function, $\sigma : [0, T] \times \mathbb{R} \to \mathbb{R}$ is a Lipstchiz continuous function
and $K$ a deterministic kernel modeling the memory or hereditary structure of the system.  The process $(W_t)_{t \geq 0}$ is an $\mathbb{R}$-valued Brownian motion independent of $X_0$, both defined on a probability space $(\Omega, \mathcal{F}, \P)$ and $\mathcal{F}_t \supset \mathcal{F}_{t}^{X_0, W}$ a filtration satisfying the usual conditions.
Such equations~\eqref{eq:Volterrameanrevert} naturally arise in the modeling of random systems with memory effects and irregular behaviour, including in mathematical finance, physics, and biology.
The intrinsically non-Markovian nature of~\eqref{eq:Volterrameanrevert} presents substantial challenges for the analysis of its long-time asymptotic behavior. In particular, the classical ergodic theory for Markov processes, one of the fundamental tools in the study of SDEs, cannot be applied directly. Moreover, the solution to~\eqref{eq:Volterrameanrevert} generally does not belong to the class of semimartingales, thereby precluding the direct use of the classical It\^o calculus.  
However, the affine structure of the drift term provides access to powerful analytical tools, notably the \emph{resolvent} and the solution of the associated \emph{Wiener--Hopf equation}, which play a crucial role in overcoming these difficulties.

In \cite{Pages2024} a {\em fake stationary regime} for equations of the form~\eqref{eq:Volterrameanrevert} with fractional kernels is introduced without  \(\phi\), and exists under the conditions
\(\mu(t)=\mu_0,\quad dt\text{-a.e.}, \;
\mathbb{E}[X_0]=\frac{\mu_0}{\lambda}\) together with a separable diffusion coefficient featuring a deterministic time-dependent factor called {\em stabilizer} induced by the kernel. In this paper, we treat more general kernels and we show that the condition of a constant mean-reversion function can be removed by introducing a deterministic initial profile \(\phi\). More precisely, we establish that the same effect can be achieved by deriving explicit closed-form expressions for the pair \((\phi,\mu)\) dispensing with the stabilizer (section~\ref{subsec:fakestatioIntrin}), or by introducing a deterministic multiplicative factor in the state-dependent diffusion coefficient associated with the kernel and the mean-reversion function \(\mu\) which may remain completely arbitrary (section~\ref{subsec:fakestatioStabil}), in contrast to~\cite{Pages2024}.
\subsection{Our contribution}
In this paper, we investigate a weak form of stationarity for SVIEs with affine drifts and convolution kernels of the form~\eqref{eq:Volterrameanrevert}, referred to as a \emph{fake stationary regime} in the terminology of~\cite{Pages2024}, in the sense that the solution to~\eqref{eq:Volterrameanrevert} has either constant first- and second-order moments over time (\emph{fake stationary regime of type I}) or the same marginal distribution  at each time \(t\) in the Gaussian case (typically pseudo-Ornstein-Uhlenbeck process, called a \emph{fake stationary regime of type II}). This is achieved either through explicit closed-form specifications of both the initial function and the mean-reversion function, or, in the spirit of~\cite{Pages2024}, by introducing a deterministic multiplicative function in the Brownian convolution while, in contrast with that previous work, keeping the mean-reversion function fully flexible.
In this latter setting where the volatility coefficient is separable in time and state, the time-dependent multiplicative factor, referred to as the \emph{stabilizer}, 
serves to regulate or control the volatility of the process. This stabilizing function is characterized as the solution to an intrinsic convolution equation involving the mean-reverting function and the derivative of the resolvent associated with the  Volterra kernel.

Moreover, we establish the existence of limiting distributions. Formally, we prove that as \( u \to \infty \), the shifted process \( (X^u_t)_{t \geq 0} \), defined by \( X^u_t := X_{t+u} \), converges in law to a limiting continuous process \( X^{\infty} \). Unlike in~\cite{EGnabeyeuPR2025} (see also~\cite{friesen2022volterra}), this convergence does not imply that the limiting process is stationary (in the sense that its finite-dimensional distributions are invariant under time shifts). However, we prove that, under \emph{fake stationarity} regime, the limiting process is a weak \( L^2 \)-stationary.
Furthermore, since we do not characterize the dynamics of the limiting process, the notion of \emph{fake stationarity} provides a tractable alternative framework for analyzing both short- and long-term behaviors in settings where classical stationarity is either unavailable or analytically intractable.

In this work, we focus on one of the most prominent examples of Volterra kernels, namely the exponentially damped fractional (or Gamma) kernel \(	K(t):=K_{\alpha,\rho}(t)
=
\frac{t^{\alpha-1}}{\Gamma(\alpha)}e^{-\rho t},\quad \alpha\in\big(\frac12,\frac32\big),
\quad \rho>0,\) where \(\alpha := H-\frac12\).  This family also includes the fractional Riemann--Liouville kernel when \(\rho=0\)
and allows for a flexible incorporation of rough sample path properties (when $H\in\left(0,\frac12\right]$).
Although the debate on the empirical value of the Hurst index \(H\) remains controversial in the literature, we emphasize that the present paper however, focuses on the regime of ``rough path'' where \( H :=\alpha - \tfrac{1}{2} \in (0, \frac12] \)
and the complementary range of the Hurst coefficient, namely  \( H :=\alpha - \tfrac{1}{2} \in (\frac12, 1) \) corresponding to ``long-memory dynamics'', is addressed in a separate article~\cite{GnabeyeuPages2026b}.

From an applied perspective, this result may have important implications for volatility models widely used in mathematical finance. In particular, it suggests the possibility of introducing stabilized versions of such models (see e.g.~\cite{Gnabeyeu2026b,Gnabeyeu2026a}), where the dynamics driving the asset's (typically equity) volatility exhibit constant mean and variance over time. A key advantage of the stabilized formulation lies in its ability to overcome a well-known limitation of classical and rough Heston models \cite{Heston1993, el2019characteristic} driven by mean-reverting CIR or Volterra-CIR dynamics.
These models typically display two distinct regimes: a short-maturity regime, where the initial condition (deterministic value at the origin, often the long run mean) is prominent and the variance remains very small, and a long-term regime, which may correspond to the stationary distribution of the process.
In contrast, the stabilized model provides a unified and coherent framework that captures both short- and long-maturity behaviors within a single regime, thereby enabling robust and consistent fitting across the full term structure.

\subsection{Plan of the paper and Notations} 
The remainder of the paper is organized as follows: In Section~\ref{sec:background}, we review key properties of stochastic Volterra equations with convolution kernels, including results on existence, moment control, and a special focus on processes with affine drift. In this framework, several specific analytical tools become available, such as the \emph{resolvent} and the solution of the \emph{Wiener--Hopf equation}.
Section~\ref{sec:genfakestatio} investigates and derive the conditions ensuring that stochastic Volterra equations with affine drift of the form~\eqref{eq:Volterrameanrevert}  admit a weak stationary regime referred to as a \emph{fake stationary regime}, in the sense that the mean and variance of the solution remain constant over time. 
Next, in Section~\ref{sect:ExamplesFakeI-II}, we provide an example of a \emph{fake stationary regime } when the state-dependent diffusion coefficient is a trinomial function. Section~\ref{sect-LongRunB} is devoted to the analysis of the convergence and long-run behavior of these Volterra processes as time tends to infinity.
Specifically, we investigate, for such stabilized processes, the functional weak asymptotics of the time-shifted process \( (X_{t+s})_{s \geq 0} \) as \( t \to +\infty \), which turns out to be a weakly \( L^2 \)-stationary process.
Finally, in Section~\ref{sec:appl3}, we apply these results to the class of SVIEs driven by an \(\rho\)-exponential \(\alpha\)-fractional integration kernel with \(\alpha \in \left(\frac{1}{2},1\right)\), thereby capturing the {\em rough-path} and {\em short-memory} effects inherent to such Volterra equations.

\noindent {\sc \textbf{Notations.}} 

\smallskip
\noindent $\bullet$ Denote $\mathbb{T} = [0, T] \subset \mathbb{R}_+$, ${\rm Leb}_d$ the Lebesgue measure on $(\R^d, {\cal B}or(\R^d))$, $\mathbb H :=\R^d, $ etc.

\noindent $\bullet$ $\mathbb{X} := {\cal C}([0,T], \mathbb H) (\text{resp.} \; {\mathcal C_0}([0,T], \mathbb H))$ denotes the set of continuous functions (resp. null at 0)  from $[0,T]$ to $\mathbb H $ and ${\cal B}or(\mathbb{X})$ denotes the  Borel $\sigma$-field of $\mathbb{X}$ induced by the $\sup$-norm topology. 

\smallskip 
\noindent $\bullet$ For $p\in(0,+\infty)$, $L_{\mathbb H}^p(\P)$ or simply $L^p(\P)$ denote the set of  $\mathbb H$-valued random vectors $X$  defined on a probability space $(\Omega, {\cal A}, \P)$ such that $\|X\|_p:=(\E[\|X\|_{\mathbb H}^p])^{1/p}<+\infty$. For $f: E\to \R$, $\displaystyle \|f\|_{\sup}= \sup_{x\in E}|f(x)|$

\smallskip
\noindent$\bullet$ $\mathrm{Leb}_1$ (resp. $\mathrm{Leb}_{\mathbb{R}_+}$) denotes the Lebesgue measure on $\R$ (resp. $\mathbb{R}_+ = [0,+\infty)$).

\smallskip 
\noindent $\bullet$ For $f, g \!\in {\cal L}_{\R_+,loc}^1 (\R_+, {\rm Leb}_1)$, we define their convolution by $f*g(t) = \int_0^tf(t-s)g(s)ds$, $t\ge 0$.

\smallskip 
\noindent $\bullet$ For short, for any \(p\in[1,+\infty]\), we will denote by ${\cal L}^p(\R_+)$ (resp. ${\cal L}_{\mathrm{loc}}^p(\R_+)$) the space ${\cal L}^p(\mathbb{R}_+,{\rm Leb}_1)$ (resp. ${\cal L}_{\mathrm{loc}}^p(\mathbb{R}_+,{\rm Leb}_1)$).

\smallskip 
\noindent $\bullet$ For $f, g \!\in {\cal L}_{\R_+,loc}^2(\R_+, {\rm Leb})$ and $W$ a Brownian motion, we define  their stochastic convolution by 
\centerline{$
f\stackrel{W}{*}g = \int_0^t f(t-s)g(s) dW_s, \quad t\ge 0.
$}

\smallskip 
\noindent $\bullet$ For a random variable/vector/process $X$, we denote by $L(X)$ or $[X]$ its law or distribution. 

\smallskip 
\noindent $\bullet$ $X\perp \! \! \!\perp Y$  stands for independence of random variables, vectors or processes $X$ and $Y$.  

\smallskip 
\noindent $\bullet$ $\Gamma(a) = \int_0^{+\infty} u^{a-1} e^{-u} \, du, \; a > 0, \; 
\text{and} \quad 
B(a, b) = \int_0^1 u^{a-1} (1 - u)^{b-1} \, du, \quad a, b > 0.$
We will extensively use the classical identities:
$\Gamma(a + 1) = a \, \Gamma(a) 
\; \text{and} \; 
B(a, b) = \frac{\Gamma(a)\Gamma(b)}{\Gamma(a + b)}.$

\smallskip 
\noindent $\bullet$  Let $\mathrm{Pol}(\mathbb{R})$ denote the ring of all polynomials on $\mathbb{R}$ and $\mathrm{Pol}_n(\mathbb{R})$ the subspace consisting of polynomials of degree at most $n$.

\noindent Throughout this paper, we work on a filtered probability space $(\Omega,\mathcal{F},\mathbb{F}:=(\mathcal{F}_t^{X_0,W})_{t\in[0,T]},\mathbb{P})$
suppor-\\ting a \(1\)-dimensional Brownian motion \(W\), where \((\mathcal{F}_t^{X_0,W})_{t\in[0,T]}\) denotes the filtration generated by \(X_0\) and \(W\).
\section{Background on Volterra equations with convolution kernels}\label{sec:background}
We will assume that the process $(X_t)_{t\geq0}$ takes values in \(\R\) and we define $\mathbb{X} := {\cal C}([0,T], \mathbb R)$.
\begin{Definition}[Convolution kernel and Volterra equations] A kernel $K: \{(s,t)\!\in \R_+^2: 0\le s < t \} \to \R_+$ satisfying $
	\forall\, s,\, t\ge 0, \; s<t, \quad K(s,t)= K(0,t-s)
	$
	is called a {\em convolution kernel}. A Volterra equation based on a convolution kernel is called a {\em convolutive Volterra equation}.
\end{Definition}
To alleviate notations, we denote from now on $K(t):= K(0,t)$ so that $K(s,t) = K(t-s)$. For convenience we also extend the function $K: \R_+\to \R_+$  to  the whole real line by setting $K(t) = 0$, $t\le 0$.  

\subsection{Volterra processes with convolution kernels}\label{subsec:background}
A significant difference between regular diffusion processes and Volterra processes from a technical viewpoint  comes from the presence of the kernels which introduces some memory in the dynamics of the process, depriving us of the Markov property and usual tools of stochastic calculus.
We are interested in the convolutive stochastic Volterra equation: 
\begin{equation}\label{eq:Volterra}
	X_t=x_0(t) +\int_0^t K(t-s) b(s,X_s)ds+\int_0^t K(t-s)\sigma(s,X_s)dW_s,\quad t\ge 0.
\end{equation}
where $b:\mathbb{T}\times\R\to \R$, $\sigma:\mathbb{T}\times\R\to\R$ are Borel  measurable, $K\!\in {\cal L}^2_{loc}(\R_+)$  is a convolution kernel and $(W_t)_{t\ge 0}$ is a standard  Brownian motion independent from the $\R$-valued random variable $X_0$ both defined on a probability space $(\Omega,{\cal A}, \P)$. Let $(\mathcal F_t)_{t\ge 0}$ be a filtration (satisfying the usual conditions) such that $X_0$ is ${\cal F}_0$-measurable  and $W$ is an $({\cal F}_t)$-Brownian motion independent of $X_0$. \(X_0 \phi\) is thus a random function evolving deterministically for \(t>0\), i.e. \(X_0 \phi\) is ${\cal F}_0$-measurable.

\medskip
As will be shown in Example~\ref{Ex:FakeStatAutonome}~$(3)$, the mean-reverting function \(\mu\) appearing in~\eqref{eq:Volterrameanrevert}-\eqref{eq:Volterrameanrevert_} can be unbounded. Consequently, the associated drift does not satisfy uniform-in-time linear growth condition. This motivates a refined well-posedness conditions for the solution to~\eqref{eq:Volterra}, as highlighted below.
\begin{Assumption}[On Volterra Equations with convolution kernels]\label{assum:convol}
	Assume that  the kernel $K$ satisfies for every \(T>0\):
	{\small   
	\begin{equation}\label{eq:contKtilde}
		(\widehat {\cal K}^{cont}_{\widehat \theta})\;\;\exists\,\widehat\kappa_T< +\infty,\; \exists\, \widehat\theta\in (0,1],\;\forall\bar{\delta}\!\in (0,T],\; \widehat \eta(\delta) := \sup_{t\in [0,T]} \Big[\int_{(t-\bar{\delta})^+}^t \hskip-0,25cm K\big(t-u\big)^2 du\Big]^{\frac12}\le \widehat \kappa_T \,\bar{\delta}^{\,\widehat \theta} \; 
	\end{equation}
	\begin{equation}\label{eq:Kcont}
		({\cal K}^{cont}_{T,\theta}) \;\exists\, \kappa_{_T}< +\infty,\; \exists\, \theta >0,\; \forall\,\delta \!\in (0,T),\;
		\sup_{t\in [0,T]} \Big[\int_0^t |K(s+\delta)-K(s)|^2ds \Big]^{\frac 12} \le  \kappa_{_T}\,\delta^{\theta}.
	\end{equation}
	}
	Assume that the drift $b$ satisfies the following Lipschitz and linear growth conditions with time-dependent coefficients $h_{b}$ and $\ell_{b}$, which are locally square-integrable in time.
		\begin{align*}
	\noindent	(i) &\; \forall t \in [0,T], \; \forall x, y \in \mathbb{R}, \; |b(t, x) - b(t, y)| \le h_{b}(t)\, |x - y|\;\text{and}\; |b(t, x)| \leq \ell_{b}(t)\, (1 + |x|) \;\text{with}\; h_{b}, \ell_{b} \in {\cal L}^2([0,T]),
	\end{align*}
	Assume that the diffusion coefficient $\sigma$ satisfies the following Lipschitz-linear growth assumption uniformly in time with some real constants $ C_{\sigma,T}$ and $ C^\prime_{\sigma,T}$
	\begin{align*}
		(ii) &\; \forall t \in [0,T], \; \forall x, y \in \mathbb{R}, \; |\sigma(t, x) - \sigma(t, y)| \le C_{\sigma,T} |x - y| \;\text{and}\;  |\sigma(t, x)| \leq C^\prime_{\sigma,T} (1 + |x|), 
	\end{align*}
Finally, the initial random function \(x_0\) is $(\mathcal F_t)_{t\in[0,T]}$-adapted and satisfies for some $\delta > 0$ and $p>0$,
	\begin{align*}
		(iii) &\;  \mathbb{E}\big(\sup_{t \in [0,T]} |x_0(t)|^p \big) < +\infty, \quad \mathbb{E} |x_0(t') - x_0(t)|^p \leqslant C_{T,p} \Big( 1 + \mathbb{E}\big[ \sup_{t \in [0,T]} |x_0(t)|^p \big] \Big) |t' - t|^{\delta p}.
	\end{align*}
\end{Assumption}
It has been established in~\cite[Remark on Assumption 3.2 and Theorem 3.3]{GnabeyeuPages2026b} that, if $\ell_{b}$ and $h_{b}$ are defined on $\R_+$,
under assumption ~\eqref{assum:convol}, 
Equation~\eqref{eq:Volterra}  admits  a unique pathwise continuous solution on $\R_+$  
 satisfying (among other properties), \(\forall\, T>0, \;\forall\, a\in(0,\delta\wedge\widehat\theta\wedge\theta),\) 
\begin{equation}\label{eq:L^p-supBound}
		 \; \exists \,C_{a,p,T}>0,\quad \Big\|
	\sup_{0\le s<t\le T}
	\frac{|X_t-X_s|}{|t-s|^a}
	\Big\|_p^p + 	\Big\|
	\sup_{t\in[0,T]}|X_t|
	\Big\|_p^p
	<
	C_{a,p,T}
	\Big(
	1+ \mathbb{E}\big[\sup_{u\in[0,T]}|x_0(u)|^p\big]
	\Big),
	\vspace{-.1cm}
\end{equation}
for some positive constant
$C_{a,p, K, b,\sigma,T}$ that only depends on $a$, $p$, $K$, $b$,$\sigma$ and $T$.
This result extends to convolution-type kernels, the classical strong existence and uniqueness theorem for pathwise continuous solutions established in \cite[Theorem 1.1]{JouPag22}, itself refining \cite[Theorems 3.1 and 3.3]{ZhangXi2010}. In particular, it accommodates a general initial random function \(x_0\) and drift coefficients satisfying square-integrable, rather than uniform-in-time, Lipschitz and linear growth conditions.
\subsection{Fourier-Laplace transforms and Solvent core of a convolution kernel}\label{subsect:FLResolvent}
The Laplace transform is a valuable tool and we provide a brief overview here, as it is particularly effective for addressing the key equation $\eqref{eq:Volterra}$. We recall that the Laplace transform of a Borel function $f:\R_+\to \R_+$ is defined by

\centerline{$
\forall\, t\ge  0, \quad L_f(t)= \int_0^{+\infty} e^{-tu}f(u)du \!\in [0,\infty].
$} 
This Laplace transform is non-increasing and if $L_f(t_0)<+\infty$ for some $t_0\ge 0$, then $L_f(t)\to 0$ as $t\to +\infty$ by Lebesgue's dominated convergence theorem.
One can define the  Laplace transform of a Borel function $f:\R_+\to \R$ on $(0, +\infty)$ as soon as $L_{|f|}(t) <+\infty$ for every $t>0$  by the above formula. The Laplace transform can be extended to $\R_+$ as an $\R$-valued function  if $f \!\in {\cal L}^1_{\R_+}(\R_+)$.\\
Throughout this work, we adopt, as in \cite{Pages2024}, the following definition of the resolvent, which is also discussed in works such as \cite{Pruss1993} and offers a distinct perspective from that of the functional resolvent introduced in \cite{gripenberg1990}.

Let $K$ be a convolution kernel   satisfying ~\eqref{eq:contKtilde},~\eqref{eq:Kcont}  and $0<\int_{0}^tK(u)du$ for every $t>0$. 
For every $\lambda \!\in \R$,  the {\em resolvent or Solvent core} $R_{\lambda}$ associated to $K$ and $\lambda$ is defined as the unique solution -- if it exists --  to the deterministic Volterra equation
\begin{equation}\label{eq:Resolvent}
	\forall\,  t\ge 0,\quad R_{\lambda}(t) + \lambda \int_0^t K(t-s)R_{\lambda}(s)ds = 1,
\end{equation}
or, equivalently, written in terms of convolution, 
$R_{\lambda}+\lambda K*R_{\lambda} = 1.$
This equation is also known as \textit{resolvent equation} or \textit{renewal equation}. Its solution  always satisfies $R_{\lambda}(0)=1$ and admits the formal \textit{Neumann series expansion} \footnote{Recall that $ K^{1*} = K $ and $ K^{k*}(t) = \int_0^t K(t - s) \cdot K^{(k-1)*}(s) \, ds.$}:  
\begin{equation}\label{eq:Resolvent3}
	R_{\lambda} = \sum_{k\ge 0} (-1)^k \lambda^k (\mbox{\bf 1}*K^{k*}). 
\end{equation}
where \(K^{k*}\) denotes the \(k\)-th convolution of \(K\), with the convention in this formula, $K^{0*}= \delta_0$ (Dirac mass at $0$).
\medskip From now on we will assume that the kernel $K$ has a finite Laplace transform  \(L_K(t)<+\infty.\) Note that, as mentioned in \cite{Pages2024}, if the (non-negative) kernel $K$  satisfies 
\begin{equation}\label{eq:Kcontrol}
	0\le K(t)\le Ce^{bt }t^{a-1}\mbox{  for some }\; \; a,\, C>0, \; b\geq0. 
\end{equation}
 then, by induction  \(\mbox{\bf 1}*K^{*n}(t) \le C^n e^{bt}\frac{\Gamma(a)^n}{\Gamma(an+1)} t^{an},\) so that for such kernels, the above series~\eqref{eq:Resolvent3} is locally normally converging for every $t>0$ implying that the function   $R_{\lambda}$ is well-defined on $(0, +\infty)$. \\

\noindent{\bf Remark} 1. If $K$ is continuous, then the resolvent $R_{\lambda}$ is differentiable and one checks that $f_{\lambda}:=-R'_{\lambda}$ satisfies, for every  $t>0$, \(-f_{\lambda}(t) +\lambda \big( R_{\lambda}(0)K(t) - K *f_{\lambda}(t)\big)=0\)
		that is $f_{\lambda}$ is solution to the equation
		\begin{equation}\label{eq:flambda-eq}
			f_{\lambda} +\lambda K *f_{\lambda}=\lambda   K.
		\end{equation}
		2. Taking the Laplace transform from both side of the above equality~\eqref{eq:flambda-eq}, we have that :
		$L_{f_{\lambda}}(t)(1+\lambda L_K(t))=\lambda L_K(t) $, $t>0$. Consequently, \(L_{f_{\lambda}}(t) = \frac{\lambda L_K(t)}{1+\lambda L_{K}(t)}\)
		so that, for $\lambda \geq 0, $ $L_{f_{\lambda}}(t) \equiv 0$ if and only if 
		$L_K(t) \equiv 0$ i.e. if and only if $K=0$ by the injectivity of Laplace transform. Since $L_K(t) >0$ for any $t\in \R_+$, we deduce that $L_{f_{\lambda}}$ is $(0,1)$-valued.
		
		3. If $\displaystyle \lim_{t\to +\infty} R_{\lambda}(t) =0$ then, one also has that $\int_0^{+\infty} f_{\lambda}(t)dt = 1 -R_{\lambda}(+\infty) = 1$.Moreover, if $R_{\lambda}$ turns out to be non-increasing, then $f_{\lambda}$ is non-negative and satisfies $0\le f_{\lambda} \le \lambda K$, so that $f_{\lambda} $ {\em is a probability density}.

We also recall that, the one-dimensional \textit{functional resolvent of the first kind} for $f_\lambda$ in the terminology of \cite[Def. 5.5.1]{gripenberg1990} is defined as the unique measure $r_\lambda$ on $\R_+$ with values in $\R$ of locally finite variation satisfying 
\begin{equation}\label{eq:resolvent_first_kind}
	(f_\lambda \ast r_\lambda)(t) = (r_\lambda \ast f_\lambda)(t) = 1
\end{equation} for all $t > 0$. Remark that, if $f_\lambda$ is completely monotone such that $f_\lambda(t_0) > 0$ for some $t_0 > 0$, the existence of such a resolvent is guaranteed by \cite[Chapter 5, Theorem 5.4]{gripenberg1990}. In particular, this resolvent is of the form
\begin{align}\label{eq:frac_resolvent_first_kind}
	r_\lambda(dt) = \frac{\delta_0(dt)}{f_\lambda(0^+)} + \ell_\lambda(t)dt
\end{align}
with the convention $1/+\infty = 0$ and
where $\ell_\lambda \in L_{loc}^1(\R_+; \R)$ is completely monotone. By convoluting the equation~\eqref{eq:flambda-eq} defining \( f_\lambda \) with \( r_\lambda \), we immediately obtain that \( \lambda K * (r_\lambda-1)=1 \). Taking the Laplace transform yields
\begin{equation}\label{eq:LaplaceTr_lambda}
	L_{r_\lambda}(t)
	=
	\frac{1}{\lambda t L_K(t)}
	+\frac{1}{t}.
\end{equation}
\begin{Example}[Laplace transform and $\lambda-$ resolvent associated to the Exponential-fractional Kernel]\label{Ex:SolventGammaKernel}
	The Laplace transform associated to a kernel $K$ 
	reads, for $t>0$ $L_K(t):=\int_0^{+\infty}e^{-t u}K(u)du.$
	When K is the Gamma kernel $K_{b, \alpha , \rho }(t):= b e^{-\rho t}\frac{ t^{\alpha - 1}}{\Gamma(\alpha)}  \cdot \mathbf{1}_{(0,\infty)}(t)
	$, for \(b > 0, \alpha, \rho >0 \), then by
	introducing $v= u( t + \rho ) , \textit{ we have}$
	{\small
		\begin{equation}\label{eq:LTGammaKernel}
			L_{K_{b,\alpha,\rho}}(t)=\int_0^\infty be^{-(t+\rho)u}\frac{u^{\alpha-1}}{\Gamma(\alpha)}du=\frac{b(t+\rho)^{-\alpha}}{\Gamma(\alpha)}\int_0^\infty e^{-v}v^{\alpha-1}dv=b(t+\rho)^{-\alpha}.
		\end{equation}
	}
	Moreover, one checks that these kernels also satisfy~\eqref{eq:contKtilde} and~\eqref{eq:Kcont} for $\alpha >1/2$ (with $\widehat \theta=\theta= (\alpha-\frac 12)\wedge 1$, see~\cite{RiTaYa2020} or~\cite{JouPag22} or~\cite{GnabeyeuPages2026} along many others) and trivially~\eqref{eq:Kcontrol}.
	It follows from the easy identity \(K_{1,\alpha, \rho} * K_{1,\alpha^\prime, \rho} = K_{1,\alpha + \alpha^\prime, \rho}\) that \(K_{b,\alpha,\rho}^{(k*)}(t)	=
	b^k e^{-\rho t}
	\frac{t^{k\alpha-1}}{\Gamma(k\alpha)}.\) Therefore, the Neumann series expansion provided in equation~\eqref{eq:Resolvent3} shows that the resolvent
	reads:
	{\small
		\begin{equation}\label{eq:SolventGammaKernel}
			R_{\alpha, \rho, \lambda}(t) =(1*\delta_0)(t) + \sum_{k\ge 1} (-1)^k \lambda^k (\mbox{\bf 1}*K_{b,\alpha, \rho}^{(k*)}) = \mathbf{1}_{\mathbb{R}_+}(t) + \sum_{k \geq 1} (-1)^k (\lambda b)^k \int_0^t \frac{e^{-\rho s} s^{k\alpha -1}}{\Gamma(k\alpha)} \,ds.
		\end{equation}
	}
	Hence, if $\lambda > 0$, the function $f_{\alpha, \rho, \lambda}:= - R_{\alpha, \rho, \lambda} $ on $(0, +\infty)$ is given by:
	{\small
		\begin{equation}\label{eq:DerivSolventGammaKernel}
			f_{\alpha, \rho, \lambda}(t) = -\frac{d}{dt} R_{\alpha, \rho, \lambda}(t)
			= -\sum_{k \geq 1} (-1)^k (\lambda b)^k \frac{e^{-\rho t} t^{k\alpha -1}}{\Gamma(k\alpha)} = \lambda b e^{-\rho t} t^{\alpha - 1} \sum_{k \geq 0} (-1)^k (\lambda b)^k \frac{ t^{\alpha k }}{\Gamma(\alpha (k+1))}. 
		\end{equation}
	}
	\smallskip
	\noindent {\bf Remark:} Note that we recover the fractional kernel if $\rho = 1$ and the exponential kernel if $\alpha = 1$. In the latter case, \(K(t) = b e^{- \rho t} \mathbf{1}_{\mathbb{R}_+}(t)\), \(R_{1, \rho, \lambda}\)
	reads:
		\begin{align*}
			R_{1, \rho, \lambda}(t) 
			&= \mathbf{1}_{\mathbb{R}_+}(t) + \sum_{k \geq 1} (-1)^k (\lambda b)^k \int_0^t \frac{e^{-\rho s} s^{k - 1}}{\Gamma(k)} \, ds = \mathbf{1}_{\mathbb{R}_+}(t) + \int_0^t e^{- \rho s} \sum_{k \geq 1} (-1)^k (\lambda b)^k \frac{ s^{k - 1}}{k!} \, ds \\
			&= 1 - \lambda \int_0^t e^{-(\lambda b+\rho)s} \, ds
	 = 
		\begin{cases}
			\rho t + 1 & \text{if } \lambda b= -\rho \\
			\frac{\rho + \lambda b e^{-(\lambda b+\rho)t}}{\lambda b + \rho} & \text{if } \lambda b \ne -\rho
		\end{cases}
		\end{align*}
	If \(b=1\), \(\rho=0\) and \(\alpha\in (0,1)\), ~\cite[Proposition 5.1]{Pages2024} shows that \(f_{\alpha, \lambda }:=f_{\alpha,0, \lambda }\) is CM, so that since \(\lim_{t \to 0^+}f_{\alpha, \lambda }(t)=+\infty\), the resolvent of the first kind defined in~\eqref{eq:frac_resolvent_first_kind} reads $r_{\alpha, \lambda}(dt) =\ell_{\alpha, \lambda}(t)dt$. By Equation~\eqref{eq:LTGammaKernel}, \( L_{K_{\alpha}}(t)=t^{-\alpha} \) and
	therefore,  Equation~\eqref{eq:LaplaceTr_lambda} yields 
	$L_{r_{\alpha,\lambda}}(t)
	=
	\frac{t^{\alpha-1}}{\lambda}
	+\frac{1}{t}$. Hence, the first-kind resolvent of $f_{\alpha,\lambda}$ is given by 
	\begin{equation}\label{eq:FirstSolventFracKernel}
		r_{\alpha,\lambda}(dt)
	=
	\Big(1+
	\frac{t^{-\alpha}}
	{\lambda \Gamma(1-\alpha)}\Big) dt,
	\; t > 0.
	\end{equation}   
\end{Example}
\subsection{Application to the Wiener-Hopf equation}
\begin{Proposition}[Wiener-Hopf and Resolvent equations] \label{prop:W-H} Let $g, h: \R_+ \to \R$ be two  locally bounded Borel function, let $K \! \in \mathcal{L}^1_{loc}(\R_+)$ and let $\lambda \!\in \R$.  Assume that the $\lambda$-resolvent $R_{\lambda}$ of $K$ is differentiable on $(0, +\infty)$ with a derivative $R'_{\lambda}\!\in  \mathcal{L}^1_{loc} (\R_+) $, that both $R_{\lambda}$ and $R'_{\lambda}$ admit a finite Laplace transform on $\R_+$ and $\displaystyle \lim_{u\to +\infty} e^{-tu}R_{\lambda}(u) = 0$ for every $t>0$. Then,
	\begin{enumerate}
		\item[$(a)$]
	 The Wiener-Hopf equation
	$\forall\, t\ge 0, \quad x(t) = g(t) -\lambda \int_0^t K(t-s) x(s) ds$
	(also reading $x= g-\lambda K*x$)  has   a solution given by:
	
	\centerline{$
		\forall\, t\ge 0, \quad x(t) = g(t) +\int_0^t R'_{\lambda}(t-s)g(s)ds\; \quad \text{or equivalently,} \quad x= g- f_{\lambda}*g,
	$}
	where $f_{\lambda}= -R'_{\lambda}$. This solution is uniquely defined on $\R_+$ up to  $dt$-$a.e.$ equality.

 \item[$(b)$]	The integral equation
	$
		\forall\, t\ge 0, \quad x(t) = h(t) - \int_0^t R'_{\lambda}(t-s) x(s) ds \quad \text{where} \quad R'_{\lambda} = - f_{\lambda}$
	(also reading $x = h - R'_{\lambda} * x$) has a solution given by:
	
	\centerline{$
		\forall\, t\ge 0, \quad x(t) = h(t) + \lambda \int_0^t K(t-s) h(s) ds \; \quad \text{or equivalently,} \quad 	x = h + \lambda K * h.
	$}
	This solution is uniquely defined on $\R_+$ up to $dt$-a.e. equality.
	
	\end{enumerate}
\end{Proposition}
In Appendix~\ref{app:lemmata}, we provide a proof of this classical result (see also~\cite{gripenberg1990}) for the reader's convenience.
\section{Fake stationarity of a scaled stochastic Volterra Integral equation}\label{sec:genfakestatio}
From now we focus on  the special case of a {\em scaled} stochastic Volterra equation~\eqref{eq:Volterrameanrevert} associated to a convolution kernel $K:\R_+\to \R_+$ satisfying~\eqref{eq:contKtilde} and ~\eqref{eq:Kcont}:
We will assume that $\mu :\R_+ \to \R$ is a  
Borel function in \( {\cal L}^2_{loc}(\R_+)\) so that the drift $b(t,x) = \mu(t)-\lambda  x$ in~\eqref{eq:Volterra} is clearly Lipschitz continuous in $x$, and satisfy Assumption~\ref{assum:convol}~$(i)$.
 We also assume that the diffusion $\sigma: \mathbb{T}\times \mathbb{R} \to  \R$ is Lipschitz continuous in $x$, locally uniformly in $t\!\in \mathbb{T}$. 
Then, Equation~\eqref{eq:Volterrameanrevert} has a unique solution \( (X_t)_{t\geq0} \) adapted to \( \mathcal{F}^{X_0, W}_t \), starting from \( X_0 \in L^p(\mathbb{P}), p>0 \). This follows by applying the existence result of section~\ref{subsec:background}
to each time interval \( [0,T] \), \( T \in \mathbb{N} \), and gluing the solutions together under Assumption~\ref{assum:convol}. 

Note that under our assumptions, if \( p > 0 \) and \( \mathbb{E}[|X_0|^p] < +\infty \), then Equation~\eqref{eq:L^p-supBound}, 
 combined with the linear growth of \(\sigma\) uniformly in time, implies
that \( \mathbb{E}[ \sup_{t \in [0,T]}  |\sigma(t, X_{ t})|^p] < C'_T (1 + \|\phi\|^p_T\mathbb{E}[|X_0|^p]) < +\infty \) for every \( T > 0 \), enabling the unrestricted use of both regular and stochastic 
Fubini's theorems.
 
\noindent Sufficient conditions for interchanging the order of ordinary integration (with respect to a finite measure) and stochastic integration (with respect to a square integrable martingale) are provided in  \cite[Thm. 1]{Kailath_Segall},  and further details can be found in \cite[Thm. IV.65]{Protter} (see also \cite[ Theorem 2.6]{Walsh1986}, \cite[ Theorem 2.6]{Veraar2012}).

\medskip
We will always work under the following assumption:
\begin{assumption}[$\lambda$-resolvent $R_{\lambda}$ of the kernel]\label{ass:resolvent} Throughout the paper, we assume that the $\lambda$-resolvent $R_{\lambda}$ of the kernel $K$ is well-defined
	on $(0, +\infty)$ and satisfies the following for some $\lambda > 0$:
	\begin{equation}\label{eq:hypoRlambda}
		({\cal R}_\lambda)\quad
		\left\{
		\begin{array}{ll}
			(i) & R_{\lambda}(t) \text{ is } \text{differentiable on } \mathbb{R}_+,\; R_{\lambda}(0)=1 \text{ and } \lim_{t \to +\infty}R_{\lambda}(t) =a \in [0,1[,, \\
			(ii) &  
			f_{\lambda} \in {\cal L}^2(\mathbb{R}_+), f_\lambda(t) \neq 0\; dt-a.e., \text{ where we set } \quad f_{\lambda} := -R'_{\lambda} \text{ for } t > 0.
		\end{array} 
		\right.
	\end{equation}
\end{assumption}
\noindent Note that $L_{f_\lambda}(t) \in (0,1) $ and under assumptions \(({\cal R}_\lambda)\) $(i)$ and $(ii)$, \( f_{\lambda} \) is a \((1-a)\)-sum measure, i.e., \( \int_0^{+\infty} f_{\lambda}(s) \, ds = 1-a \).
In fact, 	

\centerline{$
	\int_0^{+\infty} f_{\lambda}(s)\,ds = [1 - R_{\lambda}(s)]_{s=0}^{s=+\infty} = -\lim_{s\to +\infty}R_{\lambda}(s) + R_{\lambda}(0) = 1-a
	$}
Moreover, condition \(({\cal R}_\lambda)~$(i)$\) also allows us to treat the case where the resolvent, bounded on \((0, +\infty)\), has a non-zero limit at infinity, unlike the framework of \cite{Pages2024}, thereby accommodating a broader class of kernels.
\begin{Proposition}[Wiener-Hopf transform]\label{prop:wiener_hopf}
	Let $\lambda > 0$ and let $\mu : \mathbb{R} \to \mathbb{R}$ be a Borel function. Assume the kernel $K$ satisfies the above assumptions $({\cal R}_\lambda)$, conditions~\eqref{eq:contKtilde} and \eqref{eq:Kcont} from Assumption \ref{assum:convol}. Then, the solution $(X_t)_{t \geq 0}$ of the Volterra equation ~\eqref{eq:Volterrameanrevert} also satisfies:
	\begin{equation}\label{eq:Volterrameanrevert_}
		X_t = X_0\big(\phi(t) - \int_0^t f_{\lambda}(t-s) \phi(s) \, ds \big) + \frac{1}{\lambda} \int_0^t f_{\lambda}(t-s) \mu(s) \, ds + \frac{1}{\lambda} \int_0^t f_{\lambda}(t-s) \sigma(s, X_{s}) \, dW_s.
	\end{equation}
	Conversely, any process satisfying \eqref{eq:Volterrameanrevert_} also satisfies the original Volterra equation \eqref{eq:Volterrameanrevert}. Thus, the two formulations are equivalent.
\end{Proposition}
\noindent {\bf Proof.}
	Equation~\eqref{eq:Volterrameanrevert} can be interpreted pathwise as a Wiener-Hopf equation with $x(t) = X_t(\omega)$ and
	
	\centerline{$
	g(t) = X_0(\omega)\phi(t) + (\mu*K)_t + \Big(K \stackrel{W}{*} \sigma(., X_{\cdot}(\omega))\Big)_t.
	$}
	This leads to the following expression for $X_t$:
	\begin{align*}
		X_t &= g(t) + \int_0^t R'_{\lambda}(t-s) g(s) \, ds =  X_0\phi(t) +  (\mu*K)_t + \big(K\stackrel{W}{*} \sigma(.,X_.)\big)_t \\
		&+ \int_0^t R'_{\lambda}(t-s)\Big[ X_0\phi(s)+ (\mu*K)_s+ \big(K\stackrel{W}{*} \sigma(.,X_.)\big)_s\Big] ds =X_0\phi(t) + (\mu*K)_t + \Big(K \stackrel{W}{*} \sigma(., X_{\cdot})\Big)_t \\
		&\quad + X_0 \int_0^t R'_{\lambda}(t-s)\phi(s) \, ds + \underbrace{\int_0^t R'_{\lambda}(t-s) (\mu*K)_s \, ds}_{(a)} + \underbrace{\int_0^t R'_{\lambda}(t-s) \Big(K \stackrel{W}{*} \sigma(., X_{\cdot})\Big)_s \, ds}_{(b)}.
	\end{align*}
	Using commutativity and associativity (via regular Fubini's theorem) of convolution, we obtain for $(a)$:
	{\small
	\begin{equation}\label{eq:(a)init}
		(a) = -f_{\lambda}*(\mu*K)_t = -\left((f_{\lambda}*K)*\mu\right)_t.
	\end{equation}
    }	
	Differentiating Equation~\eqref{eq:Resolvent} yields the identity \(-f_{\lambda} * K = \frac{1}{\lambda} f_{\lambda} - K,\) which, upon substitution into~\eqref{eq:(a)init}, leads to the following expression in \eqref{eq:a_and_b} for term~\((a)\).  
	For term $(b)$, owing to stochastic Fubini's theorem,  equation~\eqref{eq:flambda-eq} provides the below expression in \eqref{eq:a_and_b}.  
	\begin{equation}\label{eq:a_and_b}
		\begin{aligned}
			(a) &= \frac{1}{\lambda} (f_{\lambda} * \mu)_t - (K * \mu)_t, \quad
			(b) = \frac{1}{\lambda} \Big(f_{\lambda} \stackrel{W}{*} \sigma(\cdot, X_{\cdot})\Big)_t - \Big(K \stackrel{W}{*} \sigma(\cdot, X_{\cdot})\Big)_t.
		\end{aligned}
	\end{equation}
	Substituting~\eqref{eq:a_and_b} into~\eqref{eq:Volterrameanrevert} finally yields
	\[X_t = X_0(\phi(t)-(f_{\lambda} * \phi)_t) + \frac{1}{\lambda}(f_{\lambda} * \mu)_t + \frac{1}{\lambda} \Big(f_{\lambda} \stackrel{W}{*} \sigma(\cdot, X_{\cdot})\Big)_t,\]
	
One recognises hereinabove the Equation~\eqref{eq:Volterrameanrevert_} given in the Proposition.
The controverse is obtained by solving the corresponding Wiener-Hopf equation. We convolve both sides of Equation~\eqref{eq:Volterrameanrevert_} with the kernel \( K \), using regular and stochastic Fubini's theorem. Details are left to the reader. \hfill $\square$


\smallskip
\noindent {\bf Remark:} In the Markovian case, the Wiener--Hopf equation amounts to applying It\^o's lemma to the transformed process \( e^{\lambda t} X_t \). In fact, if $K(t) = \mathbf{1}$ in the volterra equation, then \(
	R_\lambda(t) = e^{-\lambda t} \) and  \( f_\lambda(t) = \lambda e^{-\lambda t},\)
	so that  the above computation corresponds to  It\^o' s Lemma applied to $e^{\lambda t} X_t$.
	
	
\begin{Definition}[Stationary of Order $p\geq1$ and Fake stationary regime of type I and II (see. \cite{Pages2024})]\label{def:stationarity_order}.
	\begin{enumerate}
		\item The process $(X_t)_{t\ge 0}$ starting from $X_0\!\in L^p(\P)$ for $p\geq1$,  exhibit a stationary regime of order p if:
		\centerline{$
			\forall\, t\ge 0,\quad \forall\, k \in \{1,...,p\}, \quad \E\,[X_t^k] = \E\,[X_0^k].
			$}
		
		\item The process $(X_t)_{t\ge 0}$ starting from $X_0\!\in L^2(\P)$,  exhibits a fake stationary regime of type I if:
			\begin{equation}\label{eq:fstatCond}
			\forall\, t\ge 0, \quad \E\, X_t = m_0,\quad   {\rm Var}(X_t)=v_0\ge 0\quad \mbox{and} \quad  \bar \sigma(t):= \E\, \sigma^2(X_t) = \bar \sigma_0^2\ge 0.
			\end{equation}
		\item The process $(X_t)_{t\ge 0}$ starting from $X_0$ 
		has a fake stationary regime of type II if \( (X_t)_{t \geq 0} \) has the same marginal distribution, i.e., \( X_t \overset{d}{=}X_0 \) for every \( t \geq 0 \). 
	\end{enumerate}
\end{Definition}

Before investigating the fake stationary regime of the ``scaled" stochastic Volterra equation~\eqref{eq:Volterrameanrevert}, we first determine under which conditions this equation has a constant first moment.

\subsection{Towards stationarity of the first moment.}\label{sec:statFirstMoments}
We begin by identifying the conditions under which the Volterra Equation~\eqref{eq:Volterrameanrevert} exhibits a constant mean; that is, when $\mathbb{E}[X_t] = \mathbb{E}[X_0]$ for all $t \geq 0$, assuming that $X_0 \in L^2(\mathbb{P})$.
We know that:
\centerline{$ 
\mathbb{E}\left[\left(\int_0^t f_{\lambda}(t-s) \sigma(s, X_{s}) dW_s\right)^2\right] = \int_0^t f_{\lambda}^2(t-s) \mathbb{E}[|\sigma(s,X_s)|^2] \, ds \leq C(1 + \| \phi\|_T^2\mathbb{E}[|X_0|^2]) \int_0^t f_{\lambda}^2(u) \, du < +\infty,
$}
which implies \(\mathbb{E}\left[\int_0^t f_{\lambda}(t-s) \sigma(s, X_{s}) dW_s \right] = 0.\)
Thus, we have
\begin{equation}\label{eq:espVolterramean}
	\forall t \geq 0, \quad \mathbb{E}[X_t] = (\phi(t)- (f_{\lambda} * \phi)_t) \mathbb{E}[X_0] + \frac{1}{\lambda} (f_{\lambda} * \mu)_t.
\end{equation}
Thus, $\mathbb{E}[X_t]$ is constant if and only if the following condition holds:
\begin{equation}\label{eq:condConstVolterramean}
	\forall t \geq 0, \quad \mathbb{E}[X_0] \Big(1 - \phi(t) + (f_{\lambda} * \phi)_t\Big) = \frac{1}{\lambda} \int_0^t f_{\lambda}(t - s)\, \mu(s)\, ds.
\end{equation}

\begin{Proposition}[Stationarity of the first moment]\label{prop:StationarityMean} 
	Let $\lambda > 0$ and \( (X_t)_{t \geq 0} \) be a solution to the scaled Volterra equation ~\eqref{eq:Volterrameanrevert} starting from $X_0\in L^1(\mathbb{P})$. Then the Volterra process \( (X_t)_{t \geq 0} \)
	has constant first moment $m_0$, if and only if either \(\mathbb{E}[X_0] =m_0= 0\) and \(\mu\equiv 0\; a.e.\) or  \(\mathbb{E}[X_0]=m_0 \neq 0\) and 
	{\small
	\begin{equation}\label{eq:CondMean}
	\forall\, t\ge 0\quad 	\phi(t) - (f_{\lambda} * \phi)(t) = 1 - \frac{(f_{\lambda} * \mu)(t)}{\lambda m_0} 
		\quad \text{so that}\quad  \phi(t) = 1 - \lambda \int_0^t K(t-s) \left( \frac{\mu(s)}{\lambda m_0} - 1 \right) \, ds.
	\end{equation} 
	}
\end{Proposition}
\noindent {\bf Proof.}
	\smallskip
	\noindent  {\sc Case~1} {\(m_0=0\):} In this case, we must first have \(\mathbb{E}[X_0] = 0\) and equation~\eqref{eq:condConstVolterramean} reads:
	\( (f_{\lambda} *\mu)_t= 0.\) 
	 Taking the Laplace transform and owing to assumption $({\cal R}_\lambda)~(ii)$, we \(\mu(t) = 0\; dt-a.e.\) 
	 In this case, from equation~\eqref{eq:espVolterramean}, we deduce that \(\forall t \geq 0, \quad \mathbb{E}[X_t] = 0 = m_0.\)
	
	\smallskip
	\noindent  {\sc Case~2} { \(m_0\neq0\):} In this case, \(X_0\) must be such that \(\mathbb{E}[X_0]=m_0 \neq 0\) equation~\eqref{eq:condConstVolterramean} reads:
	\(\phi(t)- (f_{\lambda} * \phi)(t) = 1 - \frac{1}{\lambda} (f_{\lambda} * \frac{\mu}{\mathbb{E}[X_0]})_t.\)
	We may read the above equation as a Wiener-Hopf equation with \(x(t)= \phi(t)\) and  \(h(t) =  1 - \frac{1}{\lambda} (f_{\lambda} * \frac{\mu}{\mathbb{E}[X_0]})_t.\)
	Then, applying the claim (b) of Proposition \ref{prop:W-H}, we get: \(\phi(t) = h(t) +\lambda (K*h)_t.\) That is:
	\begin{align*}
		\phi(t)  &=  1 - (f_{\lambda} * \frac{\mu}{\lambda\mathbb{E}[X_0]})_t + \lambda (K * 1)_t -  (K * f_{\lambda} * \frac{\mu}{\mathbb{E}[X_0]})_t = 1 -  \Big(\left(f_{\lambda} + \lambda K * f_{\lambda} \right)* \frac{\mu}{\lambda\mathbb{E}[X_0]} \Big)_t + \lambda (K * 1)_t \\
		&\overset{\eqref{eq:flambda-eq}}{=} 1 - (K * \frac{\mu}{\mathbb{E}[X_0]})_t + \lambda (K * 1)_t = 1 - \lambda \int_0^t K(t-s) \left( \frac{\mu(s)}{\lambda \mathbb{E}[X_0]} - 1 \right) \, ds.
	\end{align*} 
	\noindent Conversely, as \(\;\phi(t)- (f_{\lambda} * \phi)(t) = 1 - \frac{1}{\lambda} (f_{\lambda} * \frac{\mu}{\mathbb{E}[X_0]})_t,\) equation~\eqref{eq:espVolterramean} gives:
	\[	\forall t \geq 0, \quad \mathbb{E}[X_t] = \mathbb{E}[X_0] \Big( 1 - \frac{1}{\lambda} (f_{\lambda} * \frac{\mu}{\mathbb{E}[X_0]})_t\Big) + \frac{1}{\lambda} (f_{\lambda} * \mu)_t = \mathbb{E}[X_0]\]
\noindent Thus a necessary and sufficient condition for constant mean is given by Equation~\eqref{eq:CondMean}:
Then Equation~\eqref{eq:Volterrameanrevert_} can be rewritten as~\eqref{eq:ConstMean}
We can easily check that $\phi(0) = 1$.  \hfill $\square$

\medskip
\noindent {\bf Remark :} 
\noindent We deduce from the beginning of the section~\ref{sec:genfakestatio} that, the non-negative function defined by
\[t\longmapsto \xi^2(t):= \E\,\sigma^2(t,X_{t}), \quad t\ge 0,\]
is  locally bounded   on $\R_+$ (and thus $\xi^2\!\in L^1_{loc}( \R_+)$)  since $\sigma$ has at most linear growth in space, locally uniformly in $t\ge 0$. 
First noting that by assuming constant mean as in the above section, i.e. $  \forall\, t\ge 0,  \E\, X_t = \E\, X_0 = m_0 $, 
equation~\eqref{eq:Volterrameanrevert_} reads: 
\[X_t - m_0 = \big(X_0-m_0\big)(\phi - f_{\lambda} * \phi)(t) + \frac{1}{\lambda} \int_0^t f_{\lambda}(t-s) \sigma(s, X_{s}) \, dW_s.\]
By  It\^o's isomorphism and Fubini's Theorem
\begin{align*}
	{\rm Var}\Big( \int_0^t f_{\lambda}(t-s) \sigma(s, X_{s}) \, dW_s \Big) &= \E\Big[ \int_0^t f_{\lambda}(t-s) \sigma(s, X_{s}) \, dW_s \Big]^2 = \int_0^t f_{\lambda}(t-s)^2 \xi^2(s)ds  =(f_{\lambda}^2 * \xi^2)(t).
\end{align*}
Then, it follows from the above equation that: $\quad \forall\, t\ge 0, $ by setting $v_0  = {\rm Var}(X_0)$, we have:
{\small
	\[ 
	{\rm Var}(X_t)= \E\Big[( X_t - m_0)^2 \Big]=\E\Big[ (X_0 - m_0 )^2\Big] (\phi - f_{\lambda} * \phi)^2(t) + \frac { 1}{\lambda^2 }\int_0^t f_{\lambda}(t-s)^2 \xi^2(s)ds = v_0(\phi - f_{\lambda} * \phi)^2(t) + \frac { 1}{\lambda^2 }(f_{\lambda}^2 * \xi^2)(t) 
	\]
}
\vspace{-.5cm}
\begin{align}\label{eq:var}
	\text{i.e.} \quad	\forall\, t\ge 0, \quad   {\rm Var}(X_t)&= v_0(\phi - f_{\lambda} * \phi)^2(t) + \frac { 1}{\lambda^2 }(f_{\lambda}^2 * \xi^2)(t). 
\end{align}
\subsection{Intrinsic fake stationarity with autonomous diffusion coefficient}\label{subsec:fakestatioIntrin}
\noindent Now, assume a Lipschitz continuous time homogenous or autonomous volatility coefficient, i.e. \(\forall\, (t,x)\!\in \mathbb{T}\times \mathbb{R},\;\sigma(t,x) = \sigma(x).\)

\begin{Theorem}[Time homogenous $\sigma$]\label{thm:autonome} Set $\sigma(t,x)= \sigma(x)$ in~\eqref{eq:Volterrameanrevert}, that is the volatility coefficient is autonomous and assume $\mu$ in $L^1_{\mathrm{loc}}(\mathbb{R}_+)$ and  bounded at infinity (hence having   a well-defined finite Laplace transform on $(0,+\infty)$). Assume $X_0\!\in L^2(\P)$ and let $(X_t)_{t\ge 0}$ be a solution  to~\eqref{eq:Volterrameanrevert}  satisfying~\eqref{eq:fstatCond}.
%
\noindent Then, only two situations can occur:

\medskip 
\noindent $(i)$ If $\bar \sigma =0$ and either \(\mu\not\equiv 0\; a.e.\) or \(\mu\equiv 0\; a.e.\) and $\phi(t) \neq 1+\lambda\int_{0}^{t}K(s)ds\; dt-a.e.\;\forall\,t\geq0$, then $\sigma\big( m_0\big)=0$ and $X_t =  m_0$ $\P$-$a.s.$. However, if $\bar \sigma =0$, \(\mu\equiv 0\; a.e.\)and $\phi(t) = 1+\lambda\int_{0}^{t}K(s)ds\; dt-a.e.\;\forall\,t\geq0$, then \(\sigma(X_0)=0\) $\P$-$a.s.$ and the Volterra Equation~\eqref{eq:Volterrameanrevert_} reads \(X_t = X_0\) $\P$-$a.s.$.

\medskip 
\noindent $(ii)$ If $\bar \sigma >0$,  set $ \displaystyle c:=\frac{v_0}{\bar \sigma_0^2}$ be such that $c\geq \frac{1}{\lambda^2}
\int_0^\infty f_\lambda^2(s)\,ds$,  then, the process $(X_t)_{t \ge 0}$ has a stationary variance if and only if
the mean-reverting function $\mu$ satisfies either 
	 \(\mu\equiv 0\; a.e.\) (hence \(m_0 = 0\)) and 
\begin{equation}\label{eq:NullMeanVarAutonome}\phi(t)=h_{\lambda,c}(t)+\lambda\int_{0}^{t}K(t-s)h_{\lambda,c}(s)\,ds \quad \text{with}\quad h_{\lambda,c}(t):=\sqrt{1-\frac{1}{c\lambda^2}\int_0^t f_\lambda^2(s)\,ds}.
\end{equation}
or  \(\mu\not\equiv 0\; a.e.\) and is a solution (if any) of the first kind Volterra Equation:
	\begin{equation}\label{eq:conv_mu_general}
		(f_\lambda * \mu)(t) 
		=
		\lambda m_0
		\Big(
		1-h_{\lambda,c}(t)
		\Big),
		\qquad t \ge 0.
	\end{equation}
   In such case, there exists at most one unique solution \(\mu\) to~\ref{eq:conv_mu_general} and
	if the functional resolvent of the first kind \(r_\lambda\) in ~\eqref{eq:resolvent_first_kind} is well-defined, then $\mu$ in Equation~\eqref{eq:conv_mu_general} is uniquely determined and is given by
	\begin{equation}\label{eq:mu_general}
			\mu (t)
			:= \lambda m_0 \cdot \frac{d}{d t}
			\Big(r_\lambda
			*
			\big(
			1-h_{\lambda,c}
			\big)\Big) (t).
	\end{equation}
\end{Theorem}
\noindent{\bf Remark:} In the setting  \(\mu\not\equiv 0\; a.e.\) of Theorem~\ref{thm:autonome}~$(ii)$, the equations
\eqref{eq:Volterrameanrevert}--\eqref{eq:Volterrameanrevert_} starting from any $X_0 \in L^2(\mathbb{P})$ simplify to
{\small
	\begin{align}\label{eq:ConstMeanVar}
		X_t = X_0 -\big(X_0 - m_0\big) 	\big(
		1-h_{\lambda,c}(t)
		\big) +  \frac{1}{\lambda}\int_0^t f_{ \lambda}(t-s)\sigma(X_{s})dW_s.
	\end{align}
}
\noindent {\bf Proof of Theorem~\ref{thm:autonome} }
Assuming that the solution \( (X_t)_{t \geq 0} \) of the Volterra equation~\eqref{eq:Volterrameanrevert} starting from \( X_0 \in L^2(\mathbb{P}) \) satisfies~\eqref{eq:fstatCond}, then from equation ~\eqref{eq:var} together with the fact that, $\xi^2 = \bar \sigma^2 $ in this autonomous setting, we have:
   \[\forall\, t\ge 0, \quad   v_0 = {\rm Var}(X_t) =v_0(\phi - f_{\lambda} * \phi)^2(t) + \frac { 1}{\lambda^2 }(f_{\lambda}^2 * \bar \sigma^2)(t) = v_0(\phi - f_{\lambda} * \phi)^2(t) + \frac { \bar \sigma^2}{\lambda^2 }\int_0^t f_{\lambda}^2(s)ds\]
   or, equivalently
   {\small
   	\begin{equation}\label{eq:VolterraVar}
   		\forall\, t\ge 0, \quad v_0\big(1-(\phi - f_{\lambda} * \phi)^2(t) \big) =  \frac { \bar \sigma^2}{\lambda^2 }\int_0^t f_{\lambda}^2(s)ds.
   	\end{equation}
   }
   Consequently,\\
   \noindent 1. If $\bar \sigma^2=0$ and $\mu(t) \neq 0\; dt-a.e.\;\forall\,t\geq0$, we get $v_0=0$. In fact, since $	(\phi-f_\lambda * \phi)(t)=1-\frac{(f_\lambda * \mu)(t)}{\lambda m_0}$, 
   we have $ (\phi - f_{\lambda} * \phi)(t) \neq 1 $ (at least for infinitely many $t$). As a consequence, ${\rm Var}(X_t)=0$ for every $t\ge 0$. 
    If $\bar \sigma^2=0$, \(\mu\equiv 0\; a.e.\) and $\phi(t) \neq 1+\lambda\int_{0}^{t}K(s)ds\; dt-a.e.\;\forall\,t\geq0$, then 
    $(\phi-f_\lambda * \phi)(t)\neq 1 $ (at least for large enough $t$). Consequently, ${\rm Var}(X_t)=0$ for every $t\ge 0$.
   In both cases, we know by Proposition~\ref{prop:StationarityMean} that $\E\, X_t = \E\, X_0=  m_0$, it follows that  $X_t= m_0$ $\P$-$a.s.$. As a result, $\sigma^2 \big( m_0 \big)= \bar \sigma^2 = 0$.
   	\\
 However, if $\bar \sigma^2=0$, \(\mu\equiv 0\; a.e.\) and $\phi(t) = 1+\lambda\int_{0}^{t}K(s)ds\; dt-a.e.\;\forall\,t\geq0$, then $(\phi - f_{\lambda} * \phi)(t)=1\; dt-a.e.\;\forall\,t\geq0$ by Proposition~\ref{prop:W-H}~$(b)$ and \(\sigma(X_t)=0\) $\P$-$a.s.$.\\
   	\noindent 2. If $\bar \sigma >0$,
	rearranging Equation~\eqref{eq:VolterraVar} yields \((\phi - f_{\lambda} * \phi)^2(t)=1-\frac{1}{c\lambda^2}\int_0^t f_\lambda^2(s)\,ds.\)
	Since, the initial value in the left hand side is $ (\phi - f_{\lambda} * \phi)(0)=1,$
	continuity implies that the positive branch must be selected, that is:
	\begin{equation}\label{eq:ContMeanVarAutonome}
		(\phi - f_{\lambda} * \phi)(t)=\sqrt{1-\frac{1}{c\lambda^2}\int_0^t f_\lambda^2(s)\,ds}=:h_{\lambda,c}(t).
	\end{equation}
	Therefore, if \(\mu\equiv 0\; a.e.\) ( case \(m_0=0\), see Proposition~\ref{prop:StationarityMean}), then by Proposition~\ref{prop:W-H}~$(b)$, \(\phi(t)=h_{\lambda,c}(t)+\lambda\int_{0}^{t}K(t-s)h_{\lambda,c}(s)\,ds\) and consequently, Equation~\eqref{eq:NullMeanVarAutonome} holds. Moreover, if
	 $\mu(t) \neq 0\; dt-a.e.\;\forall\,t\geq0$, recalling from the constant mean condition in Proposition~\ref{prop:StationarityMean} that $	(\phi-f_\lambda * \phi)(t)=1-\frac{(f_\lambda * \mu)(t)}{\lambda m_0}$, Equation~\eqref{eq:ContMeanVarAutonome} implies that \(1-
	 \frac{(f_\lambda * \mu)(t)}{\lambda m_0} = h_{\lambda,c}(t) \)
	 leading to \((f_\lambda * \mu)(t) = \lambda m_0\big(1-h_{\lambda,c}(t)\big)\), which gives ~\eqref{eq:conv_mu_general}.
	Conversely, if either~\eqref{eq:NullMeanVarAutonome} or \eqref{eq:conv_mu_general} holds, then  substituting this
	expression into the variance identity in equation ~\eqref{eq:var} together with the fact that, here $\xi^2 = \bar \sigma^2 $,
	shows immediately that
	$\operatorname{Var}(X_t)=v_0$ for all $t \ge 0$.
	
	Suppose that \( \mu_1 \) and \( \mu_2 \) are two solutions of Equation~\eqref{eq:conv_mu_general}. Define \( \delta \mu := \mu_1-\mu_2 \). By subtracting the two equations satisfied by \( \mu_1 \) and \( \mu_2 \), we obtain \( (f_\lambda * \delta \mu)(t)=0 \) for all \( t\ge 0 \). Hence, \( f_\lambda * \delta \mu =0 \) in \( \mathcal{L}^1_{\mathrm{loc}}(\mathbb{R}_+) \).
	Applying the Laplace transform yields \( L_{f_\lambda}(z)\,L_{\delta\mu}(z)=0 \) for all \( z>0 \). Since \(L_{f_\lambda}(z)>0\) for every \(z>0\) (by Assumption~\(({\cal R}_\lambda)(ii)\)), it follows that \( L_{\delta\mu}(z)=0 \) for all \( z>0 \).
	By uniqueness of the Laplace transform, we conclude that \( \delta\mu =0 \) in \( \mathcal{L}^1_{\mathrm{loc}}(\mathbb{R}_+) \). Therefore, \( \mu_1=\mu_2 \), which proves that Equation~\eqref{eq:conv_mu_general} admits at most one solution.

	Finally, if \(r_\lambda \) defined in~\eqref{eq:resolvent_first_kind} exists, using the resolvent identity $(f_\lambda \ast r_\lambda)(t) = (r_\lambda \ast f_\lambda)(t) = 1$, 
	we convolve both sides of \eqref{eq:conv_mu_general} with $r_\lambda$ and
	obtain
	\begin{equation*}
	1 * \mu 
	=
	r_\lambda * (f_\lambda * \mu) =  \lambda m_0 \cdot r_\lambda
	*
	\Big(
	1-
	\sqrt{
		1-
		\frac{1}{c \, \lambda^2 }
		\int_0^{\,\cdot} f_\lambda^2(s)\,ds
	}
	\Big) ,
	\end{equation*}
	which gives \eqref{eq:mu_general}. At this stage, uniqueness also follows from the fact that
	the convolution operator $g \mapsto f_\lambda * g$ is inverted by
	convolution with $r_\lambda$. \hfill $\Box$

\begin{Example}[Markovian, Exponential and Fractional kernels]\label{Ex:FakeStatAutonome}

\smallskip	
\noindent {\em 1. Markovian case.}
		When $K\equiv1$, i.e \(\alpha=1\) in Example~\ref{Ex:SolventGammaKernel}, one checks using Equation~\eqref{eq:LaplaceTr_lambda} that the first-kind resolvent of \(f_{1,\lambda}:= \lambda e^{-\lambda t},
		\; t \ge 0\) reads $r_{1,\lambda}(dt)
		=\frac{\delta_0(dt)}{\lambda} + dt $. Therefore, Theorem~\ref{thm:autonome} yields
		\begin{equation}
		\mu(t)
		=
		g_{1,\lambda}(t)
		+
		\frac{1}{\lambda}g_{1,\lambda}'(t), \quad 
		\text{where}\quad
		g_{1,\lambda}(t)
		=
		\lambda m_0
		\Big(
		1-
		\sqrt{
			1-
			\frac{1}{2c\,\lambda}
			\bigl(1-e^{-2\lambda t}\bigr)
		}
		\Big), \quad c\geq \frac{1}{2\lambda}.
	   \end{equation}
	   One checks straightforwardly that \(\mu\) is continuous on \(\mathbb{R}_+\), that  \(\mu(0)=\frac{m_0}{2c}\),
	   and that $\lim_{t\to+\infty}\mu(t)
	   =
	   \lambda m_0
	   \big(
	   1-\sqrt{1-\frac{1}{2c\lambda}}
	   \big).$

		\smallskip	
\noindent {\em 2. Exponential Kernel case.} Let $\rho>0$ and consider the exponential kernel $K_{\rho}(t) = e^{-\rho t}\mbox{\bf 1}_{\R_+}(t)$, with
$f_{\rho, \lambda}$ defined by Equation~\eqref{eq:DerivSolventGammaKernel} of Example~\ref{Ex:SolventGammaKernel}.
We have \( L_K(t)=(t+ \rho)^{-1} \). Therefore by  Equation~\eqref{eq:LaplaceTr_lambda} $L_{r_{\rho, \lambda}}(t)
=
\frac{t+\rho}{\lambda\, t}
+\frac{1}{t}$ and hence, the first-kind resolvent of $f_{\rho,\lambda}$ is given by:
\begin{align*}r_{\rho, \lambda}(dt)&:= \Big(\int_{0}^{t} e^{-\rho u} L^{-1}_{1+
	\frac{s}
	{\lambda } }(u) du\Big) dt
	=\Big(\int_{0}^{t} e^{-\rho u}\delta_0(du) + \frac{1}
	{\lambda }\int_{0}^{t} e^{-\rho u} \delta_0^{\prime} du\Big) dt= \big(1+
	\frac{\rho}
	{\lambda }\big) dt,
	\; t > 0
\end{align*}
Consequently, by Theorem~\ref{thm:autonome}, which applies under the condition $c\geq \frac{1}{2(\lambda + \rho)}$, we obtain that the solution takes the form	$\mu(t)
=\big(1+
\frac{\rho}
{\lambda }\big)g_{\rho,\lambda}(t) $, with \(g_{\rho,\lambda}(t)
=
\lambda m_0
\Big(
1-
\sqrt{
	1-
	\frac{1-e^{-2(\lambda+\rho)t}}{2c(\lambda+\rho)}
}
\Big).\)
It is straightforward that \(\mu\) is continuous on \(\mathbb{R}_+\). Moreover, one checks that \(\mu(0)=0\), and that  $\mu_\infty:=\lim_{t\to+\infty}\mu(t)
=
(\lambda+\rho)m_0
\Big(
1-
\sqrt{
	1-
	\frac{1}{2c(\lambda+\rho)}
}
\Big).$

 \smallskip	
 \noindent {\em 3. Fractional case.}
		Let $\alpha \in (0,1)$ and consider the fractional integration kernel $K_{\alpha}(t) = \frac{u^{\alpha-1}}{\Gamma(\alpha)} \mbox{\bf 1}_{\R_+}(t)$, with
		$f_{\alpha,\lambda}:=f_{\alpha,0,\lambda}$ defined by Equation~\eqref{eq:DerivSolventGammaKernel} of Example~\ref{Ex:SolventGammaKernel}. 
Consequently, Theorem~\ref{thm:autonome}
		applies with $f_\lambda := f_{\alpha,\lambda}$, provided that $c\geq \frac{1}{\lambda^2} \int_0^\infty f_{\alpha,\lambda}^2(s)\,ds$.
		In view of Equation~\eqref{eq:FirstSolventFracKernel}, it follows that
		\begin{equation*}
			\mu(t)
			=
			g_{\alpha,\lambda}(t)
			+ \frac{1}{\lambda} D^{\alpha} g_{\alpha,\lambda}(t)
			 \quad \text{with}\quad
			g_{\alpha,\lambda}(t)
			=
			\lambda m_0
			\Big(
			1-
			\sqrt{
				1-
				\frac{1}{c\,\lambda^2}
				\int_0^t f_{\alpha,\lambda}^2(s)\,ds
			}
			\Big)
		\end{equation*}
		 and \( D^{\alpha} := \frac{d}{dt} I^{1 - \alpha} \) where \( D^\alpha \) and \( I^{1 - \alpha} \) denote, respectively, the Riemann--Liouville fractional derivative of order \( \alpha \), and the Riemann--Liouville fractional integral of order \( 1 - \alpha \) (see \cite[Chapter~2]{samko}). \footnote{Recall that the fractional integral of order \( r \in (0, 1] \) of a function \( f \) is \(I^r f(t) = \frac{1}{\Gamma(r)} \int_0^t (t - s)^{r - 1} f(s) \, ds,\)
			while the fractional derivative of order \( r \in [0, 1) \) is defined as \(	D^r f(t) = \frac{1}{\Gamma(1 - r)} \frac{d}{dt} \int_0^t (t - s)^{-r} f(s) \, ds,\)
			whenever the integrals exist.}
			
			\smallskip\noindent
		Moreover, defining \(C_{\alpha,\lambda}:= \frac{m_0 \lambda}{2c\,(2\alpha-1)\,\Gamma(\alpha)^2}\) and using $1-\sqrt{1-x}=\frac{x}{2}+ O\!\left(x^{2}\right),
		\; x\to 0,$  one checks that 
		$	g_{\alpha,\lambda}(t)=
		C_{\alpha,\lambda} t^{2\alpha-1}
		+ O\!\left(t^{4\alpha-2}\right),
		\; t\to 0.$  it follows that \(D^\alpha g_{\alpha,\lambda}(t)
		=
		C_{\alpha,\lambda}
		\frac{\Gamma(2\alpha)}{\Gamma(\alpha)}
		t^{\alpha-1}
		+
		O\!\left(t^{3\alpha-2}\right),
		\; t\to 0.\) As \(2\alpha-1>0\), one has $O\!\left(t^{3\alpha-2}\right)
		=
		o\!\left(t^{\alpha-1}\right),$
		so that $D^\alpha g_{\alpha,\lambda}(t)
		\sim
		C_{\alpha,\lambda}
		\frac{\Gamma(2\alpha)}{\Gamma(\alpha)}
		t^{\alpha-1}$ and consequently
	$\mu(t)
		\underset{0^+}{\sim} 
		\frac{C_{\alpha,\lambda}\,\Gamma(2\alpha)}{\lambda\,\Gamma(\alpha)}
		\, t^{\alpha-1}=\frac{ \alpha\,m_0\, B(2\alpha,1-\alpha)}{2c\,(2\alpha-1)\,\Gamma(\alpha)^2\,\Gamma(1-\alpha)}
		\, t^{\alpha-1}.$ 
		One checks that  \(\mu\) is continuous on \((0,+\infty)\) and that $\lim_{t\to+\infty}\mu(t)
			=
			\lambda m_0
			\big(
			1-\sqrt{1-\frac{1}{c\,\lambda^2}
				\int_0^\infty f_{\alpha,\lambda}^2(s)\,ds}
			\big)$. Moreover, one checks that, and \(\mu \in L_{\mathrm{loc}}^{2}(\mathbb{R}_+)\) for every \(\alpha>\frac{1}{2}\).
\end{Example}
\smallskip
The discussions above (Theorem~\ref{thm:autonome} and Example~\ref{Ex:FakeStatAutonome}) show that enforcing both constant mean and constant variance for Equation~\eqref{eq:Volterrameanrevert} is rather restrictive, as it essentially forces explicit closed-form structures on both the initial function \(\phi\) and the mean-reverting function \(\mu\) in the drift term. To relax this rigidity, and following the idea of~\cite{Pages2024}, we introduce a deterministic stabilising factor in the diffusion coefficient while allowing the mean-reversion function \(\mu\) to remain fully flexible. 

\subsection{Fake stationarity using a stabilizer}\label{subsec:fakestatioStabil}
In this setting, we shall therefore work with Equations~\eqref{eq:Volterrameanrevert}--\eqref{eq:Volterrameanrevert_} with \(\phi\) given by~\eqref{eq:CondMean} as our baseline Volterra Equation with constant mean.
From this point onward, we assume that the volatility coefficient $\sigma(t,x)$ is time-dependent or inhomogenous and factorises as:
\begin{equation}\label{eq:InhomogeneousCoefs}
\forall\, (t,x)\!\in \mathbb{T}\times \mathbb{R},\qquad\sigma(t,x) = \varsigma(t)\sigma(x)\quad \varsigma(t), \sigma(x) > 0.
\end{equation}
where \(\varsigma\) is a Borel function to be specified later, such that $\varsigma^2\in L^1_{\mathrm{loc}}(\mathbb{R}_+)$ and  bounded at infinity (hence having   a well-defined finite Laplace transform on $(0,+\infty)$) and \(\sigma\) a Lipschitz continuous function.

We assume that the kernel $K$ satisfies assumptions ~\eqref{eq:contKtilde} and ~\eqref{eq:Kcont} of Assumption \ref{assum:convol} and \(\sigma\) is {\em Lipschitz continuous}. If  \(\varsigma\) is a bounded Borel function, the scaled Volterra equation ~\eqref{eq:Volterrameanrevert} has a unique \((\mathcal{F}^{X_0, W}_t)_{t>0}\)-adapted pathwise continuous solution starting from \(X_0 \in L^2(\P)\) independent of W (still owing to the result in section~\ref{subsec:background}).

\medskip
In this time-inhomogeneous diffusion setting, we investigate the conditions under which the fake stationary property holds for the solution \( (X_t)_{t \geq 0} \) of the Volterra equation~\eqref{eq:Volterrameanrevert} starting by \( X_0 \in L^2(\mathbb{P}) \), namely
\begin{equation}\label{eq:FakeStats}
	\forall\, t\ge 0, \;  \mathbb{E}[X_t] = m_0, \quad \text{Var}(X_t) = v_0 \geq 0, \quad \text{and} \quad \bar{\sigma}^2(t) = \mathbb{E}[\sigma^2(X_{t})]=\bar \sigma^2_0 \geq 0.
\end{equation}
The theorem below shows what are the consequences of these three constraints in this settings.
	Suppose the following conditions hold for all \( t \geq 0 \):
\begin{Theorem}[Time-dependent or inhomogenous diffusion coefficient $\sigma$] \label{thm:timeautonome}
	Let Assumption~\eqref{ass:resolvent} be in force. Set \(\phi\) as in~\eqref{eq:CondMean} and let the diffusion coefficient \( \sigma(t, x)\) in equation~\eqref{eq:Volterrameanrevert} be multiplicatively time-inhomogeneous, namely of the form~\eqref{eq:InhomogeneousCoefs}. Assume that \( X_0 \in L^2(\mathbb{P}) \) with and let \( (X_t)_{t \geq 0} \) be a solution to the scaled Volterra equation ~\eqref{eq:Volterrameanrevert}-\eqref{eq:Volterrameanrevert_} starting from $X_0$ and satisfying~\eqref{eq:fstatCond}.
	Then, a necessary condition for the identities~\eqref{eq:fstatCond} to hold is that \(\mathbb{E}[X_0]=m_0 \neq 0\) and the triplet \( (v_0, \bar{\sigma}_0^2, \varsigma) \) satisfies the following {\em functional equation }:
		\begin{equation}\label{eq:VolterraStabilizer}
		\textit{($E_{\lambda, c}$)}: \quad\forall\, t\ge 0, \quad c \lambda^2\Big(1-\big(1 - \frac{(f_{\lambda} * \mu)_t}{\lambda m_0}\big)^2 \Big) =  (f_{\lambda}^2 * \varsigma^2)(t) \quad  \textit{with} \quad c = \frac { v_0 }{\bar \sigma_0^2} \quad  \textit{and} \quad \varsigma := \varsigma_{\lambda,c}.
	\end{equation}
\end{Theorem}
\noindent {\bf Remark:} If for every \(t>0\), \(\mu(t)=\mu(0)=\lambda m_0\), then owing to Equation~\eqref{eq:CondMean}, $\phi(t) = 1, \, \forall t \geq 0$ and we recover the framework studied in \cite{Pages2024} with a simplified equation $\textit{($E_{\lambda, c}$)}$.


\smallskip
\noindent {\bf Proof of Theorem \ref{thm:timeautonome}.}
 Assuming that conditions~\eqref{eq:FakeStats} hold, from Equation ~\eqref{eq:var} with $\xi^2 = \varsigma^2\bar \sigma^2 $ and the assumption of the theorem :
	\[
	\forall\, t\ge 0, \quad v_0 = {\rm Var}(X_t) =v_0\big(1 - \frac{(f_{\lambda} * \mu)_t}{\lambda m_0}\big)^2 + \frac { 1}{\lambda^2 }(f_{\lambda}^2 * \bar \sigma^2 \varsigma^2)(t) = v_0\big(1 - \frac{(f_{\lambda} * \mu)_t}{\lambda m_0}\big)^2 + \frac { \bar \sigma^2}{\lambda^2 }(f_{\lambda}^2 * \varsigma^2)(t)\]
	or, equivalently
	
	\centerline{$
		\forall\, t\ge 0, \quad v_0\big(1-\big(1 - \frac{(f_{\lambda} * \mu)_t}{\lambda m_0}\big)^2 \big) =  \frac { \bar \sigma^2}{\lambda^2 }(f_{\lambda}^2 * \varsigma^2)(t).
	$} 
\noindent This completes the proof. \hfill $\square$


\begin{Definition}
	We will call the stabilizer (or corrector) of the scaled stochastic Volterra equation ~\eqref{eq:Volterrameanrevert} a (locally) bounded Borel function \( \varsigma = \varsigma_{\lambda, c} \), which is a solution(if any) to the functional equation~\eqref{eq:VolterraStabilizer}
\end{Definition}
\noindent{\bf Remarks}
	1. Note that ~\eqref{eq:VolterraStabilizer} has a solution \( \varsigma_{\lambda, c} \) for some \( c > 0 \) if and only if it has a solution \( \varsigma_{\lambda, 1} \) when \( c = 1 \), and \( \varsigma_{\lambda, c} = \sqrt{c} \varsigma_{\lambda, 1} \). Hence, ($E_{\lambda, c}$) can be replaced by ($E_{\lambda, 1}$) denoted ($E_{\lambda}$) for simplicity.
2. In practice, uniqueness is not the delicate issue (see, for instance, the Lemma below). The crucial requirement is the existence and positivity of the solution.

\begin{Assumption}[On the stabilizer]\label{ass:on_stabilizer}
	Let \(\lambda>0\). There exists a  
	positive bounded Borel solution \( \varsigma_{\lambda} \) on \( [0, +\infty) \) of the equation $ \textit{($E_{\lambda}$)}: \quad
	\forall\, t> 0, \quad \lambda^2\big(1-\big(1 - \frac{(f_{\lambda} * \mu)_t}{\lambda m_0}\big)^2 \big) =  (f_{\lambda}^2 * \varsigma^2)(t)$.
\end{Assumption}

\begin{Lemma}[Uniqueness  and limit of the solution $\varsigma_{\lambda, c}^2$ of Equation \textit{($E_{\lambda, c}$)}]\label{lem:asymptotique}
For fixed $c>0$, under assumptions~\ref{ass:resolvent} and~\ref{ass:on_stabilizer}, there exists at most one solution 
$\varsigma^2_{\lambda, c}$ in $\mathcal{L}^1_{\text{loc}}(\R_+)$ to Equation \textit{($E_{\lambda, c}$)} in~\eqref{eq:VolterraStabilizer}.
Furthermore, if $ \ell_\infty:=\lim_{t\to +\infty} \big(1-\frac{(f_\lambda * \mu)(t)}{\lambda m_0}\big)$ exists and is finite, then the unique solution of \textit{($E_{\lambda, c}$)} in~\eqref{eq:VolterraStabilizer} satisfies

	\centerline{$
	\lim_{t\to+\infty}\varsigma^2_{\lambda, c}(t)
	= \frac{ c \lambda^2(1-\ell_\infty^2) }{\|f_{\lambda}\|^2_{\mathcal{L}^2(\R_+)} }.$}

	
	

\end{Lemma}
For clarity and conciseness, the proof of Lemma~\ref{lem:asymptotique} is deferred to Appendix \ref{app:lemmata}, where the main technical results are presented.
\subsection{Examples of fake stationary regimes of type~I and~II}\label{sect:ExamplesFakeI-II}
Motivated by Example~\ref{Ex:FakeStatAutonome}, in which the drift function \(\mu\) admits a finite limit at infinity, we shall work under the following assumption on the mean-reversion function.
\begin{assumption}[On the mean-reversion function \(\mu\)]\label{ass:MeanRevert}  Let \(m_0\neq0\) and assume that the mean-reverting function \(\mu\) is a 
	$ C^1$-function such that  \(\mu\not\equiv 0\; a.e.\),  
	$\lim_{t\to +\infty} \mu(t) = \mu_{\infty} \in \mathbb{R}$ and $\frac{\mu_\infty}{m_0}\geq 0$.
\end{assumption}
\begin{Lemma}\label{lem:on_mu}
	Assume that both assumptions \(({\cal R}_\lambda)\) and ~\ref{ass:MeanRevert} hold, then \(\lim_{t\to +\infty} \int_0^t f_{\lambda}(t-s) \mu (s)ds = \mu_{\infty} (1-a) \).
\end{Lemma}
\noindent {\bf Remarks:} 1. Let $\phi$ be given by Equation~\eqref{eq:CondMean}. Under Assumption~\ref{ass:MeanRevert} and
following Lemma~\ref{lem:on_mu}, we obtain that \(\ell_\infty\) defined in Lemma~\ref{lem:asymptotique} is given by
\begin{equation}\label{eq:liml}
	\ell_\infty=\lim_{t\to +\infty}  (\phi-f_\lambda * \phi)(t) =\lim_{t\to +\infty} 
	1-
	\frac{(f_\lambda * \mu)(t)}{\lambda m_0}=
	1 - \frac{1-a}{\lambda m_0} \mu_\infty.
\end{equation}
Note that by taking the limit in both side of the equality~\eqref{eq:condConstVolterramean} together with Lemma~\ref{lem:on_mu}, we have:
\[ \mathbb{E}[X_0](1-\ell_\infty) := \lim_{t \to +\infty} \mathbb{E}[X_0] \Big(1 -  (\phi(t)- (f_{\lambda} * \phi)_t) \Big) = \lim_{t \to +\infty} \frac{1}{\lambda}(f_{\lambda} * \mu)_t =: \frac{\mu_\infty}{\lambda}(1-a) \]
 so that, \(\mathbb{E}[X_0] = \frac{1-a}{1-\ell_\infty}\frac{\mu_\infty}{\lambda}=m_0\) in the {\em fake stationary regime}.
Moreover, under Assumption~\ref{ass:MeanRevert}, we have $\ell_\infty < 1$.
The condition \(\ell_\infty=0\) (i.e. \(m_0:
=(1-a)\frac{\mu_\infty}{\lambda}\)) will ensure a \(L^p-\)confluence property in the long-time regime, as established in Proposition~\ref{prop:confluenceLp} further on. 

\smallskip
\noindent 2. In contrast with the fractional kernel case, if \(K=K_{\alpha,\rho},\;\rho>0\) (Gamma kernel), then \(a:=\lim_{t \to +\infty} R_{\alpha,\rho,\lambda}\neq 0\) and the limit of the function \(\phi\) in Equation~\eqref{eq:CondMean} can be computed explicitly under assumption~\ref{ass:MeanRevert}. In fact, owing to the Tauberian final value theorem, combined with the Laplace transform associated with~\eqref{eq:Resolvent}, one verifies that
\[\phi_\infty := \lim_{t\to +\infty} \phi(t)=1- \lambda \lim_{s \to 0} L_{K_{\alpha, \rho}}(s) \lim_{t\to +\infty}\left( \frac{\mu(s)}{\lambda m_0} - 1 \right) = 1- \lambda \frac{1}{\lambda} (\frac{1}{a}-1) \big(\frac{\mu_\infty}{\lambda m_0}-1\big)= \frac{1}{a}-\frac{1-a}{a} \frac{\mu_\infty}{\lambda m_0}\]
3. It then follows from Lemma~\ref{lem:on_mu} that $
\lim_{t \to +\infty} \phi(t) - (f_{\lambda} * \phi)(t) = \phi_\infty \,a = 1 - \frac{1-a}{\lambda m_0} \mu_\infty=:\ell_\infty
.$
Consequently, the mean of the initial random variable \(X_0\) is fully determined by \(m_0:=\frac{1-a}{1-\ell_\infty}\frac{\mu_\infty}{\lambda}\) in the fake stationarity regime, once \(\ell_\infty\) is specified,
in agreement with the framework introduced in~\cite[Section~2.1]{Pages2024} where \(\ell_\infty=a=0\).


\medskip
In order to treat both the autonomous and time-dependent diffusion cases within a unified framework, we assume 
from this point onward in the paper that Assumption~\ref{ass:MeanRevert} holds and $\phi$ is given by Equation~\eqref{eq:CondMean} so that the Volterra SDE~\eqref{eq:Volterrameanrevert}-\eqref{eq:Volterrameanrevert_} starting from any $X_0 \in L^2(\mathbb{P})$ is rewritten as
\begin{equation}
	X_t
	=
	X_0
	-
	\frac{X_0-m_0}{\lambda m_0}
	\int_0^t f_{\lambda}(t-s)\mu(s)\,ds
	+
	\frac{1}{\lambda}
	\int_0^t f_{\lambda}(t-s)\varsigma(s)\sigma(X_s)\,dW_s,
	\quad X_0\perp\!\!\!\perp W.
	\label{eq:ConstMean}
\end{equation}
and is governed by the following master relation
\begin{equation}\label{eq:MasterEq}
	\Big(1-\frac{(f_\lambda * \mu)(t)}{\lambda m_0}\Big)^2=h_{\lambda,c}^2(t) \quad \text{with $h_{\lambda,c}$ redefined as} \quad  h_{\lambda,c}(t):=\sqrt{1-\frac{1}{c\lambda^2}(f_{\lambda}^2 * \varsigma^2)(t)}.
\end{equation}
where we furthermore adopt the convention
\begin{equation}\label{eq:VolterraStabilizerCont}
	\varsigma = \begin{cases}
		1, & \text{in the autonomous diffusion setting , i.e., } \quad \sigma(t,x) =\sigma(x), \text{ Theorem~\ref{thm:autonome}} \\
		\varsigma_{\lambda,c}, & \text{in the time-dependent diffusion setting, }  \; \sigma(t,x) = \varsigma_{\lambda,c}(t)\,\sigma(x),\;  \text{ Theorem~\ref{thm:timeautonome}}.
	\end{cases}.
\end{equation}

In the spirit of~\cite{Pages2024}, the following proposition shows that having a constant variance is equivalent to requiring $\mathbb{E}[\sigma^2(X_t)]$ to remain constant. As the proof closely resembles that of~\cite{Pages2024}, it is left to the reader.
\begin{Proposition}[Equivalence]\label{prop:equiv}
	Let \( \lambda > 0 \), and let the state-dependent diffusion coefficient \( \sigma : \mathbb{R} \to \mathbb{R} \) be a Lipschitz continuous function. Let \( X_0 \in L^2(\Omega, \mathcal{F}, \P) \) be such that \( \mathbb{E}[X_0] = m_0 \) and \( \text{Var}(X_0) = v_0 \geq 0 \). Set \( \bar{\sigma}_0^2 = \mathbb{E}[\sigma^2(X_0)] > 0 \) and set \( c = \frac{v_0}{\bar{\sigma}_0^2}\in \mathbb{R}_+ \). Assume the kernel \( K \) satisfies ~\eqref{eq:contKtilde}, ~\eqref{eq:Kcont}, and that Assumptions~\ref{ass:resolvent}-\eqref{ass:MeanRevert} and Equation~\eqref{eq:MasterEq} are in force.
	Then, the unique strong solution starting from \( X_0 \) of the scaled stochastic Volterra equation ~\eqref{eq:Volterrameanrevert}-\eqref{eq:ConstMean}, where \(\sigma(t,x) = \varsigma(t)\,\sigma(x)\) with
	\( \varsigma \) given by~\eqref{eq:VolterraStabilizerCont},
	has constant mean and satisfies the following equivalence :
	
	\begin{itemize}
		\item[(i)] \( \forall t \geq 0, \, \text{Var}(X_t) = \text{Var}(X_0) = v_0 \),
		\item[(ii)] \( \forall t \geq 0, \, \mathbb{E}[\sigma^2(X_t)] = \mathbb{E}[\sigma^2(X_0)] = \bar{\sigma}_0^2 \).
	\end{itemize}
\end{Proposition}
\noindent{\bf Remark :}	As a consequence, in~\eqref{def:stationarity_order}, Definition~(2) is equivalent to Definition~(1) when \(p=2\). The proof of this Proposition is deferred to Appendix \ref{app:lemmata}.

\bigskip
In this section we specify a family of scaled models where $b(t , x) = \mu(t)-\lambda \, x$ and the state diffusion coefficient $\sigma$ (to be specified later) satisfies the usual conditions (Lipschitz continuous) and is  sufficiently regular or smooth, specifically,  $\sigma \in C^3(\mathbb{R})$. 

\begin{Proposition}[Fake stationary regimes  (types~I and II) and asymptotics]\label{prop:GeneralSigma} Assume that assumptions~\ref{ass:resolvent}, ~\ref{ass:MeanRevert} and~\ref{ass:on_stabilizer} are in force.
	Let \(X = (X_t)_{t \geq 0}\) be a one-dimensional solution of the stabilized Volterra equation~\eqref{eq:ConstMeanVar} or ~\eqref{eq:ConstMean} starting from any random variable $X_0$ defined on
	$(\Omega, \mathcal{F}, \mathbb{P})$, with $\lambda > 0$, $\mu_\infty \in \mathbb{R}$, 
	and $\varsigma = \varsigma_{\lambda,c}$ 
	 be the unique continuous solution to Equation~\eqref{eq:VolterraStabilizer}
	for some $c \in (0, \frac{1}{[ \sigma ]^2_{\text{Lip}}}) $ (so that condition \textit{($E_{\lambda, c}$)} is satisfied) with the convention~\eqref{eq:VolterraStabilizerCont}.
	If  \( X_0 \in L^2(\mathbb{P}) \) is such that \(\mathbb{E}[X_0] = m_0, \; \text{ and} \;
	\mathrm{Var}(X_0) = v_0\) to be specified, then the following holds.
	\begin{enumerate}
		\item {\bf Case $\sigma(x) = \sigma$ is constant.}
		The solution $(X_t)_{t \geq 0}$ has a constant mean $m_0$ and variance $v_0$.
		\begin{itemize}
			\item The process exhibits a fake stationary regime of type I i.e.
			\[
			\forall t \geq 0, \quad \mathbb{E}[X_t] = m_0, \quad \text{Var}(X_t) = v_0 =c\sigma^2.\vspace{-.1cm}\]
			\item Furthermore,  if $X_0 \sim \nu^* :=  \mathcal{N}\left( m_0, v_0\right)$, this represents a fake stationary regime of type II, since in this case,  $X_t \sim X_0$ for all $t \geq 0$. ($(X_t)_{t \geq 0}$ is a Gaussian process with a fake stationary regime of type II. anyway.). $\nu^*$ is the 1-marginal distribution.
		\end{itemize}
		
		\item  {\bf Case where $\sigma$ is not constant.} Assume that the mean, variance process \((v_t := \mathbb{E}[\left(X_t- m_0\right)^2])_{t \geq 0}\) and expected squared diffusion process \((\bar\sigma^2(t):=\mathbb{E}[\sigma^2(X_t)])_{t \geq 0}\) are constant, i.e.
		
		\centerline{$
		\forall t \geq 0, \quad
		\mathbb{E}[X_t] = m_0, \quad
		\mathrm{Var}(X_t) = v_0 , \quad
		\mathbb{E}[\sigma^2(X_t)] = \bar{\sigma}_0^2.
		$}
		
		Then, a necessary and sufficient condition for this fake stationary regime of type I to hold is that there exists a function \(f\) not depending on \(t\) such that:
		{\small
		\begin{equation}\label{eq:func}
			\forall u \in [0,1] \, , \quad \forall t \geq 0 \quad \mathbb{E}\left[(X_t - m_0)^3 \partial_x^3\sigma^2(m_0 +(X_t - m_0)u)\right] = f(u).
			\vspace{-.1cm}
		\end{equation}
		}
		In this case,
		the solution \( (X_t)_{t \geq 0} \) to the Volterra equation ~\eqref{eq:Volterrameanrevert}  starting from  $X_0$ has a fake stationary regime of type I in the sense i.e. for all \( t \geq 0 \),
		{\small
		\begin{equation}\label{eq:funcFake}
			\mathbb{E}[X_t] = m_0, \quad \text{Var}(X_t) = v_0 = \frac{c(\sigma^2(m_0)+r_\infty)}{1 - c\kappa}, \quad \text{and} \quad \mathbb{E}[\sigma^2(X_t)] = \bar{\sigma}_0^2= \frac{(\sigma^2(m_0)+r_\infty)}{1 - c\kappa}.
			\vspace{-.1cm}
		\end{equation}
	    }
		\noindent \centerline{$\text{where} \; \kappa := \frac{1}{2}\partial_x^2\sigma^2(m_0)\;\text{is the curvature of } \sigma^2 \; \text{and}  \; r_\infty := \int_0^{1} \frac{(1-u)^2}{2} f(u) \;du \; \text{provided } \; \kappa c \neq 1.$}
	
	\end{enumerate}
	Moreover if \( m_0=(1-a)\frac{\mu_\infty}{\lambda}\) i.e. \(\ell_\infty=0\), as a consequence of the confluence properties in Proposition\ref{prop:confluenceLp}, for any starting value \( X_0 \in L^2(\P) \),
	
    \centerline{$
	\mathbb{E}[X_t] \to m_0, \quad \text{and} \quad \text{Var}(X_t) \to v_0 \quad \text{as} \quad t \to +\infty.
	$}	
\end{Proposition}
\noindent\textbf{Remark:}[Examples where ~\eqref{eq:func} holds.] 
1. If \( K \equiv 1 \), i.e. the solution \( (X_t)_{t \ge 0} \) to~\eqref{eq:Volterrameanrevert} 
is a (Markov) SDE, and if it admits an invariant distribution (see e.g. \cite{Pages2023}) 
\( \nu_\sigma(dx) = \pi_\sigma(x)\, \lambda_1(dx) \), then starting from 
\( X_0 \overset{d}{=} \nu_{\sigma} \) yields a fake stationary regime of type~II and, 
in particular, of type~I. In this case, for all \( t \ge 0 \), 
equation~\eqref{eq:func} corresponds to the expectation under the invariant 
distribution, i.e. \( \mathbb{E}_{\pi_{\sigma}}[\cdot] \), and thus the function  
\( f \) does not depend on \( t \).

\noindent 2. If \( K \neq 1 \) and the state-dependent diffusion coefficient is such that \(\sigma^2(x) \in \mathrm{Pol}_2(\mathbb{R})\) and is positive (see Example~\ref{Ex:Polynomial} and Proposition~\ref{prop:QuadraticSigma}).

\medskip
\noindent {\bf Proof.}
	Assume there exists at least a weak solution on the whole non-negative real line of the Voltera SDE with volatility term $\sigma(t,x)= \varsigma_{\lambda,c}(t)\sigma(x)$ starting from any $X_0 \in L^2(\mathbb{P})$ such that $\mathbb{E}[X_0] = m_0$ and $Var[X_0] = v_0.$
	The first claim $1.$ is obvious  once noted that $(X_t)_{t\ge 0}$ is a Gaussian process (and $[\sigma]_{\rm Lip}=0$).
	The last claim is a straightforward  consequence of the confluence property in Proposition~\ref{prop:confluenceLp}.
	We know that: $ 
		\E\, X_t = \E\, X_0 - \frac{1}{\lambda m_0}\Big( \E\,X_0 - m_0 \Big) \int_0^t f_{\lambda}(t-s) \mu(s)ds.$

	\noindent  {\sc Step~1.} {\em (Conditions for Fake stationary Regime of type I.)} 
	Using the second-Order Taylor Expansion of $\sigma^2$ around $m_0$ with Integral Remainder, we have:
	{\small
	\begin{equation}\label{eq:Taylor}
		\sigma^2(X_t) = \sigma^2(m_0) + \partial_x\sigma^2(m_0)Y_t + \frac{\partial_x^2\sigma^2(m_0)}{2}Y_t^2 + \int_0^{1} \frac{(1-u)^2}{2} \, Y_t^3 \, \partial_x^3\sigma^2\left(m_0 + u Y_t\right) \, du.
	\end{equation}
	}
	where $Y_t := X_t - m_0$ for $\sigma^2 \in C^3(\mathbb{R})$, and the change of variable \(u \to u - m_0\) in the integral term.
	
	Now, taking the expectation and invoking the standard Fubini lemma, we obtain:
	 {\small
	\begin{equation}\label{eq:expectedsigma}
		\bar{\sigma}(t)^2 := \mathbb{E}[\sigma^2(X_t)] = \sigma^2\left(m_0\right) + \kappa\, \mathrm{Var}(X_t) + r_t, \;
		\text{with } \; r_t = \int_0^{1} \frac{(1-u)^2}{2} \, \mathbb{E}\left[Y_t^3 \, \partial_x^3\sigma^2\left(m_0 + u Y_t\right)\right] \, du.
	\end{equation}
    }
	By the equivalence property (Proposition~\ref{prop:equiv}), the fake stationarity regime of type I holds whenever \( \bar{\sigma}(t) \) is constant in which case \({\rm Var}(X_t)\) remains constant as well (see Proposition \ref{prop:equiv}).
	It is thus necessary and sufficient that \(r_t \)  be constant (hence denoted \(r_\infty\)) or equivalently, a necessary and sufficient condition is that \(\mathbb{E}\left[Y_t^3 \partial_x^3\sigma^2(m_0+Y_tu)\right] \) is independent of \( t \) for any fixed \(u \in [0,1]\) i.e. Equation~\eqref{eq:funcFake} holds.

	\smallskip
	\noindent  {\sc Step~2.} {\em (Fake stationary Regime of type I.} 
	If \(X_0 \in L^1(\P)\) is such that \eqref{eq:CondMean} holds,
	then from equation~\eqref{eq:ConstMean} we have constant mean for every $t\ge 0$ i.e. $\E\, X_t = m_0$ 
	
	Assume that the  condition $c\kappa \neq 1$ is satisfied and as for the variance we have, from equation~\eqref{eq:var}
	\[
	\forall\, t\ge 0, \quad  {\rm Var}(X_t) ={\rm Var}(X_0) \big(1 - \frac{(f_{\lambda} * \mu)_t}{\lambda m_0}\big)^2 + \frac{1}{\lambda^2}f^2_{\lambda}* \big( \varsigma^2\,\E\, \sigma^2(X_{\cdot})\big)_t\]
	Which become, with equation~\eqref{eq:expectedsigma} in mind:
	{\small 
		\begin{equation}\label{eq:var2_2}
			{\rm Var}(X_t) = {\rm Var}(X_0) \big(1 - \frac{(f_{\lambda} * \mu)_t}{\lambda m_0}\big)^2 + \frac { 1 }{\lambda^2 }(f_{\lambda}^2 *(\varsigma\bar \sigma)^2)_t = {\rm Var}(X_0) \big(1 - \frac{(f_{\lambda} * \mu)_t}{\lambda m_0}\big)^2 + \frac { 1 }{\lambda^2 }(f_{\lambda}^2 *(\varsigma^2(\sigma^2(m_0) + \kappa {\rm Var}(X_\cdot)+r_t )))_t.
		\end{equation}
	}
Now, assuming constant variance, \(\mathrm{Var}(X_t) = v_0 \quad \text{for every } t \geq 0,\) equation~\eqref{eq:func} holds and equation~\ref{eq:expectedsigma} becomes:

\centerline{$\bar \sigma(t) ^2 := \mathbb{E}[\sigma^2(X_t)] =  \sigma^2(m_0) + \kappa {\rm Var}(X_t) + r_\infty =  \sigma^2(m_0) + \kappa v_0 + r_\infty =: \bar \sigma_0 ^2$.\\}
	
And then, the equation ~\eqref{eq:var2_2} above becomes ( where in the second line, we take advantage of the
	identity~\eqref{eq:MasterEq} satisfied by \(\varsigma \)),
	\begin{align}\label{eq:var2_2}
		\forall\, t\ge 0, \quad  v_0 = {\rm Var}(X_t)
		&= {\rm Var}(X_0) \big(1 - \frac{(f_{\lambda} * \mu)_t}{\lambda m_0}\big)^2 + \frac { \bar\sigma_0 ^2 }{\lambda^2 }(f_{\lambda}^2 *\varsigma^2)(t)\\
		&=v_0\big(1 - \frac{(f_{\lambda} * \mu)_t}{\lambda m_0}\big)^2 + (\sigma^2(m_0) + \kappa v_0 + r_\infty ) c(1-\big(1 - \frac{(f_{\lambda} * \mu)_t}{\lambda m_0}\big)^2)
	\end{align}
	Which also reads: \(v_0(1-\big(1 - \frac{(f_{\lambda} * \mu)_t}{\lambda m_0}\big)^2) = c(\sigma^2(m_0) + \kappa v_0 + r_\infty)(1-\big(1 - \frac{(f_{\lambda} * \mu)_t}{\lambda m_0}\big)^2), \; t\ge 0 \)
	i.e. the variance becomes \(v_0 = \frac{c\left(\sigma^2(m_0)+ r_\infty\right)}{1-c\kappa}>0 \)
	which is clearly solution to the equation. 
	
	Conversely one checks that this constant value for the variance solves the above equation. Let us prove that it is the only one. 
	Assume that there exist two solutions to Equation~\eqref{eq:var2_2} starting from a unique initial value ${\rm Var}(X_0) = v_0$, and let $x \in {\cal C}(\R_+, \R)$ represent the discrepancy over time between those solutions. By the linearity of Equation~\eqref{eq:var2_2}, it suffices to show  that  the equation in $x\!\in {\cal C}(\R_+,\R)$
	\[
	x(t) =  \frac{\kappa}{\lambda^2} \big(f^2_{\lambda}*(\varsigma^2.\, x)\big)_t ,\quad x(0)=0
	\]
	only has  the null function as solution. If $x$ solves the above equation, then 
	\[
	\forall\,t>0,\; |x(t)|\le  \frac{\kappa}{\lambda^2} (f^2_{\lambda}*\varsigma^2)_t  \sup_{0\le s \le t} |x(s)| \, = \kappa c |1-\big(1 - \frac{(f_{\lambda} * \mu)_t}{\lambda m_0}\big)^2|  \sup_{0\le s \le t} |x(s)| \, < c\kappa    \sup_{0\le s \le t} |x(s)| .
	\]
	where the last inequality comes from condition~\eqref{eq:MasterEq}.
	If $x\equiv \!\!\!\!\! / \,0$, there exist $\varepsilon>0$ such that $\tau_{\varepsilon} = \inf\{t: |x(t)|>\varepsilon \} <+\infty$. By continuity of $x$ it is clear that $\tau_{\varepsilon}  >0$ and $ |x(\tau_{\varepsilon})|= \sup_{0:\le s \le t} |x(s)| = \varepsilon$ which is impossible since $\kappa c \neq 1$. Consequently $x\equiv 0$. 
	We also have:
	
	\centerline{$\bar \sigma_0 ^2 = \sigma^2(m_0) + \kappa v_0 + r_\infty = \sigma^2(m_0) + \kappa \frac{c\left(\sigma^2(m_0)+ r_\infty\right)}{1-c\kappa} + r_\infty = \frac{\sigma^2(m_0)+ r_\infty}{1-c\kappa}.$}
	
	Hence $(X_t)_{t\ge 0}$ is  a fake stationary regime of type~I with the above mean and variance.  \hfill $\square$

\begin{Example}[Polynomial of degree 2]\label{Ex:Polynomial}
	Consider as in \cite{Pages2024} a square-root trinomial diffusion coefficient:
	\begin{equation}\label{eq:sigmafakeI&II}
		\sigma(x) = \sqrt{ \kappa_0 +\kappa_1\,x +\kappa_2\,x^2}\quad \mbox{ with }\quad  \kappa_i\ge 0,\;i=0,2, \;\kappa^2_1 \le 4\kappa_2\kappa_0. 
	\end{equation}
	
	\noindent 
	$\bullet$ This type of vol-vol term also appears in the quadratic rough volatility dynamic introduced in~\cite{GaJuRo2020} ($V_t = \sigma^2(X_t)$). In that model, the asset and its volatility are driven by {\em the same Brownian motion}, aiming to jointly calibrate the the S\&P 500 and VIX  smile, accounting for the so-called Zumbach effect, which links the evolution of the asset (here, an index) with its volatility.
\end{Example}
\begin{Proposition}\label{prop:QuadraticSigma}
	Under the same assumptions as Proposition \ref{prop:GeneralSigma}, We have the following claims:
	\begin{enumerate}
		\item  If the diffusion coefficient  $ \sigma$ is degenerated in the sense that $\sigma(m_0) =0$, (in particular \( \bar{\sigma}_0^2 =0\) and $v_0=0$) then the solution $X_t = m_0$ $\mathbb{P}$-$a.s.$  represents a fake stationary regime (of type I).
		\item If \(\sigma^2\) is 
		not degenerated given by \eqref{eq:sigmafakeI&II} i.e. \(\sigma^2(x) \in \mathrm{Pol}_2(\mathbb{R})\):
		
		\noindent $(a)$ Case $\kappa_2 > 0$. The solution \( (X_t)_{t \geq 0} \) to the Volterra equation~\eqref{eq:ConstMeanVar} or ~\eqref{eq:ConstMean} has a fake stationary regime of type I, in the sense that
		\[
		\forall\, t\ge 0, \quad 
		\mathbb{E}[X_t] = m_0, \; \text{Var}(X_t) = v_0= \frac{c\sigma^2(m_0)}{1 - c\kappa_2},\; \text{and}\;\; \mathbb{E}[\sigma^2(X_t)] = \bar{\sigma}_0^2 =\frac{\sigma^2(m_0)}{1 - c\kappa_2}.\]
		Moreover if $m_0= (1-a) \frac{\mu_\infty}{\lambda}$ i.e. $\ell_\infty=0$, for any starting value \( X_0 \in L^2(\P) \),
		
		\centerline{$
			\mathbb{E}[X_t] \to m_0, \quad \text{and} \quad \text{Var}(X_t) \to \frac{c\sigma^2(m_0)}{1 - c\kappa_2} \quad \text{as} \quad t \to +\infty.
			$}
	\noindent $(a)$ Case $\kappa_2 = 0$. Then $\kappa_1=0$ and \((X_t)_{t\geq0}\) is a Gaussian process with a fake stationary regime of type~II, with \( \nu \) as one-dimensional marginal distribution, provided that if $X_0 \sim \nu := \mathcal{N}\!\left(m_0,c\kappa_0^2\right)$.
	\end{enumerate}
\end{Proposition}

\noindent {\bf Proof.} {\em (Applicability of equation~\eqref{eq:func}).} 
	\begin{enumerate}
		\item First, in the degenerate setting  $\sigma(m_0) =0$, one has $\bar{\sigma}_0^2 = \mathbb{E}[\sigma^2(X_0)]=0$, we get $v_0=0$ owing to equation~\eqref{eq:VolterraStabilizer} since $\lim_{t\to+\infty} \big(1 - \frac{(f_{\lambda} * \mu)_t}{\lambda m_0}\big) = \ell_\infty < 1 \Rightarrow \big(1 - \frac{(f_{\lambda} * \mu)_t}{\lambda m_0}\big) \neq 1 $ (at least for $t$ large enough). As a consequence, ${\rm Var}(X_t)=0$ for every $t\ge 0$. But, we know that $\E\, X_t = \E\, X_0=  m_0$, it follows that  $X_t= m_0$ $\P$-$a.s.$ and \(\forall t \geq 0,  \quad \mathcal{L}(Y_t)(dy) = \delta_{0}(dy) \) so that 
		\centerline{$  \forall t \geq 0,  \quad 
		\mathbb{E}\left[Y_t^3 \, \partial_x^3 \sigma^2\left(m_0 + u Y_t\right)\right] 
		= \int_{\mathbb{R}} y^3 \, \partial_x^3 \sigma^2\left(m_0 + u y\right) \, \mathcal{L}(Y_t)(dy) = 0 \; \text{and}\; r_\infty =0.
		$}
		\item Secondly, if \( \partial_x^3 \sigma^2(v) = 0,\; \forall v \in \mathbb{R} \), then the integral reminder in~\eqref{eq:Taylor} necessarily vanishes.  
		This corresponds to the \emph{trinomial setting}, which has already been extensively studied in \cite{Pages2024} and in which case if $ \kappa_2>0$,  $[\sigma]_{\rm Lip} = \sqrt{\kappa_2}$, the curvature \(\kappa = \kappa_2\) and \(r_\infty = 0\) (since \(f\equiv 0\) in~\eqref{eq:func}).	\hfill$\Box$	
	\end{enumerate}
\medskip
\noindent {\bf Practitioner's corner:} 1. The constraint \( c \in \left( 0, \frac{1}{\kappa} \right) \) implies that we treat \( c \) as a free parameter, from which we can deduce \( v_0 \) and \( \bar{\sigma}_0^2 \). 2. The presence of the {\em stabilizer} $\varsigma_{\lambda,c}$ allows a better control of the behaviour of the  equation since  it induces an $L^2$-confluence and a stability of first two moments if needed.
		3. Note that \([\sigma]^2_{\rm Lip} = \kappa_2=\kappa\) so that, in practice, if we rather fix the value of $v_0$, then $c = \frac{v_0}{\sigma^2(m_0)+ v_0\kappa}$ so that, $\sigma$ being $\sqrt{\kappa}$-Lipschitz continuous, one has $ c \kappa = \frac{v_0\kappa}{\sigma^2(m_0)+ v_0\kappa} <1$ provided $\sigma^2(m_0)>0$ which ensures the $L^2$-confluence of the paths of the solution (Proposition \ref{prop:confluenceLp} further on). 

\section{Towards Long run behaviour: asymptotics and confluence}\label{sect-LongRunB}
\begin{Remark}\label{rem:on_sigma}
	Let $ m_0 \in \R$, by the Lipschitz assumption on the state-dependent \(\sigma\), one has for every \(x \in \mathbb{R}\),
	
	\centerline{$
	\sigma^2(x) \leq \Big(\sigma(m_0) + [\sigma]_{\text{Lip}} |x - m_0|\Big)^2 \leq \kappa_0 + \kappa_2 |x - m_0|^2
	$}
	where \(k_0 = k_0 (\epsilon):=(1+\epsilon) |\sigma(m_0)|^2\) and \(k_2= k_2(\epsilon):=(1+\frac{1}{\epsilon})[\sigma]^2_{\text{Lip}}\), owing to Young's inequality
	$(a+b)^2 \le (1+\epsilon) a^2+(1+\frac1\epsilon )b^2$. Therefore, we can always assume that {\em $\sigma$ is sublinear} i.e.:
	\begin{equation}\label{eq:on_sigma}
		(SL_{\sigma}) : \exists k_0 = k_0 (\epsilon) \in \R_+, k_2= k_2(\epsilon) \in \R_+ \quad \textit{such that} \quad \forall x \in \mathbb{R}, |\sigma(x)|^2 \leq \kappa_0 + \kappa_2 |x -m_0|^2.  
	\end{equation}
\end{Remark}
\subsection{Moments control.}\label{subsect-Moments}
\begin{Lemma}[Best constant in a BDG inequality (see Remark 2 in \cite{carlen1991lp})]\label{lm:best_bdg_constant}
	
	Let $M$ be a continuous local martingale null at $t = 0$. Then, for every $p \geq 1$,
	
	\centerline{$
	\frac{1}{\sqrt{p}} \| M_t \|_p \leq 2 \| \langle M \rangle_t^{1/2} \|_p.
	$}
\end{Lemma}
\begin{Proposition}[Moment control] \label{prop:Momentctrl} Assume ~\eqref{assum:convol} (ii)  and \(({\cal R}_\lambda)\) $(ii)$ hold. Let $\sigma(t,x):= \varsigma(t)\sigma(x)$  where $\varsigma $ is a non-negative and continuous solution to~\eqref{eq:VolterraStabilizer} for some fixed $\lambda,c >0$ (i.e.($E_{\lambda}$) is in force).
	Let $(X_t)_{t\ge 0}$ be  the solutions to the Volterra equation~\eqref{eq:ConstMean}--\eqref{eq:MasterEq} starting  from any random variable $X_0$.\\
	
	$(a)$ {\em First two moments}. Assume $X_0\!\in L^2(\P)$ and $c\!\in \big(0,\frac{1}{[\sigma]^2_{\text{Lip}}}\big)$. Then, one has:
	\[\Big| \E\, \big(X_t\big)-m_0\Big| \le \Big|1 - (f_{\lambda} * \frac{\mu}{\lambda m_0})_t\Big| \Big| \E\,\big(X_0\big)-m_0 \Big|,\;\; t\geq 0\]

	\[
	\sup_{t\ge 0}\Big\||X_t-m_0|\Big\|_2 \le   \left[\frac{\sqrt{c}}{1-[\sigma]_{\text{Lip}}\sqrt{c}}|\sigma(m_0)|\right]\vee   \Big\||X_0-m_0|\Big\|_2< +\infty.
	\]
    
	$(b)$ {\em $L^p$-moments}. Let $p\!\in (2,+\infty)$. If $X_0\!\in L^p(\P)$ and $c$ is such that  $\rho_p  := 4c\,p\,[\sigma]^2_{\text{Lip}}<1$, then 
	\[
	\sup_{t\ge 0}\Big\|  |X_t-m_0|\Big\|_p\le \inf_{\epsilon\in (0, \frac{1}{\rho_p}-1)}   \left[\frac{2\sqrt{pc(1+\epsilon)}}{1 -2[\sigma]_{\text{Lip}}\sqrt{pc(1+\epsilon)}} 
	\big|\sigma(m_0)\big|\right]\vee \left[ (1+1/\epsilon)^{\frac12} \Big\||X_0-m_0|\Big\|_p\right] < +\infty
	\]
\end{Proposition} 

The proof is postponed to Appendix~\ref{app:lemmata}. It relies on techniques similar to those in~\cite{Pages2024}, which extend the methods developed for the Markovian framework, as discussed extensively in~\cite[Chapter 7]{pages2018numerical}.

	\subsection{$L^p$-Confluence or Contraction Properties}\label{subsect-confluence}
	Fix \( p>0\). Let \((X_t)_{t\ge0}\) and \((X'_t)_{t\ge0}\) be two solutions of the Volterra stochastic equation~\eqref{eq:Volterrameanrevert} with initial conditions \(X_0,X'_0\in L^p(\mathbb P)\). According to the Lipschitz assumption on the state-dependent \(\sigma\),  $\exists \; \kappa > 0$ such that for every \(t\geq0\), 
	\begin{equation}\label{eq:lipschitz} 
		\Big\| |\sigma(X_{s}) - \sigma(X^\prime_{s})|^2 \Big\|_{\frac{p}{2}} 
		\leq \kappa \Big\| | X_{s} - X^\prime_{s} | \Big\|_p^2.
	\end{equation}
	The following $L^p$-confluence result extends the $L^2$-confluence established in~\cite[Proposition 4.3]{Pages2024} to arbitrary $p>0$. Unlike the approach in that reference, the proof here does not rely on similar arguments, but instead leverages a Banach fixed point arguments. 
	\begin{Proposition}[$L^p$-confluence]\label{prop:confluenceLp} 
		Assume assumptions ~\ref{assum:convol} (ii),~\ref{ass:resolvent} and~\ref{ass:MeanRevert} are in force.  Assume in~\eqref{eq:Volterrameanrevert}-\eqref{eq:Volterrameanrevert_}, that the diffusion coefficient $\sigma(t,x):= \varsigma(t)\sigma(x)$  where $\varsigma $ is a non-negative, continuous function given by Equation~\eqref{eq:VolterraStabilizerCont} for some fixed $\lambda,c >0$ and \( \sigma : \mathbb{R} \to \mathbb{R} \) is a Lipschitz continuous function.
		Let $p>0$, for $X_0,X'_0\!\in L^p({\mathbb P})$, we consider the solutions to Volterra equation~\eqref{eq:ConstMean}--\eqref{eq:MasterEq} denoted $(X_t)_{t\ge 0}$ and $(X'_t)_{t\ge 0}$ starting from $X_0$ and $X'_0$ respectively. Let $c\!\in \big(0,\frac{1}{\kappa}\big)$, where $\kappa$ is defined in \ref{eq:lipschitz}, set \(\rho_p := 4 p\, c\,\kappa.\) and assume that \(\rho_p < 1-\ell_\infty^2\).
		Then, one has:
		\begin{itemize}
			\item[(a)] There exists a continuous non-negative function \( \varphi_{\infty,p}^{\lambda,c, K, \mu}=: \varphi_{\infty,p}: \mathbb{R}_+ \to [0, \frac{1}{(1-\sqrt{\rho_p})^2}] \), such that \( \varphi_{\infty,p}(0) = \frac{1}{1-\rho_p} \), \(0\leq \lim_{t\to\infty}\varphi_{\infty,p}(t)
			\leq \frac{\ell_\infty^2}{\big(1-\sqrt{\rho_p}\sqrt{1-\ell_\infty^2}\big)^2} \), only depending on \( \lambda, c, \mu \), and the kernel \( K \), such that :
			
			\centerline{$
				\forall t \geq 0, \quad \E\, \Big(\Big| X_t - X^\prime_t\Big|^p\Big) \leq \varphi_{\infty,p}(t) \E\, \Big(\Big| X_0 - X^\prime_0\Big|^p\Big).
				$}
			\item[(b)] This result can be written using the \(p\)-Wasserstein distance between marginals of \(X\) and \(X^\prime\):
			\[
			\forall t \geq 0, \quad W_p([X_t],[X'_t]) \le \varphi_{\infty,p}(t)^{1/2}\,W_p([X_0],[X'_0]).
			\]
			\item[(c)]  In particular, whenever \(\ell_\infty=0\), the limit yields \(\varphi_{\infty,p}(t)\to0\) and thus the process is contracting in \(W_p\) as \(t\to\infty\) i.e more generally finite-dimensional \( W_p \)-convergence of marginals.
		\end{itemize}
	\end{Proposition}
	\noindent {\bf Proof.}
	{\em By a Banach fixed point argument on the complete space \((C_b([0,\infty),\mathbb{R}),\|\cdot\|_\infty).\)\\}
	Fix \(p>0\). Let \((X_t)_{t\ge0}\) and \((X'_t)_{t\ge0}\) be two solutions of the same SVIE~\eqref{eq:ConstMean} with initial conditions \(X_0,X'_0\in L^p(\mathbb P)\). Set $\Delta_t = X_t - X^\prime_t \in L^p (\P)$ for every $t\ge 0$. one writes owing to equation ~\ref{eq:Volterrameanrevert_}:
	\begin{align*}
		X_t-X^\prime_t &= \big(X_0 - X^\prime_0\big) \Big(1-\frac{(f_\lambda * \mu)(t)}{\lambda m_0}\Big) + \frac{1}{\lambda}\int_0^tf_{\lambda}(t-s)\varsigma(s) \Big(\sigma(X_{s}) -\sigma(X^\prime_{s})\Big)dW_s 
	\end{align*}
	Let \( \bar \delta_t = \Big\| |\Delta_t|\Big\|_p \) for convenience. Under our assumptions, \( t \mapsto \bar \delta_t \) is continuous (see e.g.~\cite[Theorem 3.3]{GnabeyeuPages2026b} ). 
	Let \(C_p^{\mathrm{BDG}}>0\) denote a BDG constant in \(L^p\), which by Lemma~\ref{lm:best_bdg_constant} is given by $C_p^{\mathrm{BDG}}:=2\sqrt{p}$. Set \(\rho_p := c\,(C_p^{\mathrm{BDG}})^2\,\kappa = 4\,p\,c\,\kappa.\)
	Owing to the triangle inequality and applying the   {\em BDG} inequality to the (a priori) local martingale $M_u = \int_0^u f_{\lambda}(t-s)\varsigma(s)\sigma(X_{s})dW_s$, $0\le s\le t$, (see~\cite[Proposition~4.3]{RevuzYor}) follow by the generalized Minkowski inequality, we get:
	\begin{align*}
		\Big\| | X_t - X^\prime_t |\Big\|_p &\le \Big\|  |X_0-X^\prime_0|\Big\|_p \Big|1-\frac{(f_\lambda * \mu)(t)}{\lambda m_0}\Big| + \frac{C_p^{BDG}}{\lambda}\Big\|\left(f_{\lambda}^2 \stackrel{dt}{*} \varsigma^2(\cdot)|\sigma(X_{s}) -\sigma(X^\prime_{s})|^2\right)_t \Big\|_{\frac p2}^{\frac 12}\\
		& \le  \Big\|  |X_0-X^\prime_0|\Big\|_p \Big|1-\frac{(f_\lambda * \mu)(t)}{\lambda m_0}\Big| + \frac{C_p^{BDG}}{\lambda}\Big( \int_0^t f^2_{\lambda}(t-s)\varsigma^2(s) \big\||\sigma(X_{s}) -\sigma(X^\prime_{s})|^2\big\|_{\frac p2}\Big)^{\frac 12}
	\end{align*}
	
	Fix \(\epsilon > 0\), using the elementary inequality $(a+b)^2 \le (1+1/\epsilon) a^2+(1+\epsilon )b^2$ for $\epsilon\!\in (0,1/\rho_p-1)$ i.e. \(\beta := \rho_p(1+\varepsilon) \;<\; 1\), it follows  owing to equation ~\eqref{eq:lipschitz} that:
	\begin{align*}
		\Big\| |X_t - X^\prime_t|\Big\|_p^2 &\le \Big\|  |X_0-X^\prime_0 |\Big\|_p^2 \Big(1-\frac{(f_\lambda * \mu)(t)}{\lambda m_0}\Big)^2(1+1/\epsilon)
		+ \frac{\rho_p}{c\lambda^2}(1+\epsilon) \int_0^t f^2_{\lambda}(t-s)\varsigma^2(s) \Big\||X_{s} -X^\prime_{s}|^2\Big\|_{p}ds
	\end{align*}
	which entails:
	\begin{equation}\label{eq:lp_confl}
		\Big\| |\Delta_t| \Big\|_p^2 \le \Big\| |\Delta_0| \Big\|_p^2 \Big(1-\frac{(f_\lambda * \mu)(t)}{\lambda m_0}\Big)^2 \left(1 + \frac{1}{\epsilon}\right) 
		+ \left(1 + \epsilon\right) \frac{\rho_p}{c\lambda^2} \int_0^t f_{\lambda}^2(t - s) \varsigma^2(s) 
		\Big\| |\Delta_{s}| \Big\|_p^2 \, ds.
	\end{equation}
	i.e. we obtain, for all \(t\ge0\),
	\begin{equation}\label{eq:key_lp}
		\bar\delta_t^2 \le \bar\delta_0^2\,\Big(1-\frac{(f_\lambda * \mu)(t)}{\lambda m_0}\Big)^2\Big(1+\frac{1}{\varepsilon}\Big)
		+ \rho_p(1+\varepsilon)\frac{1}{\lambda^2 c}\int_0^t f_\lambda^2(t-s)\varsigma^2(s)\,\bar\delta_s^2\,ds.	
	\end{equation}
	\noindent  {\sc Step~1.} {\em  Non-expansivity via a deterministic stopping-time argument:}
	For the fixed \(\varepsilon>0\) such that \(\beta:=\rho_p(1+\varepsilon)<1\), let \( \eta > 0 \) such that \(  1 + \frac{1}{\epsilon} < \rho_p(1+\epsilon)(1 + \eta)^2 \) and define the stopping time
	\[
	\tau_\eta := \inf\{t\ge0:\bar\delta_t > (1+\eta)\bar\delta_0\}
	\]
	(with the convention \(\inf\varnothing=+\infty\)). If \(\tau_\eta<\infty\), then for \(s\le\tau_\eta\) we have \(\bar\delta_s\le(1+\eta)\bar\delta_0\) and by continuity \(\bar\delta_{\tau_\eta}=(1+\eta)\bar\delta_0\). Evaluating \eqref{eq:key_lp} at \(t=\tau_\eta\) and bounding \(\bar\delta_s^2\le(1+\eta)^2\bar\delta_0^2\) in the integral combined with the identity  $f^2_{\lambda}*\varsigma^2 = c \lambda^2 (1-\big(1-\frac{(f_\lambda * \mu)(t)}{\lambda m_0}\big)^2)$ yields:
	\begin{align*}
		\bar \delta_{\tau_{\eta}}^2 &\le \bar \delta_0^2 \left[ \big(1-\frac{(f_\lambda * \mu)(t)}{\lambda m_0}\big)^2(\tau_{\eta})(1+\frac1\varepsilon) + (1 - \big(1-\frac{(f_\lambda * \mu)(t)}{\lambda m_0}\big)^2(\tau_{\eta})) \rho_p(1+\varepsilon) (1 + \eta)^2 \right] \\
		&\leq \delta_0^2 \left[ \big(1-\frac{(f_\lambda * \mu)(t)}{\lambda m_0}\big)^2(\tau_{\eta})\Big(1+\frac1\varepsilon - \rho_p(1+\varepsilon) (1 + \eta)^2\Big) + \rho_p(1+\varepsilon) (1 + \eta)^2 \right] \\
		&< \rho_p (1+\varepsilon) (1 + \eta)^2 \bar \delta_0^2 <  (1 + \eta)^2 \bar \delta_0^2.
	\end{align*}
	which leads to a contradiction. Whence  \( \tau_{\eta}  = +\infty \) i.e., \( \bar \delta_s \le (1 + \eta) \bar \delta_0 \) for all \( s \ge 0 \). This holds for every \( \eta > 0 \), implying the non-expansivity bound \( \bar \delta_t \le \bar \delta_0 \) for all \( t \ge 0 \) by letting \(\eta\downarrow0\).
	
	\medskip
	\noindent  {\sc Step~2.} {\em  Iteration and the Volterra map:}
	Substituting this (i.e. \( \bar \delta_t \le \bar \delta_0 \)) into~\eqref{eq:lp_confl} combined with the stabilizer identity  $f^2_{\lambda}*\varsigma^2 = c \lambda^2 (1-\big(1-\frac{(f_\lambda * \mu)(t)}{\lambda m_0}\big)^2)$ gives, for all \( t > 0 \),
	\[
	\bar\delta_t^2 \le \bar\delta_0^2 \,\varphi_{1,p}^{\varepsilon}(t), \quad \text{where} \quad \,\varphi_{1,p}^{\varepsilon}(t) := \Big(1+\frac{1}{\varepsilon}\Big) \big(1-\frac{(f_\lambda * \mu)(t)}{\lambda m_0}\big)^2 + \rho_p (1+\varepsilon)  (1 - \big(1-\frac{(f_\lambda * \mu)(t)}{\lambda m_0}\big)^2).
	\]
	Note that \( \varphi_{1,p}^{\varepsilon}(t) = \rho_p (1+\varepsilon) + \big(1-\frac{(f_\lambda * \mu)(t)}{\lambda m_0}\big)^2(1+\frac{1}{\varepsilon} - \rho_p (1+\varepsilon)) \) satisfies:
	\[
	\varphi_{1,p}^{\varepsilon}(0) = 1 +\frac{1}{\varepsilon}, \quad \varphi_1 \text{ is continuous, } M_1:=\|\varphi_{1,p}^{\varepsilon}\|_\infty \leq 1 + \frac{1}{\epsilon} + \rho_p (1+\varepsilon)
	\]
	Substituting this upper bound for \( \delta_t^2 \) (i.e. $\bar\delta_t^2 \le \bar\delta_0^2 \varphi_{1,p}^{\varepsilon}(t)$) into~\eqref{eq:lp_confl} yields
	\[
	\bar\delta_t^2 \le \bar \bar\delta_0^2 \varphi_{2,p}^{\varepsilon}(t), \quad \text{where} \quad \varphi_{2,p}^{\varepsilon}(t) := \Big(1+\frac{1}{\varepsilon}\Big) \big(1-\frac{(f_\lambda * \mu)(t)}{\lambda m_0}\big)^2 + \rho_p (1+\varepsilon)  \int_0^t f_{\lambda}^2(t - s) \varsigma^2(s) \varphi_{1,p}^{\varepsilon}(s) \frac{ds}{\lambda^2 c}.
	\]
	and inductively for \(k\ge2\)
	\begin{equation}\label{eq:phi_iter}
		\bar\delta_t^2 \le \bar \bar\delta_0^2 \varphi_{k,p}^{\varepsilon}(t), \quad 	\varphi_{k,p}^{\varepsilon}(t)
		:= \big(1-\frac{(f_\lambda * \mu)(t)}{\lambda m_0}\big)^2 \Big(1+\frac{1}{\varepsilon}\Big) + \rho_p (1+\varepsilon) \frac{1}{\lambda^2 c}\int_0^t f_\lambda^2(t-s)\varsigma^2(s)\,\varphi_{k-1,p}^{\varepsilon}(s)\,ds.
	\end{equation}
	To obtain a uniform sup-bound, put \(M_k:=\|\varphi_{k,p}^{\varepsilon}\|_\infty\). 
	As \(-1\leq 1-\frac{(f_\lambda * \mu)(t)}{\lambda m_0}\le1\) by condition~\eqref{eq:MasterEq},we get from \eqref{eq:phi_iter} together with the equation~\eqref{eq:MasterEq} satisfies by \(\varsigma^2\) that \(M_k \le \Big(1+\frac{1}{\varepsilon}\Big) + \beta M_{k-1}.\)
	Iterating yields (since \(\beta<1\)):
	
	\(M_k \le \Big(1+\frac{1}{\varepsilon}\Big)\sum_{j=1}^{k-1}\beta^j + \beta^{k-1} M_1
	\le \max(M_1, \frac{1+\frac1\varepsilon}{1-\beta}) \le \max(1+\frac1\varepsilon + \beta, \frac{1+\frac1\varepsilon}{1-\beta}) = \frac{1+\frac1\varepsilon}{1-\beta}.\)
	Thus for every \(k\ge1\) and \(t\ge0\) one has the uniform bound
	\begin{equation}\label{eq:uniformbound}
		0 \le \varphi_{k,p}^{\varepsilon}(t) \le \frac{1+1/\varepsilon}{1-\rho_p (1+\varepsilon)} := M_*^\epsilon.
	\end{equation}
	
	\noindent  {\sc Step~3.} {\em  Define the operator  \(\mathcal T:C_b([0,\infty))\to C_b([0,\infty))\) :}
	\begin{equation}\label{eq:opT}
		(\mathcal T\psi)(t)
		:= \big(1-\frac{(f_\lambda * \mu)(t)}{\lambda m_0}\big)^2\Big(1+\frac{1}{\varepsilon}\Big) +  \rho_p (1+\varepsilon) \frac{1}{\lambda^2 c}\int_0^t f_\lambda^2(t-s)\varsigma^2(s)\,\psi(s)\,ds.
	\end{equation}
	and, for \(k\ge2\), set \(\varphi_{k,p}^{\varepsilon}=\mathcal T\varphi_{k-1,p}^{\varepsilon}\).  
	The operator \(\mathcal T\) is linear 
	and for any \(\psi_1,\psi_2\in C_b\),
	\[
	\|\mathcal T\psi_1 - \mathcal T\psi_2\|_\infty
	\le \rho_p (1+\varepsilon)  \cdot\sup_{t\ge0}\frac{1}{\lambda^2 c}\int_0^t f_\lambda^2(t-s)\varsigma^2(s)\,ds \;\cdot\;\|\psi_1-\psi_2\|_\infty
	= \rho_p (1+\varepsilon)  \|\psi_1-\psi_2\|_\infty
	\]
	because the convolution integral equals \(1-\big(1-\frac{(f_\lambda * \mu)(t)}{\lambda m_0}\big)^2\le1\). By assumption \(\rho_p (1+\varepsilon) <1\), so \(\mathcal T\) is a strict contraction in \(\|\cdot\|_\infty\) with Lipschitz constant \(\rho_p (1+\varepsilon) <1\) on the complete or Banach space \(C_b([0,\infty))\) with the sup norm. The Banach fixed point theorem therefore provides a unique fixed point \(\varphi_{\infty,p}^{\varepsilon}\in C_b([0,\infty))\) and, moreover, \(\varphi_{k,p}^{\varepsilon} = \mathcal T^{k-1}\varphi_{1,p}^{\varepsilon} \xrightarrow[k\to\infty]{\|\cdot\|_\infty} \varphi_{\infty,p}^{\varepsilon}.\)
	
	In particular the convergence is uniform on \([0,\infty)\) i.e. \(\varphi_{k,p}^{\varepsilon}=\mathcal T^{k-1}\varphi_{1,p}^{\varepsilon}\) converges uniformly (on \([0,\infty)\)) to \(\varphi_{\infty,p}^{\varepsilon}\). For every \(t\ge0\) the \(L^p\)-norm satisfies
	\(	\bar\delta_t^2 \le \bar\delta_0^2\,\varphi_{\infty,p}^{\varepsilon}(t).\)
	
	\medskip
	\noindent  {\sc Step~4.} {\em  Limit equation and \(\varepsilon\)-dependent asymptotic bound:}
	Passing to the limit in \eqref{eq:phi_iter} yields that \(\varphi_{\infty,p}^{\varepsilon}\) satisfies the Volterra or functional fixed-point equation
	\begin{equation}\label{eq:varphi_fixed}
		\varphi_{\infty,p}^{\varepsilon}(t)
		= \big(1-\frac{(f_\lambda * \mu)(t)}{\lambda m_0}\big)^2\Big(1+\frac{1}{\varepsilon}\Big) + \rho_p (1+\varepsilon) \frac{1}{\lambda^2 c}\int_0^t f_\lambda^2(t-s)\varsigma^2(s)\,\varphi_{\infty,p}^{\varepsilon}(s)\,ds.
	\end{equation}
	By the uniform bound in equation~\eqref{eq:uniformbound}, \(\varphi_{\infty,p}^{\varepsilon}\) is bounded and nonnegative on \([0,\infty)\) i.e. \(\forall t \geq 0, \quad 0 \le \varphi_{\infty,p}^{\varepsilon}(t) \le \frac{1+1/\varepsilon}{1-\rho_p (1+\varepsilon)} .\)
	Taking \(\liminf_{t \to +\infty}, \limsup_{t\to\infty}\) in \eqref{eq:varphi_fixed} and using \(\big(1-\frac{(f_\lambda * \mu)(t)}{\lambda m_0}\big)^2\to \ell_\infty^2\) and the stabilizer identity we obtain $\underline{\ell}_{\infty,p}^{\varepsilon}, \ell_{\infty,p}^{\varepsilon} := \liminf_{t \to +\infty} \varphi_{\infty,p}^{\varepsilon}(t), \limsup_{t \to +\infty} \varphi_{\infty,p}^{\varepsilon}(t) \in [0,M_*^\epsilon]$. Now,  $\underline{\ell}_{\infty,p}^{\varepsilon}, \ell_{\infty,p}^{\varepsilon} \in [0,M_*^\epsilon]$ implies that for any $\eta > 0$, there exists $t_\eta \in \R^+$ such that for $t \geq t_\eta$, $ \underline{\ell}_{\infty,p}^{\varepsilon} - \eta \leq \varphi_{\infty,p}^{\varepsilon}(t)\leq \ell_{\infty,p}^{\varepsilon} + \eta $. Then, we obtain on the first hand,
	\begin{align*}
		\int_0^t f^2_{\lambda}(t-s) \varsigma^2(s) \varphi_{\infty,p}^{\varepsilon}(s) \frac{ds}{\lambda^2 c} &\leq \frac{1}{c\lambda^2} \int_{t_\eta}^t f^2_{\lambda}(t-s) \varsigma^2(s) (\ell_{\infty,p}^{\varepsilon} + \eta) \, ds + \frac{1}{c\lambda^2} \int_{0 }^{t_\eta} f^2_{\lambda}(t-s) \varsigma^2(s)\varphi_{\infty,p}^{\varepsilon}(s)\, ds\\
		&\leq \frac{1}{c\lambda^2} \int_{t_\eta}^t f^2_{\lambda}(t-s) \varsigma^2(s) (\ell_{\infty,p}^{\varepsilon} + \eta) \, ds + \frac{1}{c\lambda^2} M_*^\epsilon \int_{t - t_\eta}^t f^2_{\lambda}(u) \varsigma^2(t-s) \, du.
	\end{align*}
	where the second term on the right-hand side of the last inequality follows from the fact that \( \varphi_{\infty,p}^{\varepsilon}(t - u) \leq M_*^\epsilon\) for all \(u \leq t \leq t_\eta\) and vanishes as t goes to infinity.
	
	Since $f_\lambda \in \mathcal{L}^2 (\R_+)$ and $\lim_{t \to +\infty} \big(1-\frac{(f_\lambda * \mu)(t)}{\lambda m_0}\big)^2 = \ell_\infty^2$ both owing to Assumption \ref{ass:resolvent}, we conclude from equation~\eqref{eq:varphi_fixed} and the identity satisfied by $\varsigma$:
	\[
	\ell_{\infty,p}^{\varepsilon} =:\limsup_{t \to +\infty} \varphi_{\infty,p}^{\varepsilon}(t)  \leq (1+\frac1\varepsilon)\,\ell_\infty^2 + \rho_p (1+\varepsilon)(\ell_{\infty,p}^{\varepsilon} + \eta)(1-\ell_\infty^2) \overset{\eta \to 0}{\implies} \ell_{\infty,p}^{\varepsilon} \leq \frac{(1+\frac1\varepsilon)\,\ell_\infty^2}{1-\rho_p (1+\varepsilon)(1-\ell_\infty^2)}
	\]	
	On the other hand, we also have:
	\begin{align*}
		\int_0^t f^2_{\lambda}(t-s) \varsigma^2(s) \varphi_{\infty,p}^{\varepsilon}(s) \frac{ds}{\lambda^2 c} &\geq \frac{1}{c\lambda^2} \int_{t_\eta}^t f^2_{\lambda}(t-s) \varsigma^2(s) (\underline{\ell}_{\infty,p}^{\varepsilon} - \eta) \, ds + \int_{t - t_\eta}^t f^2_{\lambda}(u) \, \varsigma^2(t-u)\varphi_{\infty,p}^{\varepsilon}(t-u)\, \frac{du}{c\lambda^2} \\
		&\geq \frac{1}{c\lambda^2} \int_{t_\eta}^t f^2_{\lambda}(t-s) \varsigma^2(s) (\underline{\ell}_{\infty,p}^{\varepsilon} - \eta) \, ds.
	\end{align*}
	Therefore, still with the fact that $f_\lambda \in \mathcal{L}^2 (\R_+)$ and $\lim_{t \to +\infty} \big(1-\frac{(f_\lambda * \mu)(t)}{\lambda m_0}\big)^2 = \ell_\infty^2$, we obtain from equation~\eqref{eq:varphi_fixed} and the identity satisfied by $\varsigma$:
	\[
	\underline{\ell}_{\infty,p}^{\varepsilon} =:\liminf_{t \to +\infty} \varphi_{\infty,p}^{\varepsilon}(t)  \geq (1+\frac1\varepsilon)\,\ell_\infty^2 + \rho_p (1+\varepsilon)(\underline{\ell}_{\infty,p}^{\varepsilon} - \eta)(1-\ell_\infty^2) \overset{\eta \to 0}{\implies} \underline{\ell}_{\infty,p}^{\varepsilon} \geq \frac{(1+\frac1\varepsilon)\,\ell_\infty^2}{1-\rho_p (1+\varepsilon)(1-\ell_\infty^2)}.
	\]
	Consequently, \(\underline{\ell}_{\infty,p}^{\varepsilon} = \ell_{\infty,p}^{\varepsilon} = \frac{(1+\frac1\varepsilon)\,\ell_\infty^2}{1-\rho_p (1+\varepsilon)(1-\ell_\infty^2)} := L(\epsilon)\).
	The minimizer of \(L(\epsilon)\) in \((0,1/\rho_p-1)\) is \(	\varepsilon_*=\frac{1}{\sqrt{\rho_p(1-\ell_\infty^2)}}-1,\)
	which is admissible iff \(\rho_p < 1-\ell_\infty^2\). In that admissible case one obtains the optimal asymptotic value: \(\inf_{\varepsilon\in(0,1/\rho_p-1)} L(\epsilon) = \ell_{\infty,p}^{\varepsilon_*}
	= \frac{\ell_\infty^2}{\big(1-\sqrt{\rho_p}\sqrt{1-\ell_\infty^2}\big)^2}.\)
	
	\medskip
	\noindent  {\sc Step~5.} {\em Passage to the \( \varepsilon \)-Free Control:} Finally, optimizing \(\varphi_{\infty,p}^{\varepsilon}\) over admissible \(\varepsilon\) gives the \(\varepsilon\)-free control i.e. passing to the infimum over admissible \(\varepsilon\) gives the claimed \(\varepsilon\)-free control with \(\varphi_{\infty,p}(t)\).
	\(	\bar\delta_t^2 \le \bar\delta_0^2\,\varphi_{\infty,p}(t),\qquad
	\varphi_{\infty,p}(t):=\inf_{\varepsilon\in(0,1/\rho_p-1)}\varphi_{\infty,p}^{\varepsilon}(t). \quad \)
	Now, note that \[ \forall t \geq 0, \quad 0 \le \varphi_{\infty,p}(t):=\inf_{\varepsilon\in(0,1/\rho_p-1)}\varphi_{\infty,p}^{\varepsilon}(t) \le \inf_{\varepsilon\in(0,1/\rho_p-1)}\frac{1+1/\varepsilon}{1-\rho_p (1+\varepsilon)} = M_*^{\epsilon = \frac{1}{\sqrt{\rho_p}}-1}= \frac{1}{(1-\sqrt{\rho_p})^2}.\]
	Moreover, \(\limsup_{t\to\infty}\varphi_{\infty,p}(t) = \limsup_{t\to\infty} \inf_{\varepsilon\in(0,1/\rho_p-1)} \varphi_{\infty,p}^{\varepsilon}(t) \leq \limsup_{t\to\infty} \varphi_{\infty,p}^{\varepsilon_*}(t)  =\ell_{\infty,p}^{\varepsilon_*}\) 
	so that 
	\(\lim_{t\to\infty}\varphi_{\infty,p}(t)
	\leq \frac{\ell_\infty^2}{\big(1-\sqrt{\rho_p}\sqrt{1-\ell_\infty^2}\big)^2}, \)
	with equality if the  uniform convergence  \(\sup_{\varepsilon\in(0,1/\rho_p-1)} \big|\varphi_{\infty,p}^{\varepsilon}(t)-\ell(\varepsilon)\big|\xrightarrow[t\to\infty]{}0\;\) holds.
	Hence \(\|X_t-X'_t\|_p \le \varphi_{\infty,p}(t)^{1/2}\,\|X_0-X'_0\|_p\) for every \(t\ge0\), 
	and therefore, by coupling, for the \(p\)-Wasserstein distance between marginals,\(	W_p([X_t],[X'_t]) \le \varphi_{\infty,p}(t)^{1/2}\,W_p([X_0],[X'_0]).\)
	In particular, if \(\ell_\infty=0\) the asymptotic bound above yields \(\varphi_{\infty,p}(t)\to0\) and thus the process is contracting in \(W_p\) as \(t\to\infty\).
	\noindent This completes the proof. \hfill $\square$
	
	\medskip
	\noindent {\bf Remarks.} 
	1. The function \( \varphi_{\infty,p} \) quantifies the time decay of the \(L^p\) discrepancy between two solutions of the SVIE with different initial values. Since \(\varsigma\) is bounded (i.e. \(\| \varsigma^2 \|_\infty < +\infty\) ), if both \(\kappa < \frac{\lambda^2}{4\,p\, (1+ \varepsilon_*)\,\| \varsigma^2 \|_\infty \int_0^{+\infty} f_\lambda^2(u) \, du}\) and \(\big(1-\frac{(f_\lambda * \mu)}{\lambda m_0}\big)  \in \mathcal{L}^2 (\R_+)\), 
	then one derives from equation~\eqref{eq:varphi_fixed} and using Fubini-Tonelli's theorem that:
	{\small
	\begin{align*}
		\int_0^{+\infty} \varphi_{\infty,p}(s) \, ds &\leq \int_0^{+\infty} \varphi_{\infty,p}^{\varepsilon_*}(s) \, ds \leq \frac{\lambda^2 \Big(1+\frac{1}{\varepsilon_*}\Big)}{\lambda^2 - (C_p^{\mathrm{BDG}})^2\,\kappa (1+ \varepsilon_*) \| \varsigma^2 \|_\infty \|f_{\lambda}\|^2_{\mathcal{L}^2 (\R_+)}} \int_0^{+\infty} \big(1-\frac{(f_\lambda * \mu)(t)}{\lambda m_0}\big)^2 \, dt < \infty.
	\end{align*}
	}
	Note that, the condition \(\big(1-\frac{(f_\lambda * \mu)}{\lambda m_0}\big)  \in \mathcal{L}^2 (\R_+)\)  is in particular satisfied when \(\mu(t):=\lambda m_0(1 + \theta(t))\) with  \(\theta  \in \mathcal{L}^1(\R_+)\). Indeed, in this case, using~\eqref{eq:CondMean}, we obtain \(\big(1-\frac{(f_\lambda * \mu)}{\lambda m_0}\big)= R_{\lambda}+  f_{\lambda} *\theta\). Since \(R_{\lambda}, f_{\lambda}\in \mathcal{L}^2 (\R_+)\) for the exponential-fractional kernel whenever \(\alpha>\frac12\) (see e.g. Proposition~\ref{prop:main2}~$(a)$, $(b)$), the claim follows from the arithmetic--geometric mean inequality together with Young's inequality. 
	
	\smallskip
	2. {\em $L^2$-confluence:} Under the assumption of Proposition \ref{prop:confluenceLp} with $c\!\in \big(0,\frac{1}{\kappa}\big)$, $ \rho:=c\kappa $ and $X_0,X'_0\!\in L^2({\mathbb P})$. By \cite[Proposition 5.3]{EGnabeyeuPR2025} (which use It\^o's Isometry and the first Dini Lemma), one has that the solutions to Volterra equation~\eqref{eq:Volterrameanrevert} denoted $(X_t)_{t\ge 0}$ and $(X'_t)_{t\ge 0}$ starting from $X_0$ and $X'_0$ respectively satisfies:
	\begin{itemize}
	\item[(a)] There exists a continuous non-negative function \( \varphi_{\infty}^{\lambda,c, K, \mu}=: \varphi_{\infty}: \mathbb{R}_+ \to [0, 1] \), such that \( \varphi_\infty(0) = 1 \), \( \lim_{t \to +\infty} \varphi_\infty(t) =  \frac{\ell_\infty^2}{1-\rho(1-\ell_\infty^2)} \), only depending on \( \lambda, c, \mu \), and \( K \), such that :
	\[
	\forall t \geq 0, \; \E\, \Big(\Big| X_t - X^\prime_t\Big|\Big)^2 \leq \varphi_\infty(t) \E\, \Big(\Big| X_0 - X^\prime_0\Big|\Big)^2 \quad \; W_2([X_t'], [X_t]) \leq \varphi_\infty(t)^{1/2} W_2([X_0'], [X_0]).
	\]
	\item[(b)]  In particular, if \(\ell_\infty=0\) and \( X \) has a fake stationary regime of type I with parameters $(\mu_\infty, \lambda, c,m_0, v_0)$, then \( \mathbb{E} X_t' \to m_0 \), \( \text{Var}(X_t') \to v_0 \) as \(t \to +\infty\).
	And more generally finite-dimensional \( W_2 \)-convergence.
	Thus, the process \(X^\prime\) mixes: as time increases, the random variable  $X^\prime_t$ progressively forgets its initial mean and variance and converges to those of the limiting fake stationary regime.
	While, if X has a fake stationary regime of type II, its marginal distribution is unique.
	\end{itemize}
	\subsection{Asymptotics: Long run functional weak behaviour:}\label{subsect-Asymptotics}
	In the following,  $\stackrel{(C)}{\rightarrow }$ stands for functional weak convergence on $C(\R_+,\R)$ equipped with the topology of uniform convergence  on compact sets.
	\begin{Assumption}[Integrability and Uniform H\"older continuity]\label{ass:int_holregul}
	Let \( \lambda, c > 0 \), and assume the kernel \( K \) and its \( \lambda \)-resolvent \( R_{\lambda} \) satisfy
	\vspace{-.3cm}
	\begin{equation}\label{eq:Integf}
		\int_0^{+\infty} f_{\lambda}^{2\beta}(u) \, du < +\infty \quad \text{for some } \beta > 1.
	\end{equation}
	\vspace{-.3cm}
	and there exists \( \vartheta \in (0, 1] \), and a real constant \( C < +\infty \) such that
	\begin{equation}\label{eq:Regulf}
		 \left[ \int_0^{+\infty} |f_{\lambda}(u + \bar{\delta}) - f_{\lambda}(u)|^2 \, du \right]^{\frac{1}{2}} \leq C \bar{\delta}^{\vartheta}.
	\end{equation}
	\vspace{-.5cm}
		\end{Assumption}
		\begin{Theorem}\label{Thm:funcWeak}
	
	Let \( \lambda > 0 \), let \( \mu_\infty \in \mathbb{R} \), and let  \( \mu : \mathbb{R} \to \mathbb{R} \) a bounded bornel funtion,  \( \sigma : \mathbb{R} \to \mathbb{R} \) be a Lipschitz continuous function satisfying equation ~\eqref{eq:on_sigma} $(SL_{\sigma})$. Assume Assumption \ref{ass:int_holregul} and Assumption \ref{ass:on_stabilizer} on Equation \( (E_{\lambda}) \) are in force. Let \(c> \|\varsigma\|_{\infty}^2 \frac{\|f_{\lambda}\|^2_{\mathcal{L}^2(\R_+)}}{\lambda^2}\) and \( (X_t)_{t \geq 0} \) be the solution to the scaled Stochastic Volterra Integral Equation ~\eqref{eq:Volterrameanrevert} starting from \( X_0 \in L^p(\P) \)  for some suitable \( p \).
	
	\medskip
	\noindent {\em (a) C-tightness of time-shifted processes.} Assume
	\begin{equation}\label{eq:func_weak}
		X_0 \in L^p(\mathbb{P}) \quad \text{and} \quad
		\begin{cases}
			p = 2 \quad \text{and} \quad c < \frac{1}{\kappa} \quad \text{if} \quad ( \vartheta \wedge \frac{\beta - 1}{2\beta}) > \frac{1}{2}, \\
			p > \frac{1}{ \vartheta \wedge \frac{\beta - 1}{2\beta}} \quad \text{and} \quad c < \frac{1}{4p \kappa} \quad \text{if} \quad ( \vartheta \wedge \frac{\beta - 1}{2\beta}) \leq \frac{1}{2}.
		\end{cases}
	\end{equation}
	Then, the family of shifted processes \( (X_{t+u})_{u \geq 0} \) is C-tight and uniformly integrable, square uniformly integrable if \( p > 2 \) as \( t \to +\infty \). For any limiting distribution \( P \) on \( \Omega_0 := C(\mathbb{R}_+, \mathbb{R}) \), the canonical process \( Y_t(\omega) = \omega(t) \) has a \( \left(  \vartheta \wedge \frac{\beta - 1}{2\beta} - \frac{1}{p} - \eta \right) \)-H\"older pathwise continuous \( P \)-modification for sufficiently small \( \eta > 0 \). 
	That is, along a sequence \( t_n \uparrow \infty \), there exists a process \( X^{\infty} \) with continuous sample paths such that 
	\[(X_{u+t_n})_{t \geq 0} \Rightarrow (X^{\infty}_u)_{u \geq 0}
	\quad \text{weakly in } C(\mathbb{R}_+; \mathbb{R}) \text{ as } n \to \infty.\]
	\noindent Any limiting process \( X^{\infty} \) satisfies \( \forall t \ge 0, \quad X^{\infty}_t \in L^p(\P) \)  for each $p \geq 2$ and its first moment is given by
	\(\mathbb{E}[X^{\infty}_t] =
	\ell_\infty \mathbb{E}[X_0]  + (1-a) \frac{\mu_\infty}{\lambda}.\)
	
	\noindent Moreover, if $m_0= (1-a) \frac{\mu_\infty}{\lambda}$ i.e. $\ell_\infty=0$, the shifted processes of two solutions \( (X_t)_{t \geq 0} \) and \( (X_t')_{t \geq 0} \) are \( L^p \)-confluent, i.e. there exists a non-increasing function \( \bar{\varphi}_{\infty,p} : \mathbb{R}_+ \to [0, \frac{1}{(1-\sqrt{4pc\kappa})^2}] \) with \( \lim_{t \to +\infty} \bar{\varphi}_{\infty,p}(t) = 0 \) and
	\[W_p\left( \left[ (X_{t+t_1}, \ldots, X_{t+t_N}) \right], \left[ (X'_{t+t_1}, \ldots, X'_{t+t_N}) \right] \right) \leq \sqrt{N} \bar \varphi_{\infty,p} (t) {\cal W}_p\big([X_0],[X'_0]\big)\to 0 \quad \text{as} \quad t \to +\infty.\]
	Hence, the functional weak limiting distributions of \( [X_{t+\cdot}] \) and \( [X'_{t+\cdot}] \) coincide, meaning that if \( [X_{t_n+\cdot}] \stackrel{(C)}{\rightarrow} P \) for some subsequence \( t_n \to +\infty \), then \( [X'_{t_n+\cdot}] \stackrel{(C)}{\rightarrow} P \) and vice versa.
	
	\medskip
	\noindent {\em (b)  Functional weak long-run behavior.} Assume furthermore that the solution \((X_t)_{t \geq 0}\) of the volterra equation ~\eqref{eq:Volterrameanrevert} has a fake stationary regime of type I, starting from a random variable $X_0 \in L^2(P)$ with mean \( m_0 
	\) and variance \(v_0\).Then for any limiting distribution \(X^{\infty}\),  \( \mathbb{E}[X^{\infty}_t] =  m_0
	\) 
	while its autocovariance function is, for \( t_1, t_2 \geq 0 \), \( t_1 \leq t_2 \), given by \(\text{Cov}(X_{t_1}^\infty, X_{t_2}^\infty)=C_{f_{\lambda}}(t_1,t_2)\)
	{\small
		\begin{equation}\label{eq:funclongRun}
			\text{Cov}(X_{t+t_1}, X_{t+t_2}) \overset{t \to +\infty}{\to} C_{f_{\lambda}}(t_1,t_2) := \ell_\infty^2 v_0 + \frac{\varsigma_{\infty}^2 v_0 }{c\lambda^2}\int_0^{+\infty}  f_{\lambda}(t_2-t_1+u)f_{\lambda}(u)du.
		\end{equation}
	}
	$$\varsigma_{\infty}^2 =\begin{cases}
		1, & \text{in the setting of Theorem~\ref{thm:autonome}},\\
		\frac{c\lambda^2 (1-\ell_\infty^2)}{\int_0^{+\infty}f^2_{\lambda}(s)ds}, & \text{in the setting of Theorem~\ref{thm:timeautonome}}.
	\end{cases} $$
	Thus, under any functional limiting distribution \( P \), the canonical process \( Y \) is a (weak) \( L^2 \)-stationary process
	with mean \(  m_0 \) and covariance function \( C_{f_\lambda}(s,t) \), for \( s, t \geq 0 \).
	
	\medskip
	\noindent {\em (c) Stationary Gaussian Case.}  If \( \sigma(x) = \sigma > 0 \) is constant and \(X_0\) has a Gaussian distribution, (say \( X_0 \sim \mathcal{N}\left( m_0, v_0 \right) \) with \(v_0=c\sigma^2 \)), 
	then \( (X_t)_{t \geq 0} \) satisfies
	\(
	X_{t+\cdot} \stackrel{(C)}{\longrightarrow} \mathcal{GP}(f_{\lambda}) \quad \text{as} \quad t \to +\infty,\)
	where \( \mathcal{GP}(f_\lambda) \) is the 
	stationary Gaussian process with mean \(  m_0\)
	and covariance function \( C_{f_\lambda}(\cdot,\cdot) \).
	\end{Theorem}
	
	\noindent {\bf Proof:}
	It follows from ~\eqref{eq:func_weak} that either \( p = 2 \) and \( c < \frac{1}{\kappa} \), or \( p > 2 \) and \( c < \frac{1}{4\,p\,\kappa} \). Hence, Proposition \ref{prop:Momentctrl} implies that $ \sup_{t\ge 0}\Big\|  |X_t-m_0|\Big\|_p < +\infty$.
	As a consequence of \( \sigma \) having at most affine growth, we derive that \( \sup_{t \geq 0} \Big\| |\sigma(X_t)| \Big\|_p < +\infty \).
	
	\smallskip
	\noindent  {\sc Step~1.} {\em (Kolmogorov criterion).} 
	Now, we can establish C-tightness by the Kolmogorov criterion. Let \( p\geq2 \) and $c$ satisfying constraint ~\eqref{eq:func_weak}. One writes  for $s,\,t\ge 0$ with $s\le t$ and owing to equation~\eqref{eq:ConstMean}:
	\begin{equation}\label{eq:Tight} X_t-X_s = \big(X_0 - m_0\big) \Big( I(s) -I(t) \Big) + \frac{1}{\lambda}\Big(J(t)-J(s)\Big).
	\end{equation}
	where we set $J(t):=\int_0^tf_{\lambda}(t-u)\varsigma(u)\sigma(X_{u})dW_u$ and $ I(t) :=  \frac{1}{\lambda m_0}\int_0^t f_{\lambda}(t-u) \mu(u)  \, du$. 	
	
	First, owing to the equation~\eqref{eq:MasterEq}, we have \(\big|I(t)-I(s)\big|  = 	\big|h_{\lambda,c}(t)-h_{\lambda,c}(s)\big|\). Moreover, by Young's inequality \((f_\lambda^2 * \varsigma^2)(t) \le \|f_{\lambda}\|^2_{\mathcal{L}^2(\R_+)}\,\|\varsigma\|_\infty^2,\) so that under our assumption, we obtain the uniform lower bound \(h_{\lambda,c}(t)\geq h_{*}:=
	\sqrt{
		1-
		\frac{\|f_{\lambda}\|^2_{\mathcal{L}^2(\R_+)}\,\|\varsigma\|_{\infty}^2}
		{c\lambda^2}
	}
	>0\) for every \(t\geq0\).
	Therefore, 	on the first hand,
	by using the triangle inequalities, one gets:
	\begin{align*}
		\big|I(t)&-I(s)\big|  = 	\big|h_{\lambda,c}(t)-h_{\lambda,c}(s)\big| \le \frac{\|\varsigma\|_{\infty}^2}{2 c \lambda^2 h_*}  \Big(\int_s^{t}f_{\lambda}^2(u)\, du  +  \int_0^s \big| f_{\lambda}^2(t-s + u) - f^2_{\lambda}(u) \big|\, du \Big)  \\
		&\le \frac{\|\varsigma\|_{\infty}^2}{2 c \lambda^2 h_*}   \times \Big( \big(\int_0^{+\infty} f_{\lambda}^{2\beta}(u)du \big)^{\frac{1}{\beta}}  (t-s)^{1-\frac{1}{\beta}} + 2 \|f_{\lambda}\|^2_{\mathcal{L}^2(\R_+)} \int_0^s | f_{\lambda}(t-s + u) - f_{\lambda}(u) |^2\,du \Big)\\
		&\le C_{X_0, \lambda, \beta, \varsigma, f_{\lambda}}\cdot\Big(|t-s|^{(1-\frac{1}{\beta})}+ |t-s|^{\vartheta}\Big).
	\end{align*}
	where, in the penultimate inequality, we used the elementary identity \(|a^2-b^2|\leq|a-b|\big(|a|+ |b|\big) \),
	together with the Cauchy--Schwarz inequality, while the last estimate follows from Assumption~\ref{ass:int_holregul}.
	It boilds down that the first term of the sum on the right hand side of the above equality~\eqref{eq:Tight} provide the bound:
	\[ \Big\| \big(X_0 - m_0\big) \Big( I(s) -I(t) \Big) \Big\|_p = \Big\| |X_0 - m_0 |\Big\|_p\Big|I(t)-I(s)\Big| \leq C_{X_0,p, \lambda, \beta, \varsigma, f_{\lambda}}|t-s|^{ (\vartheta\wedge\frac{\beta-1}{\beta})}\]
	where where \(C_{X_0, p,\lambda, \beta, \varsigma, f_{\lambda}}<+\infty\) owing to 	Equations~\ref{eq:Integf}-\eqref{eq:Regulf} of assumption~\ref{ass:int_holregul}. 
	On the other hand, combining the $L^p$-BDG and the generalized Minkowski inequality, one derives from~\eqref{eq:func_weak} that,
	\begin{align*}
		&\ \Big\| |J(t)-J(s)| \Big\|_p \le \left\| |\int_s^{t}f_{\lambda}(t-u)\varsigma(u)\sigma(X_{u})\mathrm{d}W_u |\right\|_p + \left\| |\int_0^s \left( f_{\lambda}(t-u) - f_{\lambda}(s-u) \right)\varsigma(u)\sigma(X_{u})dW_u | \right\|_p\\
		&\leq C_p \|\varsigma^2\|_{\infty} \left[\Big(\int_s^tf^2_{\lambda}(t-u)\big\||\sigma(X_{u})|\big\|_p^2 du\Big)^{\frac12} + \Big(\int_0^s \big(f_{\lambda}(t-u)-f_{\lambda}(s-u)\big)^2 \big\||\sigma(X_{u})|\big\|_p^2 du\Big)^{\frac12}\right]\\
		&\leq C_p \|\varsigma^2\|_{\infty} \sup_{u\ge 0}\Big\||\sigma(X_{u})|\Big\|_{p} \left[\Big(\int_0^{+\infty}f^{2\beta}_{\lambda}(u) du\Big)^{\frac{1}{2\beta}} |t-s|^{\frac{\beta-1}{2\beta}} + \left(\int_0^{+\infty} \big(f_{\lambda}(t-s+u)-f_{\lambda}(u)\big)^2du\right)^{\frac12}\right] \\
		&\hspace{2cm}\leq  C_{p, \sigma,  \varsigma, \lambda, f_{\lambda}, X_0} \cdot |t-s|^{ \vartheta\wedge\frac{\beta -1}{2\beta}} \quad \text{where} \quad C_p  \equiv C^{BDG}_p .
	\end{align*}
	Finally, putting all these estimates together,  since $ \frac{\beta-1}{2\beta}\le \frac{\beta-1}{\beta} $ we have the existence of a real constant $0<C_{p, X_0, \varsigma, \beta,\lambda, f_{\lambda}}<\infty$ such that:
	\[\E\, \left(|X_t-X_s| \right)^p \le C_{p, X_0, \varsigma, \beta,\lambda, f_{\lambda}}|t-s|^{p(\vartheta\wedge \frac{\beta-1}{2\beta})}\]
	Define for \( u \geq 0 \) the process \( X^u \) by \( X^u_t = X_{t+u} \), where \( t \geq 0 \). Then \( X^u \) has continuous sample paths and satisfies
	
	\centerline{$
		\sup_{u \geq 0} \mathbb{E}[|X^u_t - X^u_s|^p] \leq C(p)|t-s|^{p(\vartheta\wedge \frac{\beta-1}{2\beta})}
		\quad \text{for } 0 \leq t - s \leq 1.
		$}
		
	\smallskip
	As $p( \vartheta \wedge\frac{ \beta -1}{2\beta})>1$ according to equation \eqref{eq:func_weak}, it follows from  Kolmogorov's $C$-tightness criterion (see ~\cite[Theorem 2.1, p. 26, 3rd edition]{RevuzYor} 
	or \cite[Lemma 44.4, Section IV.44, p.100]{RogersWilliamsII}), that the family of  shifted processes $X_{t+\cdot}$, $t\ge 0$, is $C$-tight i.e. \( (X^u)_{u \geq 0} \) is tight on \( C(\mathbb{R}_+; \mathbb{R}) \) (hence the existence of a weak
	continuous accumulation point thanks to Prokhorov's theorem) with limiting distributions P under which the canonical process has the announced H\"older pathwise regularity. Therefore, we conclude that along a sequence \( u_k \uparrow \infty \), the process \( X^{u_k} \) converges in law to some continuous process \( X^{\infty} \).
	
	An application of Fatou's lemma shows that any limiting process (resp. the limit distribution) has a finite moment of any order, i.e., $\quad 
	\forall t>0, \quad \mathbb{E}[|X^{\infty}_t|^p] \leq \sup_{u \geq 0} \mathbb{E}[|X_u|^p] < \infty.$
	
	For the first moment formula, we note using equation~\eqref{eq:espVolterramean} and Lemma \ref{lem:on_mu} that
	
	\[
	\mathbb{E}[X_t] \longrightarrow \ell_\infty \mathbb{E}[X_0]  + (1-a) \frac{\mu_\infty}{\lambda} \quad \text{as } t \to \infty.
	\]
	Since \( \sup_{t \geq 0} \mathbb{E}[|X_t|^2] < \infty \), we easily conclude that $\quad 
	\lim_{t \to \infty} \mathbb{E}[X_t] = \mathbb{E}[X^{\infty}_t].$
	
	\bigskip
	\noindent{\sc Step 3}. {\em (b) Asymptotic weak stationarity. }
	Now let us consider the asymptotic  covariance between $X_{t+t_1}$ and $X_{t+t_2}$, $0<t_1<t_2$ when $X_t$ starts for $X_0$ with mean  \( m_0 
	\), variance $v_0$ and $\bar \sigma^2 = \E\, \sigma(X_{t})^2$, $t\ge 0$ constant over time. Let \(\ell_\lambda(t):=\big(1-\frac{(f_\lambda * \mu)(t)}{\lambda m_0}\big)\) for short.
	Using $\text{Cov}(aU + b, cV + d) = ac \, \text{Cov}(U, V)$ and equation ~\eqref{eq:Volterrameanrevert_}, we have:
	{\small
		\begin{align*}
			{\rm Cov}(X_{t+t_1}, X_{t+t_2}) &= {\rm Var}(X_0) \ell_\lambda(t+t_1)\, \ell_\lambda(t+t_2)+ \frac{1}{\lambda^2}\E\left[\int_0^{t+t_1} f_{\lambda}(t+t_2-s)f_{\lambda}(t+t_1-s)\varsigma^2(s)\sigma^2(X_{s})ds\right]\\
			&=  {\rm Var}(X_0)\ell_\lambda (t+t_1)\ell_\lambda (t+t_2) + \frac{\bar \sigma^2}{\lambda^2}\int_0^{t+t_1} f_{\lambda}(t_2-t_1+u)f_{\lambda}(u)\varsigma^2(t+t_1-u) du.
		\end{align*}
	}
	As $ f_{\lambda}(t_2-t_1+\cdot)f_{\lambda}\!\in {\cal L}^2(\R_+)$ since $f_{\lambda}\!\in {\cal L}^2(\R_+)$ , $\mbox{\bf 1}_{\{0\le u \le t+t_1\}}\varsigma^2(t+t_1-u)\to  \varsigma^2_\infty$ for every $u\!\in \R_+$ as $t\to +\infty$ (owing to Lemma \ref{lem:asymptotique}) and $\lim_{t\to+\infty}\ell_\lambda(t)= \ell_\infty$, we have:
	\[{\rm Cov}(X_{t+t_1}, X_{t+t_2})\stackrel{t\to+\infty}{\longrightarrow} \ell_\infty^2 {\rm Var}(X_0) + \frac{ \bar \sigma^2 \varsigma^2_\infty }{\lambda^2}\int_0^{+\infty}  f_{\lambda}(t_2-t_1+u)f_{\lambda}(u)du.\]
	Using \(v_0 = c \bar \sigma^2\), we obtain \(C_{f_{\lambda}}(t_1,t_2)\) in Equation~\eqref{eq:funclongRun}.
	The confluence result follows from Proposition \ref{prop:confluenceLp} with \( \bar \varphi_{\infty,p} (t) = \sup_{u \geq t} \varphi_{\infty,p}^{1/2} (u) \).
	In fact, if $X$ and $X'$ are two solutions of  Equation~\eqref{eq:Volterrameanrevert_} starting from $X_0$ and $X'_0$ respectively, both square integrable, using 
	Proposition~\ref{prop:confluenceLp}, we derive that for every $0\le t_1<t_2< \cdots < t_{_N}<+\infty$
	\[{\cal W}_p\big([(X_{t+t_1}, \cdots, X_{t+t_{_N}})], [(X'_{t+t_1}, \cdots, X'_{t+t_{_N}})])\leq \sqrt{N} \bar \varphi_{\infty,p} (t) {\cal W}_p\big([X_0],[X'_0]\big)\to 0 \mbox{ as }t\to +\infty.\]	
	\noindent As  a consequence, the  weak limiting distributions  of $[X_{t+\cdot}]$ and $[X'_{t+\cdot}]$ are the same in the sense that,  if $[X_{t_n+\cdot}]\stackrel{(C)}{\longrightarrow} P$ for some subsequence $t_n \to +\infty$ (where $P$ is a probability measure on $C(\R_+, \R)$ equipped with the Borel $\sigma$-field induced by the sup-norm topology), then  $[X'_{t_n+\cdot}]\stackrel{(C)}{\longrightarrow} P$ and  conversely.
	
	\bigskip
	\noindent{\sc Step 4}. {\em (c) Stationary Gausian case.} This result stems first from the fact that \( (X_t)_{t \geq 0} \) is a Gaussian process, implying that its limiting distributions in the functional weak sense are also Gaussian. Secondly, a Gaussian process is completely characterized by its mean and covariance functions.\\
	In fact, when \( \sigma(x) = \sigma > 0 \quad \forall x \in \mathbb{R}\) and \( X_0 \) follows a Gaussian distribution, the process \( X \) is Gaussian, which implies (at least for finite-dimensional weak convergence, i.e., weak convergence of all marginals of any order) that,$
	(X_{t + \cdot}) \stackrel{(C)}{\longrightarrow} \mathcal{GP}(f_{\lambda}) \quad \text{as} \quad t \to +\infty,
	$
	where \( \mathcal{GP}(f_{\lambda}) \) is a Gaussian process with mean \( m_0 \) and covariance function given above. \hfill $\Box$
	
	\medskip
	\noindent{\bf Remarks:} 1. Be aware that at this stage, we do not have uniqueness of the limit distributions since they are
	not characterized by their mean and covariance functions, except in Gaussian setting.\\
	2. To establish relative compactness in Theorem~\ref{Thm:funcWeak}~$(a)$, in terms of functional ${\cal W}_p$-compactness ( Wasserstein $p-$distance), we would require that \( \|\sup_{t \geq 0} |X_t|\|_q < +\infty \) for some \( q > p \).
	
	\section {Applications to $\rho$-Exponential  $\alpha$-Fractional 
		Volterra SDEs 
		}\label{sec:appl3}
	Let consider the below {\em Gamma Fractional integration kernel} or {\em Exponential-Fractional integration kernel} defined in Example \ref{Ex:SolventGammaKernel}, where $\alpha = H + \frac{1}{2}$, with $H$ denoting the Hurst coefficient:
	
	\centerline{$
		K(t) = K_{\alpha, \rho}(t) = e^{-\rho t} \frac{u^{\alpha - 1}}{\Gamma(\alpha)} \mathbf{1}_{\mathbb{R}_+}(t), \quad  \text{with} \quad \alpha, \rho > 0.
		$}
Note that, this is a generalization of the exponential and the fractional integration kernels (see e.g. the Remark following Example \ref{Ex:SolventGammaKernel}).  This case is of special interest since it corresponds to the recent introduction of fractional
Volterra SDEs to devise ``rough models'' of stochastic volatility dynamics in Finance ( see e.g., ~\cite{el2019characteristic,ElEuchR2018,GatheralJR2018,GaJuRo2020}).
	The gamma kernel is often adopted in the quadratic rough Heston model (see, e.g., \cite{BourgeyGatheral2025}) due to its numerical convenience, flexibility, and the availability of a closed-form expression for its resolvent of the second kind. 
	We show in this part that for such kernels $K_{\alpha, \rho}$, the resolvent $R_{\alpha,\rho,\lambda}$  satisfy our standing assumption $({\cal R}_\lambda)\; $ for all $\lambda > 0$, and that $f_{\alpha,\rho, \lambda} := -R_{\alpha,\rho,\lambda}$ exists and is square-integrable with respect to the Lebesgue measure on $\mathbb{R}_+$, for $\frac{1}{2} < \alpha < 1$ (``rough models''). As a consequence, the results and findings from Sections \ref{sec:genfakestatio}, \ref{sect:ExamplesFakeI-II} and \ref{sect-LongRunB}, will be applicable in the cases where $\sigma(t, x) = \sigma(x)$ and $\sigma(t, x) = \varsigma_{\lambda,c}(t)\sigma(x)$. 
	
	\subsection{ Exponential $\alpha-$ Fractional  kernels $0 <\alpha<1$}\label{subsec:Gammaalphafrac}
	Let  $\lambda > 0$. Taking the Laplace transform of \eqref{eq:Resolvent} and using Example~\ref{Ex:SolventGammaKernel}, we obtain, \(L_{R_{\alpha, \rho, \lambda}}(s) = \frac{1}{s (1 + L_{K_{\alpha, \rho}}(s))}=\frac{1}{s (1 + \lambda (s + \rho)^{-\alpha})}\). Hence, by the Tauberian Final Value Theorem,
	\(a:=\lim_{t \to \infty} R_{\alpha, \rho, \lambda}(t) = \lim_{s \to 0} sL_{R_{\alpha, \rho, \lambda}}(s) = \frac{1}{1 + \lambda \rho^{-\alpha}} \in [0, 1)\). Moreover,
	on $(0,\infty)$, the function $f_{\alpha, \rho, \lambda}:= - R_{\alpha, \rho, \lambda} $  introduced in Equation~\eqref{eq:DerivSolventGammaKernel} of Example~\ref{Ex:SolventGammaKernel}, satisfies
	{\small
		$$
		L_{f_{\alpha, \rho, \lambda}}(s) = L_{-R'_{\alpha, \rho, \lambda}}(s) = -s L_{R_{\alpha, \rho, \lambda}}(s) + R_{\alpha, \rho, \lambda}(0) = \frac{-s}{s (1 + \lambda (s + \rho)^{-\alpha})} + 1
		= \frac{\lambda}{\lambda + (s + \rho)^{\alpha}} = L_{e^{-\rho \cdot}f_{\alpha, \lambda}}(s)
		$$
	}
	Hence, by the injectivity of the Laplace transform, \(	f_{\alpha, \rho, \lambda}(t) = e^{-\rho t} f_{\alpha, \lambda}(t) =\alpha \lambda e^{-\rho t} t^{\alpha-1} E^\prime_{\alpha}( - \lambda t^{\alpha}).\) 
	Likewise, using Tauberian Final Value Theorem, \(\lim_{t \to \infty} f_{\alpha, \rho, \lambda}(t) = \lim_{s \to 0} s L_{-R^\prime_{\alpha, \rho, \lambda}}(s)\), that is
	\vspace{-.3cm}
	\begin{equation}
		\lim_{t \to \infty} f_{\alpha, \rho, \lambda}(t) = -\lim_{s \to 0} s\left(s L_{R_{\alpha, \rho, \lambda}}(s)-R_{\alpha, \rho, \lambda}(0)\right) = -\lim_{s \to 0} \frac{s}{(1 + \lambda (s + \rho)^{-\alpha})} - s = 0.
		\vspace{-.2cm}
	\end{equation}
	\begin{Proposition}\label{prop:main2} Let $\lambda>0$, \(\rho\geq0\) and $\alpha\!\in (0,1)$.
		
		\smallskip
		\noindent $(a)$ The $\lambda$-resolvent $R_{\alpha,\rho, \lambda} $ is completely monotonic (hence infinitely differentiable i.e. $\mathcal C^\infty$ on $(0, +\infty)$).
		Moreover $R_{\alpha,\rho, \lambda}(0)=1$, $R_{\alpha,\rho, \lambda}$ decreases  to $a:=\frac{1}{1 + \lambda \rho^{-\alpha}} \in [0,1[$.
		$R_{\alpha,\rho, \lambda}\!\in {\cal L}^{r}(\R_+)$ for every $r>\frac{1}{\alpha}$.  
		
		\smallskip
		\noindent $(b)$ $f_{\alpha,\rho, \lambda} := -R'_{\alpha, \rho, \lambda}$ is a completely monotonic function (hence convex), decreasing to \(0\) and satisfy :
		{\small
			\[ \forall\, t>0, \quad
			f_{\alpha,\rho, \lambda} (t):= e^{-\rho t}f_{\alpha, \lambda} (t) = \lambda^{\frac{1}{\alpha}}\int_0^{+\infty} e^{-(\rho+\lambda^{\frac{1}{\alpha}})tu}uH_{\alpha}(u)du,\; H_{\alpha}(u)= \frac{\sin(\alpha\pi)}{\pi}\frac{u^{\alpha-1}}{u^{2\alpha}+2u^{\alpha}\cos(\pi \alpha)+1}.
			\]
		}
		\smallskip
		\noindent Finally, \(1-R_{\alpha,\rho, \lambda}\) is a Bernstein function.
		
		\noindent $(c)$ If  $\alpha\!\in (\frac12, 1)$,  $ f_{\alpha,\rho, \lambda}$ is ${\cal L}^{2\beta}$-integrable $\forall \beta \in \big( 0, \frac{1}{2(1-\alpha)}\big)$.
		Moreover, for every $\vartheta\!\in\big(0,\alpha-\frac{1}{2}\big)$,  there exists a real constant $C_{\vartheta,\rho, \lambda}>0$ such that 
		
		\centerline{\(
			\forall\, \delta >0, \quad \left[\int_0^{+\infty} \big(f_{\alpha,\rho, \lambda}(t+\delta)-f_{\alpha,\rho,\lambda}(t) \big)^2 \right]^{\frac12}\le C_{\vartheta,\rho, \lambda}\delta^{\vartheta}.\)
		}
	\end{Proposition}
	\noindent \noindent In particular, the resolvent \(R_{\alpha,\rho,\lambda}\) is continuous and bounded on \([0,+\infty)\). For clarity and conciseness, the proof is postponed to Appendix \ref{app:lemmata}.
		\begin{figure}[H]
		\centering
		\includegraphics[width=.902\textwidth]{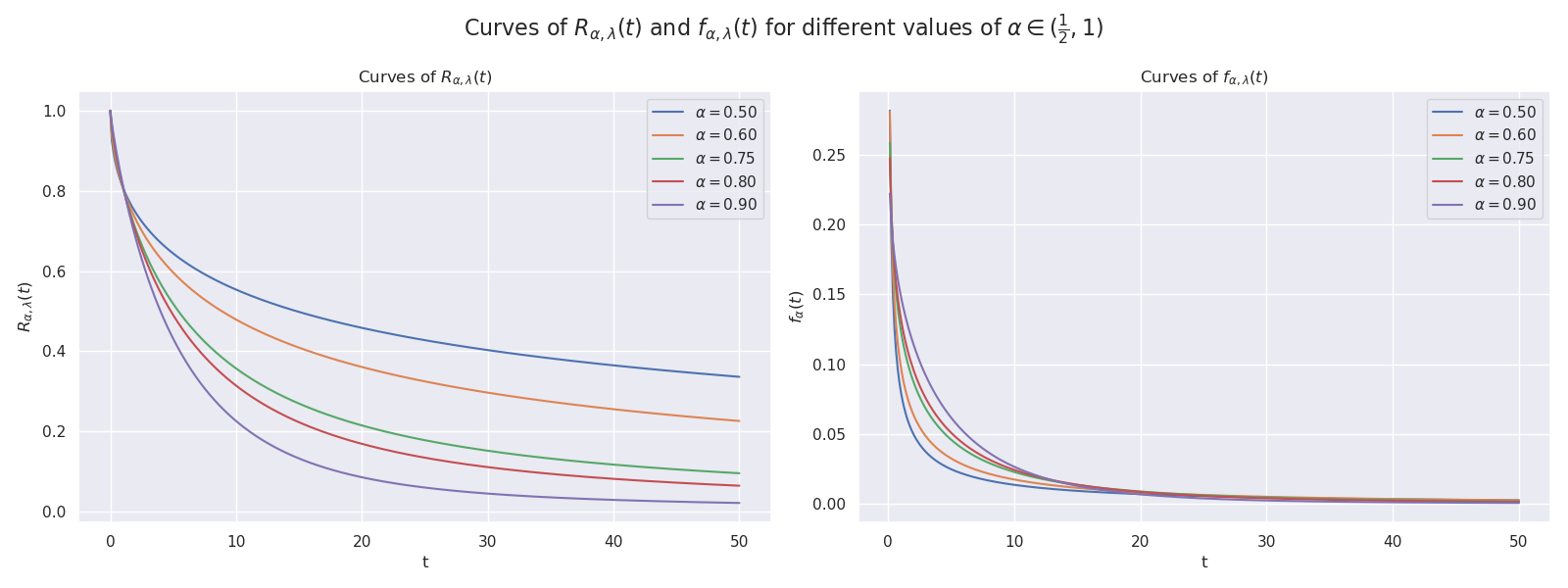} 
		\caption{Curves of $R_{\alpha,\lambda}(t)$ and $f_{\alpha,\lambda}(t)$ for different values of $\alpha \in [\frac{1}{2},1)$, $\lambda= 0.2$ and
			$\rho = 0$}\label{fig:curves_alpha_less1}
	\end{figure}
	\begin{figure}[H]
		\centering
		\includegraphics[width=.902\textwidth]{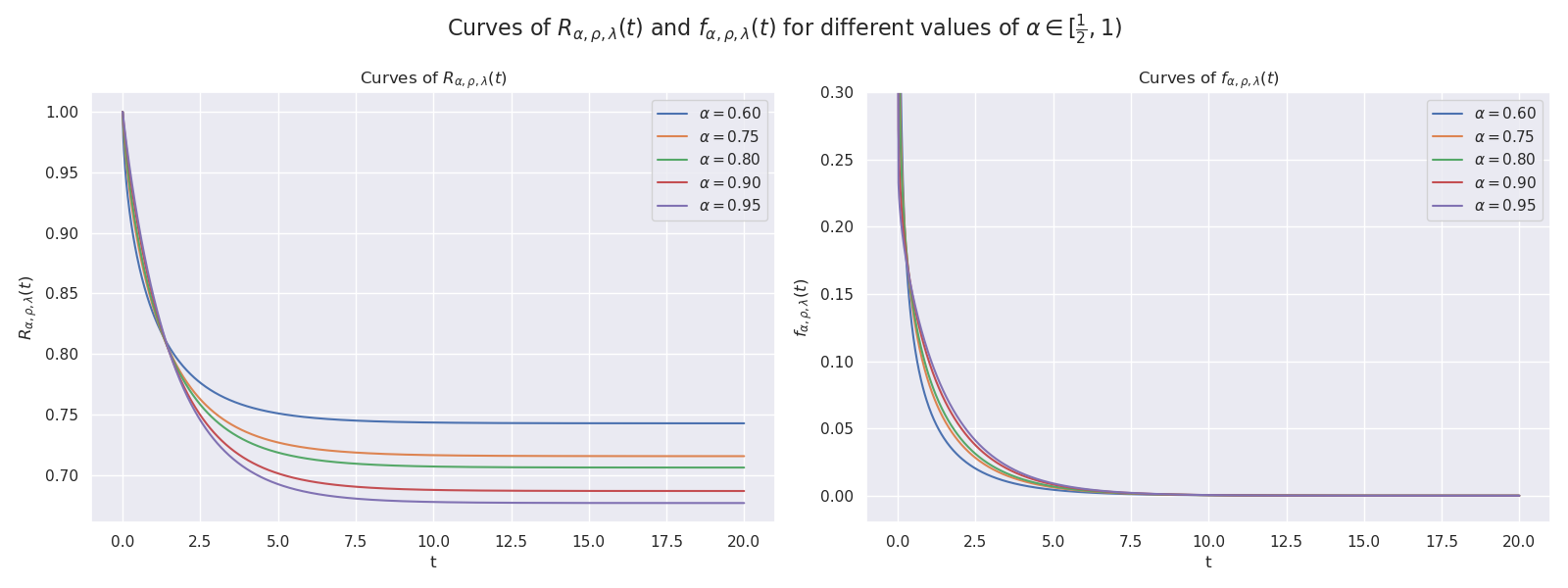} 
		\caption{Curves of $R_{\alpha,\rho,\lambda}(t)$ and $f_{\alpha,\rho,\lambda}(t)$ for different values of $\alpha \in [\frac{1}{2},1)$, $\lambda= 0.2$ and
			$\rho = 0.4$}\label{fig:curves_alpha_less1}
	\end{figure}
	\begin{Theorem}
		Let \( \alpha \in \left( \frac12, 1\right) \) , \(\rho>0\), let \( K(t) = K_{\alpha,\rho}(t) =  e^{-\rho t}\frac{t^{\alpha-1}}{\Gamma(\alpha)} \), \( t > 0 \) the Gamma fractional kernel, let $\sigma(t,x):= \varsigma(t)\sigma(x)$  with \( \sigma \) be a Lipschitz continuous function given by~\eqref{eq:sigmafakeI&II}, thus satisfying a relation of the type ~\eqref{eq:on_sigma} $(SL_{\sigma})$ with \(\kappa:=\kappa_2>0\). Let $p>\frac{2}{2\alpha-1}>2$ and  \( c \in \left( 0, \frac{2\alpha-1}{8\kappa_2} \right) \) with $\varsigma = \varsigma_{\lambda, c}$, \(\lambda>0\) and let \( X_0 \in L^p(\P) \) for some $p>\frac{2}{2\alpha-1}$ such that \( \mathbb{E}[X_0] = m_0\) and \(Var(X_0) = v_0=\frac{c\sigma^2( m_0)}{1-c\kappa_2}.\)
		Then,
		\begin{enumerate}
			\item For exponential-fractional kernels \( K_{\alpha,\rho} \) with \( \frac12 < \alpha < 1 \), the solution $(X_t )_{t\geq0}$ to the Volterra equation~\eqref{eq:ConstMeanVar} or ~\eqref{eq:ConstMean} starting from $X_0$ has a fake
			stationary regime of type I in the sense that:
			
			\centerline{$
				\forall\, t\ge 0, \qquad \E\, X_t = m_0,  \quad{\rm Var}(X_t) = v_0=\frac{c\sigma^2( m_0)}{1-c\kappa_2} \mbox{ and }\quad \E\, \sigma^2(X_t) =  \bar \sigma_0^2 =\frac{\sigma^2( m_0)}{1-c\kappa_2}.
				$}
			
			\item If $m_0= (1-a) \frac{\mu_\infty}{\lambda}$ i.e. $\ell_\infty=0$, then for every \( X_0^\prime  \in L^2(\P) \), a solution to the Volterra Equation~\eqref{eq:ConstMeanVar} or ~\eqref{eq:ConstMean} starting from \( X_0^\prime \) satisfies \( \| X_t^\prime - X_t \|_2 \to 0\) as $t \to \infty$.
			
			\item The family of shifted processes \( X_{t+\cdot}, t \geq 0 \), is \( C \)-tight as \( t \to +\infty \) and its (functional) limiting distributions are all \( L^2 \)-stationary processes with covariance function \( C_\infty \) given by ~\eqref{eq:funclongRun}.
		\end{enumerate}
	\end{Theorem}
	\noindent {\bf Proof.}	Let \( \alpha \in \left( \frac12, 1\right) \).
	By Proposition~\ref{prop:main2}~$(c)$, we let \(0<\vartheta<\alpha-\frac12\) and \(	1<\beta<\frac{1}{2(1-\alpha)}\). We have
	that conditions~\eqref{eq:func_weak} of Theorem~\ref{Thm:funcWeak} are satisfied for \(\frac{2}{2\alpha-1}<p < \frac1{4c\kappa_2}\) which we may assume, without loss of generality. In fact, in this case
	\[
	\lim_{\substack{\vartheta\to \alpha-\frac12,
			\beta\to \frac{1}{2(1-\alpha)}}}
	p\big(
	\frac{\beta-1}{2\beta}\wedge \vartheta
	\big)
	=
	\frac{p(2\alpha-1)}{2}
	>1.
	\]
	Therefore,	Theorem \ref{Thm:funcWeak} applies.
	The claim (1) is a consequence of Proposition \ref{prop:QuadraticSigma}. \hfill $\square$
	
		\subsection{
		Asymptotic expansion of $\varsigma_{\alpha,\rho,\lambda,c}^2$ solution of the stabilizer equation when $\alpha \in (\frac12,1)$}\label{sec:sigma2rough}
	In this section, we derive the asymptotic or leading-order behaviour of $\varsigma_{\alpha,\rho,\lambda,c}$ near the origin, assuming that \(\mu(t) \sim \mu(0)\) as \(t\to 0\) with \(\frac{\mu(0)}{m_0}\geq0\).  we drop the dependence in \(\alpha\) and \(\rho\) for simplicity.
	To this end we rely on the Laplace version of the equation \textit{($E_{\lambda, c}$)} in ~\eqref{eq:VolterraStabilizer} satisfied by $\varsigma^2_{\lambda,c}$, namely \(	\forall\, t\ge 0, \quad c \lambda^2 \big(1- \ell_\lambda^2(t) \big) =  (f_{\lambda}^2 * \varsigma^2)(t),\) where we set \(\ell_\lambda(t):=\big(1-\frac{(f_\lambda * \mu)(t)}{\lambda m_0}\big)\) for short. Its Laplace transform is given by:
	\(c\lambda^2 L_{1 - \ell_\lambda^2} = L_{f_{\lambda}^2} L_{\varsigma^2}.\)
	In order to get rid of the Laplace transform of \( 1 - \ell_\lambda^2 \), we apply integration by parts and obtain:
	\begin{equation}
		L_{1 - \ell_\lambda^2}(t) = L_1(t) - L_{\ell_\lambda^2}(t) = \frac{1}{t} - \Big( \frac{1}{t}\big(\ell_\lambda^2(0)- \lim_{u \to +\infty} e^{-tu}\ell_\lambda^2(u)\big) + \frac{2}{t} L_{\ell_\lambda \ell_\lambda^\prime}(t) \Big)= -\frac{2}{t} L_{\ell_\lambda \ell_\lambda^\prime}(t).
	\end{equation}
	Consequently, it follows that 
	\begin{equation}\label{eq:Laplacesigma}\forall\, t>0, \;  t\,L_{f^2_\lambda}(t).L_{\varsigma^2}(t)= - 2\,c\lambda^2 L_{\ell_\lambda\ell_\lambda^\prime}(t).
	\end{equation}
	Now, owing to Assumption~\eqref{ass:MeanRevert}, applying the Leibniz rule for differentiating an integral along with the fact that the space \( (\mathcal{L}^1(\R), +, \cdot, *) \) is a commutative algebra, one checks that $\ell_\lambda^\prime(t) = -\frac{1}{\lambda m_0} \left(\mu(0)f_{\lambda}(t) + (f_{\lambda}*\mu^\prime)_t \right)$ so that \(\ell_\lambda^\prime(t)\stackrel{0}{\sim} -\frac{\mu(0)}{\lambda m_0}f_{\lambda}(t)\).
	Since \(\ell_\lambda(t) \stackrel{0}{\sim} 1\) and given the form of the kernel $K_{\alpha, \rho}(u) = e^{-\rho u} \frac{u^{\alpha - 1}}{\Gamma(\alpha)} \mathbf{1}_{\mathbb{R}}(u), \;  \alpha, \rho > 0$ and the expansion of the resolvents $R_{\alpha,\rho, \lambda}$ and minus its derivative $-f_{\alpha,\rho, \lambda}$, in Example \ref{Ex:SolventGammaKernel}, 
	we derive  that:
	\[
	e^{2 \rho t}\ell_\lambda\ell_\lambda^\prime(t)\stackrel{0}{\sim} -\frac{\mu(0)}{\lambda m_0}\frac {\lambda t^{\alpha-1}}{\Gamma(\alpha)} \quad \mbox{ and}\quad e^{2 \rho t}f^2_{\lambda}(t)\stackrel{0}{\sim} \frac {\lambda^2 t^{2(\alpha-1)}}{\Gamma(\alpha)^2}\]
	It follows by a dual version of the  Hardy-Littlewood Tauberian Theorem (see for example \cite{FellerII} or \cite{BiGoTe1989}) that--at least heuristically-- 
	\(
	L_{e^{2 \rho \cdot}\ell_\lambda\ell_\lambda^\prime}(t)\stackrel{+\infty}{\sim} -\lambda \frac{\mu(0)}{\lambda m_0}t^{-\alpha} \; \mbox{ and}\; L_{e^{2 \rho \cdot}f^2_{\lambda}}(t)\stackrel{+\infty}{\sim} \frac {\lambda^2\Gamma(2\alpha-1) t^{-(2\alpha-1)}}{\Gamma(\alpha)^2}.
	\)
	Therefore, this implies that 
	\[
	L_{e^{2 \rho \cdot}\varsigma^2}(t) = L_{\varsigma^2}(t-2\rho) = \frac{-2c\lambda^2L_{\ell_\lambda\ell_\lambda^\prime}(t-2\rho)}{(t-2\rho)L_{e^{2 \rho \cdot}f^2_{\lambda}} (t-2\rho)} \stackrel{+\infty}{\sim} 2\lambda \,c \frac{\mu(0)}{\lambda m_0} \frac{\Gamma(\alpha)^2}{\Gamma(2\alpha-1)} t^{-(2-\alpha)} \]
	owing to Equation~\eqref{eq:Laplacesigma} evaluated at \((t-2\rho)\). 
	This in turn suggests that 
	{\small
		\begin{equation}\label{eq:varsigma2}
			\varsigma^2(t) \stackrel{t\to0}{\sim} \frac{2\lambda c \Gamma(\alpha)^2}{\Gamma(2\alpha-1) \Gamma(2-\alpha)} \frac{\mu(0)}{\lambda m_0} e^{-2 \rho t}t^{1-\alpha} \quad \text{so that}\quad  \varsigma^2(0) = 0.
		\end{equation}
	}
	
	\medskip
	\noindent {\bf Remark:\,} 
	For the computation of the function $\varsigma_{\alpha,\rho, \lambda, c}^2$, establishing a recurrence formula 
	as in~\cite{Pages2024}
	turns out to be quite challenging. This involves knowing the {\em form of the mean-reverting function \(\mu\)}. We rather solve the functional equation numerically, using for e.g. a second-order extension of the triangular Volterra discretization as established in Equation~\eqref{eq:TriangStabilScheme} further on.
	\subsection{Numerical illustration of fake stationarity for exponential $\alpha-$fractional SVIE with Lipschitz  diffusion coefficient and $\alpha \in (\frac12,1)$}
	\noindent
	In this section we specify a family of scaled volterra equations~\eqref{eq:Volterrameanrevert} for some \( \lambda > 0 \) and a diffusion coefficient $\sigma$ to be specified.
	 The reader is referred to Appendix~\ref{app:simu_det} for
	further details about the semi-integrated Euler scheme  introduced in this context for the simulation of the stochastic Volterra Equation~\eqref{eq:Volterrameanrevert_}
and to the captions of the differents figures for the numerical values of the parameters of the Stochastic Volterra equation.
	
	\subsubsection{A numerical illustration of intrinsic fake stationary SVIE with $\alpha$-fractional kernels for $\alpha \in (\frac12, 1)$ and quadratic  diffusion coefficient}
	We consider an $\alpha$-fractional kernel for $ \alpha \in (\frac12, 1)$ (``rough models ")  and a squared trinomial diffusion coefficient \(\sigma(t,x) = \sigma(x)\) with
	\(\sigma\) of the form \ref{eq:sigmafakeI&II} and given by \(\sigma(x) = \sqrt{ \kappa_0 +\kappa_1\,(x-m_0)+\kappa_2\,(x-m_0)^2}\) where \( \kappa_i\ge 0,\;i=0,2, \;\kappa^2_1 \le 4\kappa_2\kappa_0. \)
	\begin{figure}[H]
		\centering
		\includegraphics[width=0.9\linewidth]{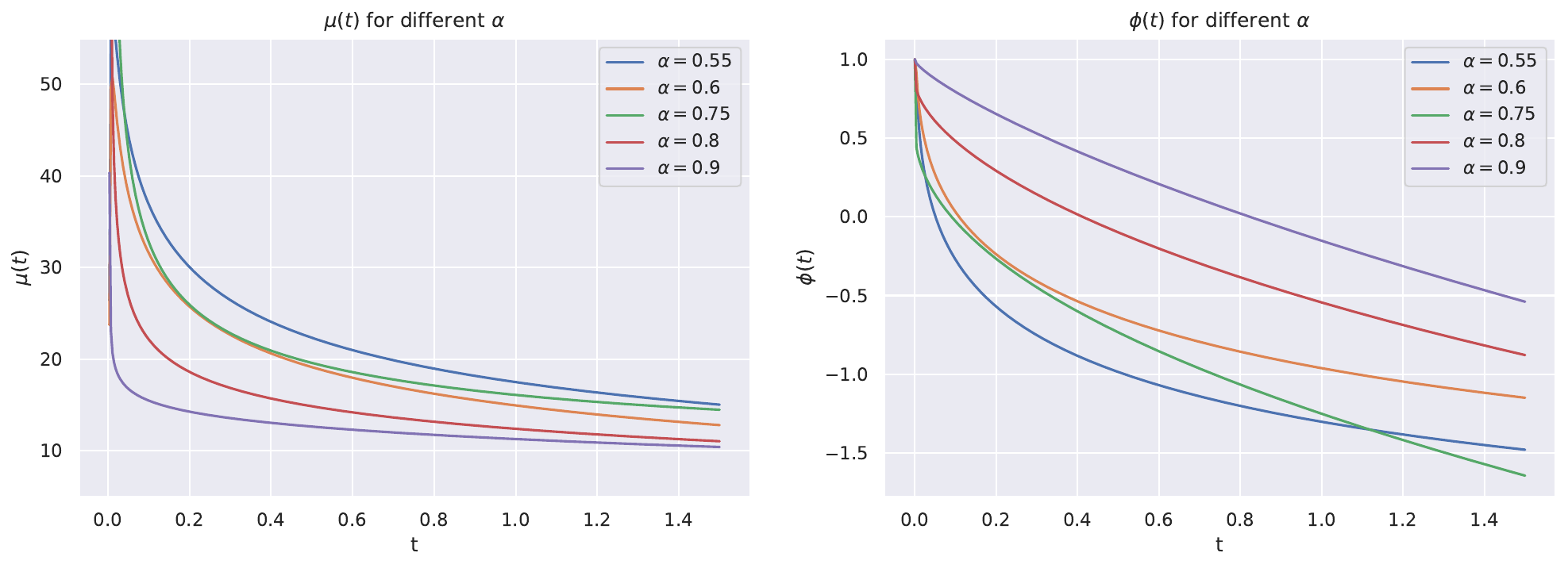}
		\caption{\textit{Graph of \( t \mapsto \mu(t) \) (left) and \( t \mapsto \phi(t) \) (right) for the fractional kernel (Example~\ref{Ex:FakeStatAutonome}~$-3.$) with different values of \(\alpha\) over the time interval \( [0, T] \), \( T = 1.5 \), \( \lambda = 0.2 \), \( m_0 = 10 \) and \(c=1.95356\).}}\label{fig:Mu_phi_autonome}
	\end{figure}
	\begin{figure}[H]
		\centering
		\includegraphics[width=0.9\linewidth]{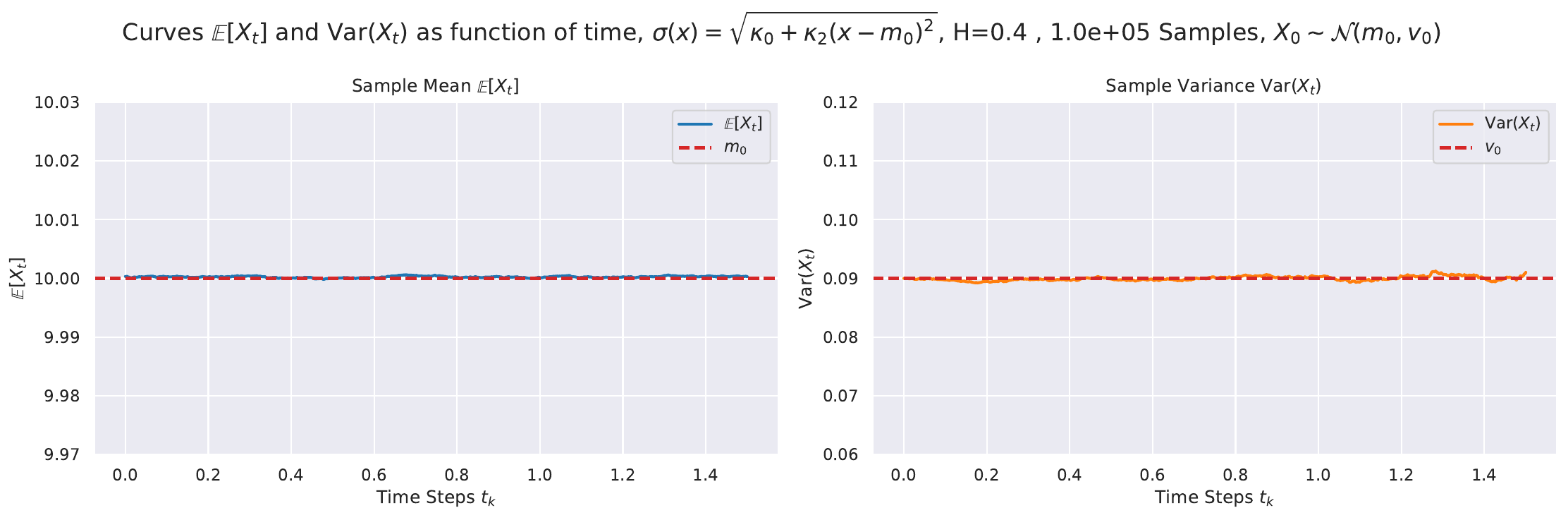}
		\caption{\textit{Graph of \( t_k \mapsto \mathbb{E}[X_{t_k}] \) and \( t_k \mapsto \text{Var}(X_{t_k}, M) \) for the process~\eqref{eq:ConstMeanVar} with fractional kernel (Example~\ref{Ex:FakeStatAutonome}~$-2.$) over the time interval \( [0, T] \), \( T = 1.5 \), \( H = 0.4 \), \( \lambda = 0.2 \), \( m_0 = 10 \), \( v_0 = 0.09 \), \(c=1.95356\), \(\kappa_0 =0.01151\) and \(\kappa_2 =0.384\). Number of steps: \( n = 800 \), Simulation size: \( M = 160000 \).}}\label{fig:autonome}
	\end{figure}
	The Figure~\ref{fig:autonome} illustrates that the Volterra process $(X_t)_{t\in[0,T]}$ in Equation~\eqref{eq:ConstMeanVar} is designed to ensure constant marginal
	mean and variance over time.
	\subsubsection{A numerical illustration of fake stationary SVIE with exponential $\alpha-$fractional kernels for $\alpha \in (\frac12, 1)$ and (stabilized) quadratic  diffusion coefficient}\label{subsubsect:Num_gammaKernel}
	We consider an $\alpha$-Gamma fractional kernel for $ \alpha \in (\frac12, 1)$ (``rough models ")  and a squared trinomial diffusion coefficient \(\sigma(t,x) = \varsigma(t)\,\sigma(x)\) where
	\(\sigma\) of the form \ref{eq:sigmafakeI&II} and given by \(	\sigma(x) = \sqrt{ \kappa_0 +\kappa_1\,(x-\frac{\mu_0}{\lambda})+\kappa_2\,(x-\frac{\mu_0}{\lambda})^2}\quad \mbox{ with }\quad  \kappa_i\ge 0,\;i=0,2, \;\kappa^2_1 \le 4\kappa_2\kappa_0.\) Let c such that $c[\sigma]^2_{Lip} < 1$.
	For the numerical illustrations, we consider the case where for every \(t>0\), \(\mu(t)=\mu_0=\lambda m_0\), so that $\phi(t) = 1, \, \forall t \geq 0$, in which case the equation with constant mean~\eqref{eq:ConstMean} reads:  
	\begin{equation}\label{eq:GammaSVIE}
		X_t= \frac{\mu_0}{\lambda} + \Big(X_0-\frac{\mu_0}{\lambda}\Big) R_{\alpha,\rho, \lambda}(t)+  \frac{1}{\lambda}\int_0^t f_{\alpha,\rho, \lambda}(t-s)\varsigma(s)\sigma(X_{s})dW_s.
	\end{equation}
	
	\medskip
	\noindent From now on, we restrict ourselves to the case where \(\mu(t)=\mu_0=\lambda m_0\) for every \(t>0\), so that \(\phi \equiv 1\).
	It follows from~\cite[Proposition 5.3]{Pages2024} that, in the case of the fractional kernel, the associated functional equations~\eqref{eq:VolterraStabilizer} admit a unique non-negative concave solution. The purpose of the following analysis is to establish an analogous result for the exponentially damped fractional kernel.
	\begin{Proposition}[Existence, positivity and Properties of the function $\varsigma_{\alpha,\rho,\lambda,c}^2$ for $\alpha \in  (\frac12,1)$]\label{prop:alphaFractKernel1_}
		\begin{enumerate}
			\item In reference to the remark on the stabilizer, consider the following equation for \(c,\lambda >0\):
			\begin{equation}\label{eq:stabil2}
				c \lambda^2 \left( 1 - R_{\alpha,\rho, \lambda}^2(t) \right) = (f_{\alpha,\rho, \lambda}^2 * g_{\alpha,\rho, \lambda})(t), \quad \forall t \geq 0.
			\end{equation}
			where \( R_{\alpha,\rho,\lambda} : \mathbb{R}_+ \to \mathbb{R} \) and \( f_{\alpha,\rho,\lambda} := -R_{\alpha,\rho,\lambda}' \) satisfy \( R_{\alpha,\rho,\lambda}(0) = 1 \), \( \lim_{t\to +\infty}R_{\alpha,\rho,\lambda} (t) = a \), and \( \lim_{t\to +\infty}f_{\alpha,\rho,\lambda} (t) = 0 \).
			\begin{enumerate}
				\item[(a)] Then equation ~\eqref{eq:stabil2} has at most one solution in \( \mathcal{L}^1_{\text{loc}}(\R_+) \) that converges to a finite limit.
				
				\item[(b)] If the equation ~\eqref{eq:stabil2} has a continuous solution \( g_{\alpha,\rho,\lambda} \) defined on \( (0, +\infty) \), then \( g_{\alpha,\rho,\lambda} \geq 0 \) on \( \mathbb{R}_+ \), so that the function \( \sqrt{g_{\alpha,\rho,\lambda}} \) is well-defined on \( \mathbb{R}_+ \). 
				Moreover, \(g_{\alpha,\rho,\lambda} \) is concave, non-decreasing and non-negative on \( [0, +\infty)\). In particular, we have \( \forall t > 0 \quad  g_{\alpha,\rho,\lambda} (t) > 0  \). 
				
				
			\end{enumerate}
			\item The stabilizer \( \varsigma^2_{\alpha,\rho,\lambda,c} \) exists as a non-negative function on \( (0, +\infty) \) and
			
			$ \lim_{t\to0} \varsigma_{\alpha,\rho, \lambda,c} = 0
			\quad \text{and} \quad 
			\lim_{t \to +\infty} \varsigma_{\alpha,\rho,\lambda,c}(t) = \frac{\sqrt{c(1-a^2)}\lambda}{\|f_{\alpha,\rho,\lambda}\|_{\mathcal{L}^2(\R_+)}}, \quad a= \frac{1}{1 + \lambda \rho^{-\alpha}} .$ 
		\end{enumerate} 
	\end{Proposition}
	\noindent {\bf Proof.} Claim 1(a) comes from Lemma \ref{lem:asymptotique} (3). To prove 1(b), we recall that the class of completely monotone (CM) functions is a convex cone, thus is stable under pointwise positive 
	summation, product, and also convolution. Differentiating both sides of equation~\eqref{eq:stabil2} and using the 
	fact that \( g_{\alpha,\rho,\lambda}(0) = 0 \) (see e.g. Equation~\eqref{eq:varsigma2}) yields
	\begin{equation}\label{eq:Deriv_stabil}
	2c \lambda^2 f_{\alpha,\rho,\lambda}(t)R_{\alpha,\rho,\lambda}(t)
	= \int_{0}^{t} f_{\alpha,\rho,\lambda}^2(t-s)\, g'_{\alpha,\rho,\lambda}(s)\,ds,
	\quad \forall\, t \geq 0.
	\vspace{-.1cm}
	\end{equation}
	Note that \( f_{\alpha,\rho,\lambda}(t) \) is CM by Proposition~\ref{prop:main2}.
	 Hence \( R_{\alpha,\rho,\lambda} \) is CM, and consequently 
	both \( 2c \lambda^2 f_{\alpha,\rho,\lambda}R_{\alpha,\rho,\lambda} \) and 
	\( f_{\alpha,\rho,\lambda}^2 \) are CM functions.  
	Since \( g_{\alpha,\rho,\lambda}(0) = 0 \), we deduce from~\cite[Theorem V.5.5]{gripenberg1990} 
	that \(g_{\alpha,\rho,\lambda}(t) = \int_{0}^{t} g'_{\alpha,\rho,\lambda}(s)\,ds = 2c \lambda^2 \int_{0}^{t} f_{\alpha,\rho,\lambda}(t-s)R_{\alpha,\rho,\lambda}(t-s) r_{f_{\alpha,\rho,\lambda}^2}(ds) \geq 0,
	\; \forall\, t \geq 0.\) where $r_{f_{\alpha,\rho,\lambda}^2}$ is the functional resolvent of the first kind of $f_{\alpha,\rho,\lambda}^2$ defined as in Equation~\eqref{eq:resolvent_first_kind}. Therefore, \(g_{\alpha,\rho,\lambda}\) is non-negative on \([0,+\infty)\), so that the function \( \sqrt{g_{\alpha,\rho,\lambda}} \)  is well-defined on \( [0, +\infty) \).
	Claim \((2)\) follows from 1(b), equation~\eqref{eq:varsigma2} and Lemma~\ref{lem:asymptotique} (4). The proof of 1(b) (concavity and non-decreasing) is postponed to Appendix \ref{app:lemmata}. \hfill $\square$

	In this exponential-fractional case, the stabilizing function $\varsigma_{\alpha,\rho, \lambda,c}$, solution of the associated functional equations~\eqref{eq:VolterraStabilizer} is computed using the scheme~\eqref{eq:TriangStabilScheme}.
	on a time grid \(t_k = k\frac{T}{n}, k = 0, . . . , n.\), using a piecewise linear interpolation on \([t_{j},t_{j+1}]\), we build the recursive
	scheme: \(\varsigma_{\alpha,\rho, \lambda,c}^2(0)=0\) and \(\forall\, k\ge 1,\)
	\begin{align}\label{eq:TriangStabilScheme}
		c \lambda^2\Big(1&- R_{\alpha,\rho, \lambda}^2 (t_k) \Big) 
		= \sum_{j=0}^{k-1} \sum_{\ell=0}^{1} \varsigma_{\alpha,\rho, \lambda,c}^2(t_{j+\ell}) \int_{t_j}^{t_{j+1}} f_{\alpha,\rho, \lambda}^2(t_k-s)(-1)^\ell \frac{t_{j+1-\ell}-s}{t_{j+1}-t_j}\, ds
	\vspace{-.2cm}
	\end{align}
	which we can solve step by step (Lower-Triangular system) to recover the values \(\varsigma_{\alpha,\rho, \lambda,c}^2(t_k), \, k\ge 1.\)
	This is the standard second-order extension of the triangular Volterra discretization.
		\begin{figure}[H]
	\centering
	\includegraphics[width=.82\linewidth]{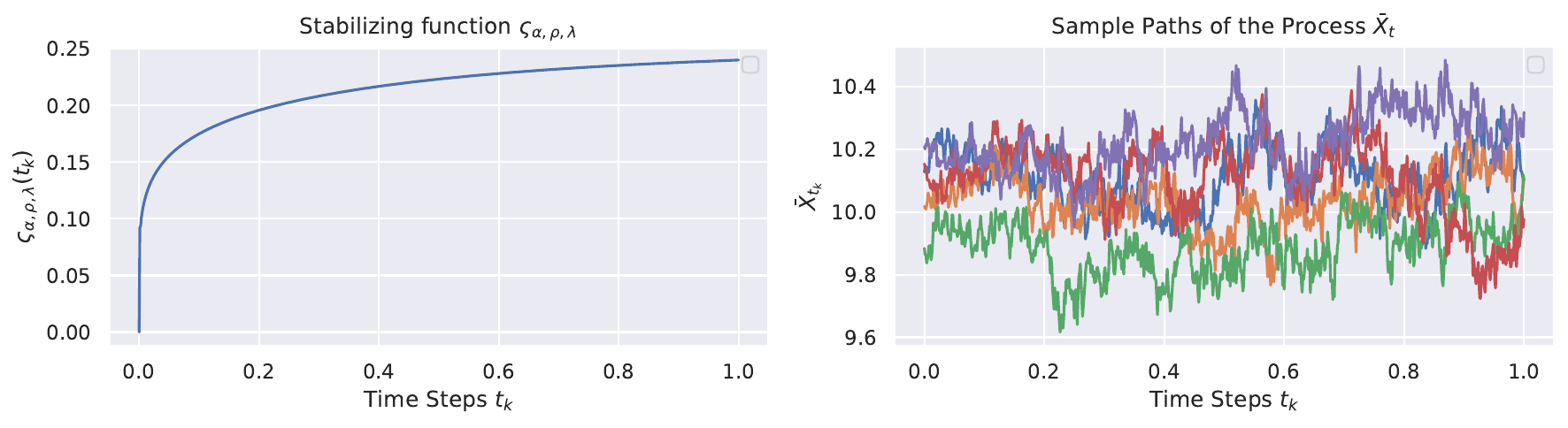}
	\caption{\textit{
			Stabilizer and sample paths of $(X_t)_{t\in[0,1] }$ with parameters: \(\mu(t)=\lambda m_0=\mu_0\), $\sigma(x) = \sqrt{\kappa_0 + \kappa_2 (x - \frac{\mu_0}{\lambda})^2}$, $H=0.1$, $c=0.316$, $\lambda=0.2$, $\rho=1.2$, $\mu_0=2$, $v_0=0.09$, $\kappa_0=0.25$ and $\kappa_2=0.384$.
	}} 
	\label{fig:illustrations_H01}
\end{figure}
	\begin{figure}[H]
		\centering
		\includegraphics[width=0.83\linewidth]{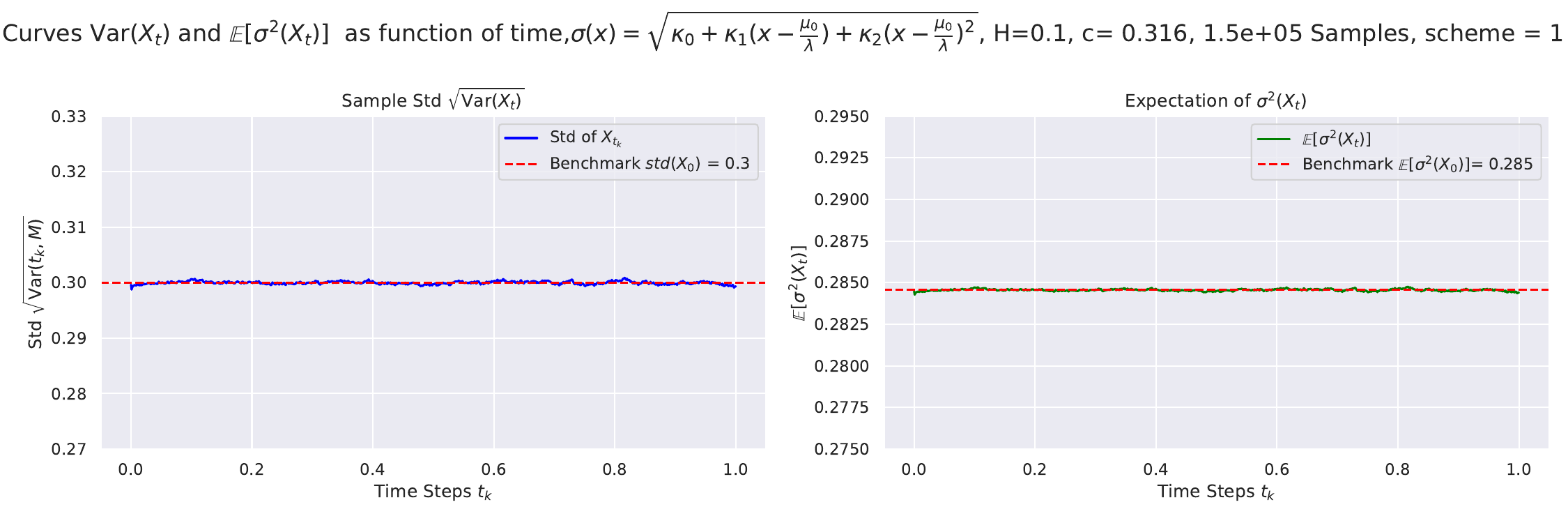}
		\caption{Graph of \( t_k \mapsto \text{Std}(t_k, M) \) and \( t_k \mapsto \mathbb{E}[\sigma^2(X_{t_k},M)] \) over the time interval \( [0, 1] \),  \( H = 0.1 \), \( \mu_0 = 2 \), \( \lambda = 0.2 \), \( v_0 = 0.09 \), \( \rho = 1.2 \), \( \text{Std}(X_0) = 0.3 \). Number of steps: \( n = 800 \), Size: \( M = 1.5\times 10^5 \).}\label{fig:inhomStabilizer1}
		\vspace{-0.5cm}
	\end{figure}
		\begin{figure}[H]
		\centering
		\includegraphics[width=0.83\linewidth]{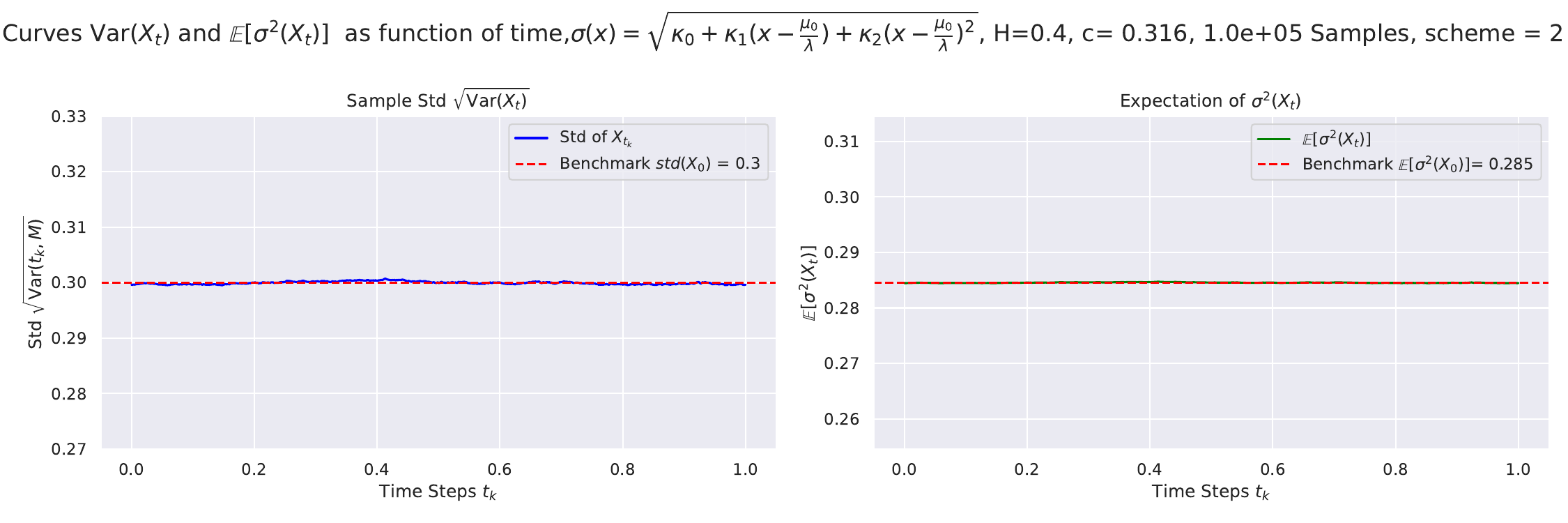}
		\caption{Graph of \( t_k \mapsto \text{Std}(t_k, M) \) and \( t_k \mapsto \mathbb{E}[\sigma^2(X_{t_k},M)] \) over the time interval \( [0, 1] \), \( H = 0.4 \), \( \mu_0 = 2 \), \( \lambda = 0.2 \), \( v_0 = 0.09 \), \( \rho = 1.2 \), \( \text{Std}(X_0) = 0.3 \). Number of steps: \( n = 800 \), Samples size: \( M = 10^5 \).}
		\label{fig:inhomStabilizer2}
		\vspace{-0.5cm}
	\end{figure}
\noindent The stabilizing functions $\varsigma$ introduced in Theorem~\ref{thm:timeautonome} and Equation~\eqref{eq:VolterraStabilizer} are designed to guarantee that the Volterra processes $(X_t)_{t\geq0}$ in~\eqref{eq:ConstMean} exhibit invariance of their first two moments under time shifts (see Figures~\ref{fig:inhomStabilizer1} and~\ref{fig:inhomStabilizer2}).
	
	\vspace{-0.3cm} 
	\subsubsection{A numerical illustration of the degenerate case of fake stationarity in SVIE with $\alpha$-fractional kernels for $\alpha \in (\frac12,1)$ and a (stabilized) tanh diffusion coefficient}
	\noindent
	In this section, we work in the framework of Section~\eqref{subsubsect:Num_gammaKernel}, but consider a family of scaled volterra equations with $\alpha$-fractional kernel (\(\rho=0\)) for $ \alpha \in (\frac12, 1)$ and a diffusion coefficient given by \(\sigma(x) = \sqrt{\frac{\tanh(x-\frac{\mu_0}{\lambda})}{2}}, \lambda > 0 \). In this fractional case, the stabilizing function $\varsigma_{\alpha, \lambda,c}$, solution of the associated functional equations~\eqref{eq:VolterraStabilizer} is computed using the power series expansion provided in~\cite[Equations (5.45) and (5.49)]{Pages2024}.
	\vspace{-0.4cm}
	\begin{figure}[H]
			\vspace{-0.3cm}
		\centering
		\begin{minipage}{0.48\linewidth}
			\centering
			\includegraphics[width=1\linewidth]{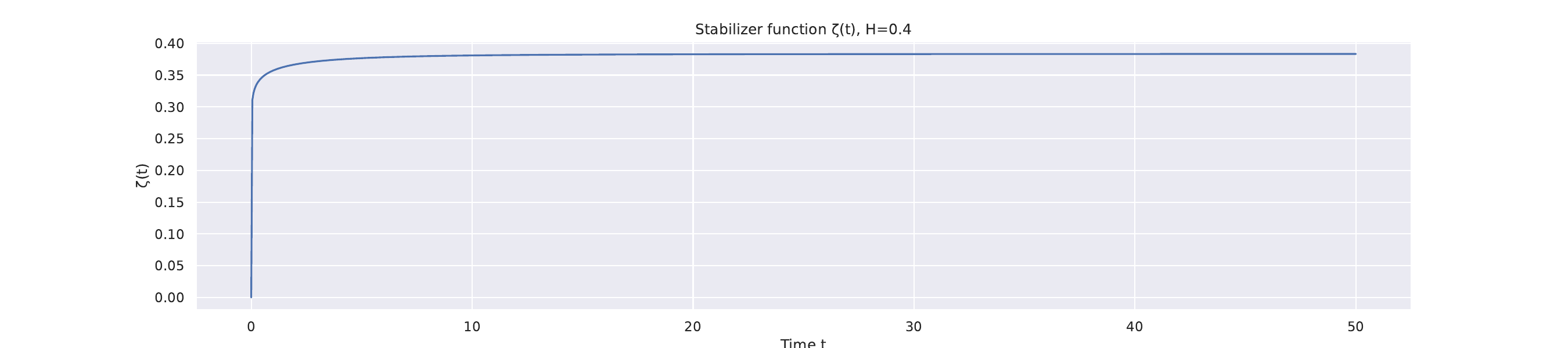}
			\caption{ Graph of the stabilizer $ t \to \varsigma_{\alpha,\lambda,c}(t)$  over time interval [0, T ], T = 50 for a value of the Hurst esponent $H=0.4$,  $\lambda = 0.2$, $c = 0.36$.}
		\end{minipage}%
		\hfill
		\begin{minipage}{0.49\linewidth}
			\centering
			\includegraphics[width=0.71\linewidth]{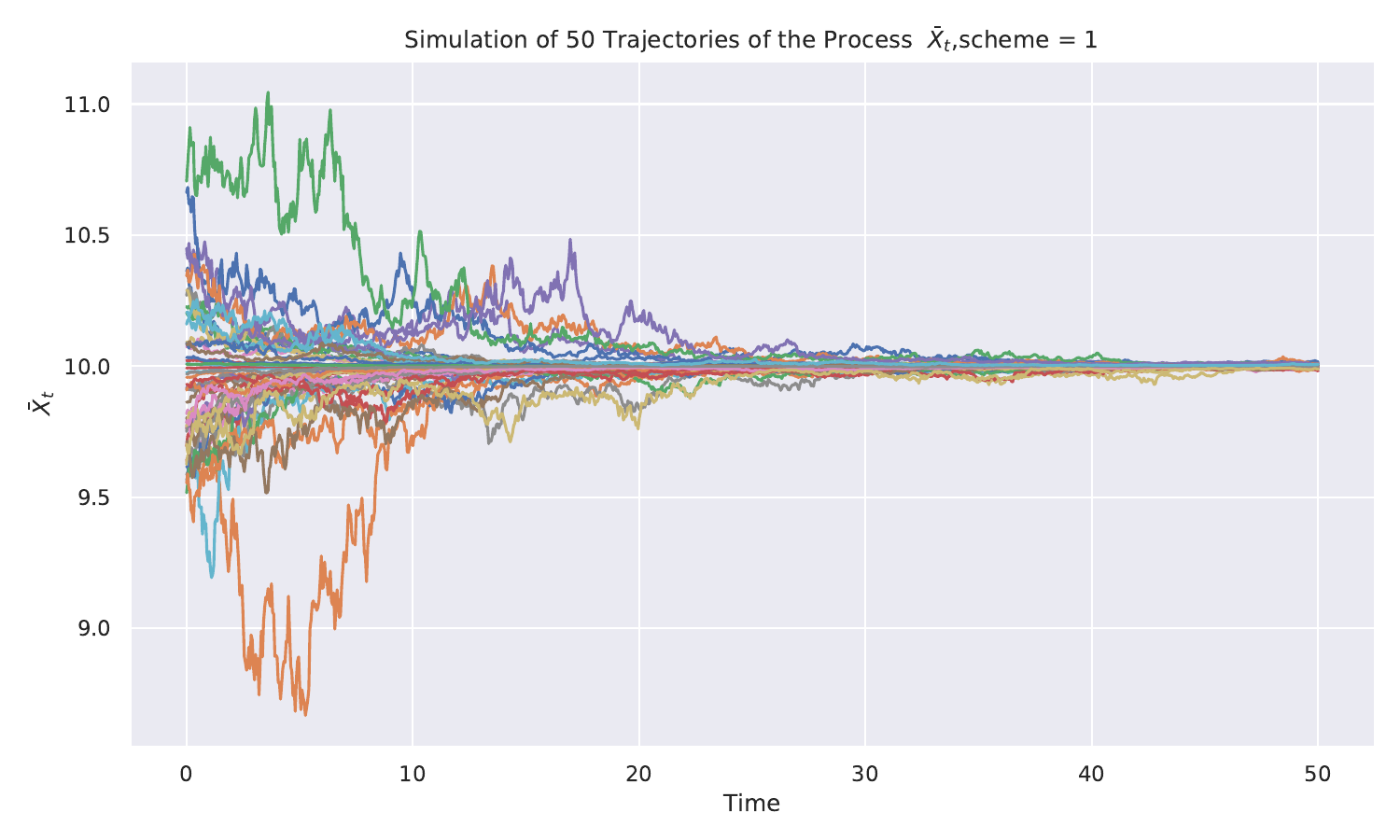}
			\caption{Confluent trajectories in the degenerate case, $T=50$, $H=0.4$,  $\lambda = 0.2$, $c = 0.36$.}\label{fig:confluent2}
		\end{minipage}
		\vspace{-0.5cm}  
	\end{figure}
	Figure~\eqref{fig:confluent2} show that that the marginals of such a process when starting with various initial values
	are confluent in $L^2$ as time goes to infinity.
	
	\smallskip	
	\noindent {\bf Acknowledgement:}  The authors thank J-F. Chassagneux and M. Rosenbaum for insightful discussions, help and comments.
	\vspace{-.5cm}
			\bibliographystyle{abbrv}
			\bibliography{StationarityVolterraEquations}
			
			\normalsize  
			\appendix
			\vspace{-.3cm}
			\section{About the Simulation of the semi-integrated scheme for stochastic Volterra integral Equations (SVIE) }\label{app:simu_det}
			We aim at providing numerical approximation for the equation~\eqref{eq:ConstMeanVar}-\eqref{eq:ConstMean} reading more generally:  
			\vspace{-.2cm}
			\begin{equation}\label{eq:VolterraSimul}
				X_t= g_0(t) +  \frac{1}{\lambda}\int_0^t f_{ \lambda}(t-s)\varsigma(s)\sigma(X_{s})dW_s,\;\; g_0(t) = X_0 - \frac{\big(X_0 - m_0\big)}{\lambda m_0} \int_0^t f_{\lambda}(t-s) \mu(s)\dd s. 
				\vspace{-.1cm}
			\end{equation}
			We introduce the following Euler-Maruyama scheme~\eqref{eq:Volterragen2} on the time grid $t_k =t^n_k =\frac{kT}{n}, k=0, \dots, n$, writing recursively:
			\vspace{-.3cm}
			{\small
				\begin{align}\label{eq:Volterragen2} 
					\overline X^n_{t_{k}^n} = g(t_k^n) + \frac{1}{\lambda}\sum_{\ell=1}^{k} \int_{t_{\ell-1}}^{t_{\ell}} f_{\lambda}(t_k^n-s) \,\varsigma(t_{\ell}^n) \sigma( \overline{X}^n_{t_{\ell-1}})dW_s= g(t_k^n) + \frac{1}{\lambda}\sum_{\ell=1}^{k} \varsigma(t_{\ell}^n) \sigma( \overline{X}^n_{t^n_{\ell-1}}) I^{n,\ell}_k
				\end{align}
			} 
			where the stochastic integrals $I^{n,\ell}_k:= \int_{t_{\ell-1}}^{t_{\ell}} f_{\lambda}(t_k^n-s) dW_s$ are defined on the time grid $t_k =t^n_k =\frac{kT}{n}, k=0, \dots, n$ and \(g_0(t_k)\) can be computed using
			quadrature formulae on different intervals ( Gauss-Legendre, Gauss-Laguerre etc.).
			Due to the lack of Markovianity, $\bar X_{t^n_k}$ is generally not a function of $(\bar X_{t^n_{k-1}}, W_{t^n_k} - W_{t^n_{k-1}})$. However, it can be computed uniquely from $(\bar X_{t^n_0}, \dots, \bar X_{t^n_{k-1}})$ and the Gaussian vector $I^{n}_k:=\Big(\int_{t_{\ell-1}}^{t_{\ell}} f_\lambda(t^n_k-s) dW_s \Big)_{\ell=1,\dots,k}$, ensuring that the Euler-Maruyama scheme is well-defined by induction.
			The exact simulation of the Euler-Maruyama scheme~\eqref{eq:Volterragen2} on the discrete  grid \((t^n_k)_{0\leq k\leq n}\) involves simulating the independent random vectors:$\left(\int_{t_{\ell-1}}^{t_{\ell}} f_{\lambda}(t_k^n-s) dW_s \right)_{\ell \leq k \leq n}, \quad \ell = 1, \dots, n.$\\
			\noindent {\bf Practitioner's corner:}.
			We aim at simulating all the $I^{n}_k:=\big(I^{n,\ell}_k \big)_{\ell=1,\dots,k}, $ together with the increment $\Delta W_{t_k}$ at once on the discrete  grid \((t_k)_{0\leq k\leq n}\). We will rather consider and generate the $n$ independent sequence of Gaussian vectors \( G^{n,\ell}, \ell=1 \cdots n\) using the Cholesky decomposition of their covariance matrix:
			\vspace{-.2cm}
			{\small
				\begin{equation}\label{eq:Gaussian_G_n}
					G^{n,\ell} = \Big(\Delta W_{t_\ell}, \big[\int_{t_{\ell-1}}^{t_{\ell}} f_{\lambda}(t_k-s) dW_s \big]_{k=\ell,\ldots,n} \Big)
					=\Big(\Delta W_{t_\ell}, \big[I^{n,\ell}_k \big]_{k=\ell,\ldots,n} \Big),\quad \ell= 1,\ldots,n.
					\vspace{-.2cm}
				\end{equation}
			}
			Let us denote by $G^n= (G_{k\ell}^{n})_{ k=1:n, \ell=1: n+1}$ the  $n\times (n+1)$ matrix  involving the random terms \( G^{n,\ell}, \ell=1 \cdots n\).
			\vspace{-.2cm}
			{\small
				\begin{equation*}
					G^n:=  \begin{array}{c@{}c}
						& \begin{array}{cccccc}
							\hspace{1.5cm}I^{n}_1 \hspace{2cm}& I^{n}_2\hspace{3cm} & \dots & I^{n}_n \\ \\
						\end{array} \\[-0.5em]
						&
						\begin{bmatrix}
							\Delta W_{t_1} 
							& \int_{t_{0}}^{t_1} f_\lambda(t_1 - s) \, dW_s
							& \int_{t_{0}}^{t_1} f_\lambda(t_2 - s) \, dW_s
							& \cdots 
							& \int_{t_{0}}^{t_1} f_\lambda(t_n-s) \, dW_s \\[0.3em]
							0 
							& \Delta W_{t_2} 
							& \int_{t_{1}}^{t_2} f_\lambda(t_2 - s) \, dW_s
							& \cdots 
							& \int_{t_{1}}^{t_2} f_\lambda(t_n-s) \, dW_s \\[0.3em]
							0 
							& 0 
							& \Delta W_{t_3} 
							& \cdots 
							& \int_{t_{2}}^{t_3} f_\lambda(t_n-s) \, dW_s \\[0.3em]
							\vdots 
							& \vdots 
							& \vdots 
							& \ddots 
							& \vdots \\[0.3em]
							0 
							& 0 
							& \cdots 
							& \Delta W_{t_n}
							& \int_{t_{n-1}}^{t_n} f_\lambda(t_n-s) \, dW_s 
						\end{bmatrix}
					\end{array}
					\vspace{-.2cm}
				\end{equation*}
			}
			\smallskip \noindent
			Then the following relation holds: $\;\overline X^{n}_{0}=g_{0}(t_0)$ and for every $k=1,\dots,n$
			\vspace{-.2cm}
			\begin{equation}\label{eq:VarianceModel2}
				\big( \overline X^{n}_{t_{k}} \big)_{k=1:n}= \big(g_0(t_{k})\big)_{k=1:n} +G^n_{\cdot,k+1}\big(\varsigma(t_{\ell})\sigma( \overline{X}^{n}_{t_{\ell-1}})\mbox{\bf 1}_{\{\ell \le k\}}\big)_{\ell=1:n}
			\end{equation}
			\noindent {\bf Remark:}
			Note that we consider the Brownian increment in the vector~\eqref{eq:Gaussian_G_n} because, in applications to volatility model dynamics, the dynamics of the traded asset and its volatility process can be jointly driven by the same Brownian motion (see for exemple the quadratic rough volatility dynamic introduced in~\cite{GaJuRo2020} or also ~\cite{Gnabeyeu2026b,Gnabeyeu2026a, dro2026optimal}). This approach takes into account, among other factors, the so-called Zumbach effect, which links the evolution of the asset or an index with its volatility.
			
			\smallskip
			\noindent For $ \ell= 1,\ldots,n$, the covariance matrix of $\big(\big[I^{n,\ell}_k \big]_{k=\ell,\ldots,n} \big)$ is symmetric and of size $(n-\ell+1) \times (n-\ell+1)$, given (where we use the change of variables $u= \frac Tn (\ell-v)$, $v\!\in [0,1]$) by:
			{\small
				\begin{align*}
					\Sigma^{n,\ell} & = \Big[ Cov(I^{n,\ell}_k, I^{n,\ell}_{k'}) \Big]_{\ell\le k,k'\le n} =  \Big[ \int_{t_{\ell-1}}^{t_\ell} f_{\lambda}(t_k - u) f_{\lambda}(t_{k'}-u) du \Big]_{\ell\le k,k'\le n} , \\
					& = \Big(\frac Tn\Big)\Big[ \int_{0}^{1} f_{\lambda}(\frac Tn (k-\ell+v) ) f_{\lambda}(\frac Tn(k'-\ell+v)) dv \Big]_{\ell\le k,k'\le n} = \Big(\frac Tn\Big) \big[\Omega^{n}_{k-\ell, k'-\ell}\big]_{\ell\le k,k'\le n},  
				\end{align*}
			}
			where the  symmetric matrix $\Omega^{n}$ is defined by \(\Omega^{n} := \Big[\int_0^1 f_{\lambda}(\frac Tn(i+v))f_{\lambda}(\frac Tn(j+v)) dv \Big]_{0\le i,j\le n-1}\). We also set \(\Omega^{0} := \Big[\int_0^1 f_{\lambda}(\frac Tn(i+v))dv \Big]_{0\le i\le n-1}\) and define for $ \ell= 1,\ldots,n$,
			{\small
				\begin{align*}
					\Sigma^{0,\ell}  
					&= \Big[\mathrm{Cov}\Big( \Delta W_{t_\ell}\text{,\ } I^{n,\ell}_{k} \Big)\Big]_{\ell \le k \le n} 
					= \Big[ \int_{t_{\ell-1}}^{t_\ell} f_{\lambda}(t_k - u) \, \mathrm{d}u \Big]_{\ell \le k \le n} \\
					&= \Big[ R_{\lambda}(t_k - t_\ell) - R_{\lambda}(t_k - t_{\ell-1}) \Big]_{\ell \le k \le n} 
					= \Big(\frac{T}{n}\Big) \Big[ \int_0^1 f_{\lambda}\Big(\frac{T}{n}(k-\ell+v)\Big) \, \mathrm{d}v \Big]_{\ell \le k \le n} =: \Big(\frac{T}{n}\Big) \Big[ \Omega^{0}_{k-\ell} \Big]_{\ell \le k \le n}.
				\end{align*}
			}
			Therefore, the covariance matrix of the Gaussian vector $G^{n,\ell}$ is a symmetric $(n-\ell+2)\times(n-\ell+2)$ sub-matrix of the $(n+1)\times(n+1)$ covariance matrix \(C\), given by
			{\small
				\[C =
				\begin{pmatrix}
					\frac{T}{n} & \begin{array}{ccc} ^t \Sigma^{0,1} \end{array} \\[10pt]
					\begin{array}{c} \Sigma^{0,1} \end{array} & \Sigma^{n,1}
				\end{pmatrix} = \Big(\frac Tn\Big)
				\begin{pmatrix}
					1 & \begin{array}{ccc} \Omega^0 \end{array} \\[10pt]
					\begin{array}{c} \Omega^0 \end{array} & \Omega^n
				\end{pmatrix}.
				\]
			}
			At this stage, we can compute any fixed sub-matrix  of $C$ by a cubature formula (such as Trapezoid, Midpoint, Simpson, higher-order Newton-Cote, or Gauss-Legendre integration formulas) and then perform a numerically stable extended Cholesky transform. This results in the decomposition:
			\centerline{$ [C_{ij}]_{1\le i,j\le n+1}= T^{(n)}D^{(n)}T^{(n)*} 
				\quad T^{(n)} \text{ lower triangular}.
				$}
			
			\noindent $T^{(n)}$ is a lower triangular matrix with diagonal entries $T^{(n)}_{ii} = 1$, and $D^{(n)}$ is a diagonal matrix with non-negative entries. 
			Then, taking advantage of the telescopic feature and the structure of this Cholesky transform one has: 
			\vspace{-.2cm}
			\begin{equation}
				[C_{ij}]_{1\le i,j\le n+2-\ell}=  [T_{ij}^{(n)}]_{1\le i,j\le n+2-\ell}[D^{(n)}_{ij}]_{1\le i,j\le n+2-\ell} [T_{ij}^{(n)}]^*_{1\le i,j\le n+2-\ell},\; \ell= 1,\ldots,n.
				\vspace{-.1cm}
			\end{equation}
			Finally, for each $\ell = 1, \ldots, n$, we have: $
			(G^{n,\ell})_{\ell=1, \dots, n} \stackrel{d}{=} (\tilde{T}^{(n+2-\ell)} Z^{(\ell)})_{\ell=1, \dots, n},$ where 
			\vspace{-.1cm}
			\begin{equation}Z^{(\ell)} \sim \mathcal{N}(0, I_{n-\ell+2}) \quad \text{and} \quad \tilde{T}^{(n+2-\ell)} = [T_{ij}^{(n)}]_{1 \leq i,j \leq n +2 -\ell} [\sqrt{D^{(n)}_{ij}}]_{1 \leq i,j \leq n +2 -\ell}.
				\vspace{-.2cm}
			\end{equation}
			\noindent {\bf Remarks:} 1. This Cholesky matrix is usually quite sparse ({\em when $H$ is small} in the case of fractional kernel for example) since, all entries beyond the fourth column are numerically $0$ (in fact smaller than $10^{-4}$). This is due to the fact that such singular kernels have essentially no memory for small $H$. This feature quickly disappears when running the procedure with $H>1/2$.
			
			2. If the fonction $f_{\lambda}$ is a monone (case where we replace it mutantis mutandis by the fractional integration kernel $K (u)=K_{1,\alpha,0}(u) = \frac{u^{\alpha-1}}{\Gamma(\alpha)}, u \in [0,T]$, where $\alpha \!\in (1/2,1)$), we will have the fact that $C^{n}$ is a certain power factor of $\Big(\frac Tn\Big)$, say $\Psi \Big(\frac Tn\Big)^\beta$,  times an infinite  symmetric matrix ($\bar\Gamma$) (not depending on n anyway) defined by \(\bar\Gamma := \Big[\int_0^1 \big((i+v)(j+v)\big)^{(\alpha-1)} dv \Big]_{i,j\ge 0}.\)
			In this case, the diagonal entries of $\bar\Gamma$ have  closed form formular and the matrices of interest  $[C_{ij}]_{0\le, i,j\le n-1}$, $n\ge 1$  are telescopic  sub-matrices of $\bar\Gamma$ times the factor $\Psi\Big(\frac Tn\Big)^\beta$.
			
			3. For comprehensive results concerning convergence rates and \emph{a priori} error estimates related to the approximation of the stochastic Volterra process~\eqref{eq:VolterraSimul} by the semi-integrated Euler--Maruyama scheme~\eqref{eq:Volterragen2} , as well as its continuous-time (or ``genuine") extension, the reader is referred to~\cite{JouPag22}.
			
					\vspace{-.1cm}
			\section{Supplementary material and Proofs. }\label{app:lemmata}
			\noindent {\bf Proof of Proposition \ref{prop:W-H}:} 
			Convoluting $x(t) + \lambda\int_0^t K(t-s) x(s) ds$ with $R'_{\lambda}$ together with the fact that $\lambda K*R'_{\lambda} = -\lambda K -R'_{\lambda}$(see equation ~\ref{eq:flambda-eq} ), we obtain:
			
			\vspace{-.5cm}
			\begin{align*}
				\int_0^t g(s) R'_{\lambda}(t-s)\,& ds = \int_0^t x(s) R'_{\lambda}(t-s) ds + \lambda\int_0^t \int_0^s K(s-u) x(u) du R'_{\lambda}(t-s) ds\\ 
				&= \int_0^t x(s) R'_{\lambda}(t-s) ds + \lambda\int_0^t \int_0^{t-u} K(t-u-s) R'_{\lambda}(s)  ds x(u) du \\
				&= \int_0^t x(s) R'_{\lambda}(t-s) ds - \int_0^t \big(\lambda K(t-u)+R'_{\lambda}(t-u)\big) x(u) du  
				= - \lambda \int_0^t K(t-s) x(s) ds.
			\end{align*}
			Inserting this in the Wiener-Hopf equation gives the results.
			For the second claim, we can use the Laplace transform in the integral equation and deduced that:
			
			\centerline{$
				L_x(t) = \frac{L_h(t)}{1 + L_{R'_{\lambda}}(t)} = L_h(t) \big(1 + \lambda L_K(t)\big) = L_{h + \lambda K * h}(t)
				$}
			\noindent
			where the penultimate equality comes from applying the Laplace transform to $R'_{\lambda} * K = - K - \frac{R'_{\lambda}}{\lambda}$ (see Equation~\eqref{eq:flambda-eq}). We then conclude by the injectivity of the Laplace transform. \hfill $\Box$
			
						\bigskip
			\noindent {\bf Proof of Lemma \ref{lem:on_mu}.} 
			To prove the second equality, for our convenience, we will consider two cases:
			
			\smallskip
			{\em ($f_{\lambda}$ is a positive density).} If \(f_{\lambda} > 0 \) on \((0,+\infty)\) (i.e. \(R_{\lambda} \) decreases), then the function \(f_{\lambda} \) is a positive \((1-a)-\) sum density.
			{\small
				\begin{align*} 
					\int_0^t f_{\lambda}(t-s) \mu (s)ds - \mu_{\infty} (1-a)
					&= \int_0^t f_{\lambda}(t-s) \mu (s)ds - \mu_{\infty}\int_0^{+\infty} f_{\lambda}(s) \, ds \\
					&= \int_0^t f_{\lambda}(t-s)(\mu(s)-\mu_\infty)\,ds - \mu_\infty \int_t^{+\infty} f_{\lambda}(s)\,ds 
				\end{align*} 
			}
			so that, by the triangle inequality, we have:
			\vspace{-.3cm}
			\begin{equation}
				\Big|\int_0^t f_{\lambda}(t-s) \mu (s)ds - \mu_{\infty} (1-a)\Big| \leq \int_0^t f_{\lambda}(t-s)\Big|\mu(s)-\mu_\infty\Big|\,ds + |\mu_\infty| \int_t^{+\infty} f_{\lambda}(s)\,ds
				\vspace{-.2cm}
			\end{equation}
			First note that we can split the first integral as follows:
			{\small
				\begin{align*}
					\int_0^t f_{\lambda}(t-s)\left|\mu(s)-\mu_\infty\right|\,ds &=\int_0^{t-A_{\epsilon}} f_{\lambda}(t-s)\left|\mu(s)-\mu_\infty\right|\,ds
					+ \int_{t-A_{\epsilon}}^t f_{\lambda}(t-s)\left|\mu(s)-\mu_\infty\right|\,ds \\
					&= \int_{A_{\epsilon}}^t f_{\lambda}(s) \left|\mu(t-s)-\mu_\infty\right|\,ds + \int_0^{A_{\epsilon}} f_{\lambda}(s) \left|\mu(t-s)-\mu_\infty\right|\,ds.
				\end{align*}
			}
			where $A_{\epsilon}$ is chosen such that for all $s \ge A_{\epsilon}$, we have $|\mu(s) - \mu_\infty| \le \epsilon$. Moreover $\forall s \in ]0,A_{\epsilon}[, \quad t-s \geq t-A_{\epsilon}\geq A_{\epsilon}$ for t large enough, (say $t \geq 2A_{\epsilon}$), and hence, this implies that $\left|\mu(t-s)-\mu_\infty\right| \leq \epsilon, \quad \forall s \in ]0,A_{\epsilon}[$.
			
			We thus have: $
			\int_0^t f_{\lambda}(t-s)\left|\mu(s)-\mu_\infty\right|\,ds \le \|\mu - \mu_\infty\|_{\sup} \int_{A_{\epsilon}}^{t} f_{\lambda}(s)\,ds + \epsilon \int_{0}^{A_{\epsilon}} f_{\lambda}(s)\,ds.
			$
			And hence,
			{\small
				\begin{align*} 
					\lim_{t \to \infty}\left|\int_0^t f_{\lambda}(t-s) \mu (s)ds - \mu_{\infty} (1-a)\right|
					&\leq \|\mu - \mu_\infty\|_{\infty} \lim_{t \to \infty} \int_{A_{\epsilon}}^{t} f_{\lambda}(s)\,ds + \epsilon \int_{0}^{A_{\epsilon}} f_{\lambda}(s)\,ds\\ 
					&+ |\mu_\infty|  \lim_{t \to \infty} \int_t^{+\infty} f_{\lambda}(s)\,ds \leq \epsilon (1-a)  \leq \epsilon, \quad \text{since} \int_0^{\infty} f_{\lambda}(s)\,ds = 1-a.
				\end{align*}
			}
			\smallskip
			{\sc Case~2} {\em (\( f_{\lambda} \) is just a $(1-a)$-sum measure).} If \( \int_0^{+\infty} f_{\lambda}(s) \, ds = 1-a \) , a more rigorous Approach to prove the above Lemma make used of Laplace Transforms - and Tauberian Final Value Theorem. Let's assume the \( \mathcal{L}^1\)-integrability of \( f_{\lambda} \), i.e., $\int_0^\infty |f_\lambda(s)| ds < \infty$ so that $L_{| f_{\lambda}|}(t) <+\infty$ for every $t>0$: \( f_{\lambda} \)  has subsequently
			a finite Laplace transform defined (at least) on \(\R^+\). 
			
			Since $\lim_{t \to +\infty} \mu(t) = \mu_\infty$, we have by Tauberian Final Value Theorem $\lim_{z \to 0} z L_{\mu}(z)  = \mu_\infty$. As the Laplace transform of the convolution is the product of the Laplace transforms, we have:
			
			\centerline{$\mathcal{L}\left(\int_0^t f_\lambda(t-s) \mu(s) ds\right)  (z) = L_{f_\lambda}(z) L_{\mu}(z)$}
			
			Therefore,  by Tauberian Final Value Theorem if $\lim_{t \to +\infty} \int_0^t f_\lambda(t-s) \mu(s) ds$ exists, then
			$$\lim_{t \to +\infty} \int_0^t f_\lambda(t-s) \mu(s) ds = \lim_{z \to 0} z \mathcal{L}\left(\int_0^t f_\lambda(t-s) \mu(s) ds\right)  (z) = \lim_{z \to 0} z  L_{\mu}(z) L_{f_\lambda}(z) = \mu_\infty \lim_{z \to 0} L_{f_\lambda}(z)$$
			
			However, by our assumption $\int_0^\infty f_\lambda(s) ds = 1 - a$, we have $\lim_{z \to 0} L_{f_\lambda}(z) = \lim_{z \to 0}\int_0^\infty e^{-zs} f_\lambda(s) ds= L_{f_\lambda}(0)= \int_0^\infty f_\lambda(s) ds = 1 - a$.
			This completes the proof. \hfill$\Box$
			
			\bigskip
			\noindent {\bf Proof of Lemma \ref{lem:asymptotique}.} 
			\medskip
			\noindent  {\sc Step~1. }  If \(  \varsigma_{\lambda, c}^1 \) and \(  \varsigma_{\lambda, c}^2 \) are two solutions of the equation \textit{($E_{\lambda, c}$)} in ~\eqref{eq:VolterraStabilizer}, then \( f_{\lambda}^2 * \delta \varsigma_{\lambda, c} = 0 \) in \( \mathcal{L}^1_{\text{loc}}(\mathbb{R}_+) \) with \( \delta \varsigma_{\lambda, c} = \varsigma_{\lambda, c}^1 - \varsigma_{\lambda, c}^2 \). As $L_{f^2_{\lambda}}(t) > 0$ for $t > 0$ (by Assumption~(${(\cal K)}_\lambda$)(ii)), then \( \delta g = 0 \), which implies \( \delta \varsigma_{\lambda, c} = 0 \) in \(\mathcal{L}^1_{\text{loc}}(\mathbb{R}_+) \). Thus, the solution \( \varsigma^2_{\lambda, c} \) of Equation~\eqref{eq:VolterraStabilizer} if any, is unique. 
			
			We would also note that, \( c \) being fixed, \( L_{f^2_{\lambda}}(t) > 0 \) for \( t > 0 \) (by Assumption~(\({\cal K}_\lambda\))(ii)). Then \( L_{\varsigma^2_{\lambda, c}} \) is uniquely determined by \eqref{eq:Laplacesigma}, which in turn implies the uniqueness of \( \varsigma^2_{\lambda, c} \) (at least \( \text{dt} \)-a.e.).
			
			\medskip
			\noindent  {\sc Step~2. } $\varsigma^2_{\lambda, c}$
			is non-negative and has a limit \( l_\infty \in (0, +\infty] \) as \( t \to +\infty \). If \( l_\infty = +\infty \), then for every \( A > 0 \), there exists \( t_A \) such that for \( t \geq t_A \), \( \varsigma^2_{\lambda, c}(t) \geq A \). Hence
			
			\centerline{$
				(f_{\lambda}^2 * \varsigma^2_{\lambda, c})(t) = \int_{0}^{t_A}  f_{\lambda}^2(t - s) \varsigma^2_{\lambda, c}(s) \, ds + \int_{t_A}^{t} f_{\lambda}^2(t - s) \varsigma^2_{\lambda, c}(s) \, ds \geq
				\int_{0}^{t_A}  f_{\lambda}^2(t - s) \varsigma^2_{\lambda, c}(s) \, ds + A \int^t_{t_A} f_{\lambda}^2(t - s) \, ds$}
			i.e.
			\[
			(f_{\lambda}^2 * \varsigma^2_{\lambda, c})(t) \geq
			\int_{0}^{t_A}  f_{\lambda}^2(t - s) \varsigma^2_{\lambda, c}(s) \, ds + A \int_0^{t -t_A} f_{\lambda}^2(s) \, ds.
			\]
			Consequently, as \( (f_{\lambda}^2 \ast \varsigma^2_{\lambda, c})(t) = c\lambda^2(1 - \big(1-\frac{(f_\lambda * \mu)(t)}{\lambda m_0}\big)^2) \to c\lambda^2(1-\ell_\infty^2) \) as \( t \to +\infty \) owing to Lemma \ref{lem:on_mu}, we have: \(c\lambda^2(1-\ell_\infty^2)  = \lim_{t \to +\infty} (f_{\lambda}^2 \ast \varsigma^2_{\lambda, c})(t) \geq A \int_0^{+\infty} f_{\lambda}^2(u) \, du.\)
			As \( f_{\lambda} \in \mathcal{L}^2 (\R_+) \), this yields a contradiction by letting \( A \to \infty \). Hence, \( l_\infty < +\infty \).
			Now, still by the same arguments,  \(\lim_{t \to +\infty} \varsigma^2_{\lambda, c}(t) = l_\infty\in (0, +\infty[ \implies \forall \eta >0, \exists t_\eta \in \R^+\) such that \(\forall t > t_\eta \quad l_\infty- \eta \leq \varsigma^2_{\lambda, c}(t)\leq l_\infty+ \eta\)
			
			On the first hand, we have:
			{\small
				\[
				(f_{\lambda}^2 * \varsigma^2_{\lambda, c})(t) \geq
				\int_{0}^{t_\eta}  f_{\lambda}^2(t - s) \varsigma^2_{\lambda, c}(s) \, ds + (l_\infty - \eta) \int^t_{t_\eta} f_{\lambda}^2(t - s) \, ds = \int_{0}^{t_\eta}  f_{\lambda}^2(t - s) \varsigma^2_{\lambda, c}(s) \, ds + (l_\infty - \eta)\int_0^{t -t_\eta} f_{\lambda}^2(s) \, ds.
				\]
			}
			Hence, we obtain:
			
			\centerline{$
				c\lambda^2 (1-\ell_\infty^2) = \lim_{t \to +\infty} (f_{\lambda}^2 \ast \varsigma^2_{\lambda, c})(t) \geq (l_\infty - \eta) \int_0^{+\infty} f_{\lambda}^2(u) \, du \implies l_\infty \leq \frac{c\lambda^2(1-\ell_\infty^2)}{\int_0^{+\infty} f^2_{\lambda}(s) \, ds} \; \textit{by letting} \; \eta \to 0.
				$}
			On the other hand, we also have:
			{\small
				\[
				(f_{\lambda}^2 * \varsigma^2_{\lambda, c})(t) \leq
				\int_{0}^{t_\eta}  f_{\lambda}^2(t - s) \varsigma^2_{\lambda, c}(s) \, ds + (\ell_\infty + \eta) \int^t_{t_\eta} f_{\lambda}^2(t - s) \, ds = \int_{0}^{t_\eta}  f_{\lambda}^2(t - s) \varsigma^2_{\lambda, c}(s) \, ds + (\ell_\infty + \eta)\int_0^{t -t_\eta} f_{\lambda}^2(s) \, ds.
				\]
			}
			Therefore, we obtain:
			
			\centerline{$
				c\lambda^2(1-\ell_\infty^2)  = \lim_{t \to +\infty} (f_{\lambda}^2 \ast \varsigma^2_{\lambda, c})(t) \leq (l_\infty + \eta) \int_0^{+\infty} f_{\lambda}^2(u) \, du \implies l_\infty \geq \frac{c\lambda^2(1-\ell_\infty^2)}{\int_0^{+\infty} f^2_{\lambda}(s) \, ds} \quad \textit{as} \quad \eta \to 0.
				$}
			\noindent This completes the proof and we are done. \hfill $\square$
			
			\bigskip
			\noindent {\bf Proof of Proposition~\ref{prop:equiv}.}
			We adapt the proof of the corresponding Proposition in \cite{Pages2024} in order to prove the equivalence of the two statements or claims above.
			
			\noindent \fbox{$(i) \Rightarrow (ii)$} Assume  $  {\rm Var}(X_t) = {\rm Var}(X_0) = v_0$ for every  $t\ge 0$. If $v_0=0$, then $X_t = m_0$ $a.s.$ for every $t\ge 0$. Consequently $\E\, \sigma^2(X_t ) =  \E\, \sigma^2\big(\frac{\mu_\infty}{\lambda}\big)$ which is constant over time anyway. 
			Assume that \( v_0 > 0 \). Then, since the constant \( c \) is finite, it follows that \( \mathbb{E}[\sigma^2(X_0)] > 0 \). Define the function
			\begin{equation}
				g(t) = \varsigma^2_{\lambda,c} \left( \frac{\mathbb{E}[\sigma^2(X_t)]}{\mathbb{E}[\sigma^2(X_0)]} - 1 \right).
			\end{equation}
			We can check using equation~\eqref{eq:var} and \((E_{\lambda,c})\) that this function satisfies the convolution equation \( f^2_{\lambda} * g = 0 \) with the initial condition \( g(0) = 0 \).
			Furthermore, under the assumption that \( \sigma \) has linear growth, the expectation \( \mathbb{E}[\sigma^2(X_t)] \) remains bounded due to the boundedness of \( \mathbb{E}[X_t^2] \). As a consequence, the function \( g \), along with its positive and negative parts \( g^+ \) and \( g^- \), admits a Laplace transform.
			
			Since the Laplace transform of \( f^2_{\lambda} \), denoted \( L_{f^2_{\lambda}} \), is not identically zero and is strictly positive on \( (0, +\infty) \), we obtain:
			
			\centerline{$
				L_{f^2_{\lambda}} \cdot L_{g^+} = L_{f^2_{\lambda}} \cdot L_{g^-}.
				$}
			This implies \( L_{g^+} = L_{g^-} \), hence \( g^+ = g^- \), and consequently \( g = 0 \).
			\fbox{$(ii)\Rightarrow (i)$}  First, we have that $\bar \sigma^2_0 =\bar \sigma^2_t = \E\, \sigma^2(X_t)$, $t\ge 0$, so that  it follows from Equation~\eqref{eq:var} and \( (E_{\lambda,c}) \).
			\begin{align*}
				{\rm Var}(X_t ) &= {\rm Var}(X_0 )\big(1-\frac{(f_\lambda * \mu)(t)}{\lambda m_0}\big)^2 + \frac{\bar \sigma^2_0}{\lambda^2} (f^2_{\lambda}*\varsigma_{\lambda,c}^2)(t)  = {\rm Var}(X_0 )\big(1-\frac{(f_\lambda * \mu)(t)}{\lambda m_0}\big)^2 + \frac{v_0}{c\lambda^2} (f^2_{\lambda}*\varsigma_{\lambda,c}^2)(t) \\
				&\overset{(E_{\lambda,c})}{\operatorname*{=}}  v_0 \big(1-\frac{(f_\lambda * \mu)(t)}{\lambda m_0}\big)^2 + v_0 \big(1- \big(1-\frac{(f_\lambda * \mu)(t)}{\lambda m_0}\big)^2 \big)  = v_0. 
			\end{align*}
			This completes the proof.  \hfill $\Box$
				
		\medskip
		\noindent {\bf Proof of Proposition \ref{prop:Momentctrl}.} Using equation~\ref{eq:Volterrameanrevert_}
	and owing to~\eqref{eq:CondMean}, we have: $  \forall\, t\ge 0$,
		{\small
			\begin{align*}
				X_t - m_0 &= \Big(X_0-m_0\Big)\big(1-\frac{(f_\lambda * \mu)(t)}{\lambda m_0}\big) + m_0\Big(\big(1-\frac{(f_\lambda * \mu)(t)}{\lambda m_0}\big)-1\Big)+  \frac{1}{\lambda}(f_{\lambda} * \mu)_t + \frac{1}{\lambda} \left(f_{\lambda} \stackrel{W}{*} \varsigma(\cdot)\sigma(X_{\cdot})\right)_t\\
				&= \Big(X_0-m_0\Big)\big(1-\frac{(f_\lambda * \mu)(t)}{\lambda m_0}\big) + \frac{1}{\lambda} \int_0^t f_{\lambda}(t-s)\varsigma(s) \sigma(X_{ s}) \, dW_s.
			\end{align*}
		}
		so that in particular, $\quad \Big| \E\, \big(X_t\big)-m_0\Big| \le   \Big|1 - (f_{\lambda} * \frac{\mu}{\lambda m_0})_t\Big| \Big| \E\,\big(X_0\big)-m_0 \Big|$.\\
		$(a)$
		Using elementary computations and It\^o's Isometry show that for every $t\ge 0$
		\begin{align*}
			\E\, \Big(\Big|X_t-m_0\Big|\Big)^2 &\le \E\, \Big(\Big|(X_0-m_0)\big(1-\frac{(f_\lambda * \mu)(t)}{\lambda m_0}\big) \Big|^2\Big) + \frac{1}{\lambda^2} \int_0^t f^2_{\lambda}(t-s) \varsigma^2(s) \E\, \Big( \sigma^2(X_{s})\Big) ds.
		\end{align*}
		Set \( \rho := c\, [\sigma]^2_{\text{Lip}} \in (0,\, 1)\) and let \(\epsilon\) in Remark \ref{rem:on_sigma} be equal to \(\epsilon = \frac{\rho}{\eta}\) where \(\eta \in (0,\, 1 - \rho)\) is  a free parameter such that \( \rho + \eta \in (0,\, 1) \).
		From equation ~\eqref{eq:on_sigma},
		The real constants \( \kappa_i,\, i = 0, 2 \) depending on \( \eta \) and given by
		\(k_0 = k_0 (\eta):=(1+\frac{\rho}{\eta}) |\sigma(m_0)|^2\) and \(k_2= k_2(\eta):=(1+\frac{\eta}{\rho})[\sigma]^2_{\text{Lip}}\) so that \(c\, \kappa_2 = \rho + \eta < 1.\)
		
		Next, we have using equations ~\eqref{eq:on_sigma} and ~\eqref{eq:VolterraStabilizer} ($f^2_{\lambda}*\varsigma^2 = c \lambda^2(1-\big(1-\frac{(f_\lambda * \mu)}{\lambda m_0}\big)^2)$):
		\begin{align*}
			\E\, \Big( |X_t-m_0|\Big)^2 &\le \E\, \Big( |X_0-m_0|\Big)^2 \big(1-\frac{(f_\lambda * \mu)(t)}{\lambda m_0}\big)^2 + \kappa_0c\big( 1-\big(1-\frac{(f_\lambda * \mu)(t)}{\lambda m_0}\big)^2\big)\\
			&+ \frac{\kappa_2}{\lambda^2} \int_0^t f^2_{\lambda}(t-s) \varsigma^2(s) \E\, \Big(| X_{s} -m_0|^2 \Big) ds.
		\end{align*}\label{eq:variance bound}
		Now let $A> \bar A_\eta:=\frac{\kappa_0c}{1-\kappa_2c}\vee  \E\, \Big( | X_0-m_0|\Big)^2$, $\delta>0$ and 
		$
		t_{\delta} = \inf\Big\{ t \geq 0 :  \E\, \Big( |X_t-m_0|\Big)^2 \ge A + \delta\Big\}.
		$
		As $t\mapsto \E\, \Big( |X_t-m_0|\Big)^2$  is continuous  and $A>\E\, \Big( |X_0-m_0|\Big)^2$ it follows from the above inequality  and the identity satisfied by $\varsigma$ (Equation~\eqref{eq:VolterraStabilizer}) 
		 that, if $t_{\delta}<+\infty$,  then $ \E\, \Big(| X_s -m_0| \Big)^2 < A + \delta \quad  \forall s \leq t_{\delta} $ and we have:
		\(	A+\delta = \E\, \Big(| X_{t_\delta}-\frac{\mu_\infty}{\lambda}|\Big)^2< A \big(1-\frac{(f_\lambda * \mu)(t_\delta)}{\lambda m_0}\big)^2(t_\delta) + \big(\kappa_0 c +\kappa_2c (A+\delta)\big)\big(1-\big(1-\frac{(f_\lambda * \mu)(t_\delta)}{\lambda m_0}\big)^2(t_\delta)\big).\)
		
		Now, as $\kappa_0 c +\kappa_2c A<A$ by construction of $A$, we have:
		{\small
		\begin{align*}
		A+\delta = \E\, \Big(| X_{t_{\delta}}-m_0|\Big)^2 &< A\big(1-\frac{(f_\lambda * \mu)(t_\delta)}{\lambda m_0}\big)^2(t_\delta) +A(1-\big(1-\frac{(f_\lambda * \mu)(t_\delta)}{\lambda m_0}\big)^2(t_\delta)) + \kappa_2c \delta\big((1-\big(1-\frac{(f_\lambda * \mu)(t_\delta)}{\lambda m_0}\big)^2\big)\\
		 &< A + \delta.
		\end{align*}
		}
		As c is
		so that $c\kappa_2<1$. This yields a contradiction. Consequently, $t_{\delta}= +\infty$ which implies that $\E\, \Big( |X_t-m_0|\Big)^2 \le A+\delta$ for every $t\ge 0$. Letting $\delta \to 0$ and 
		$A\to \bar A_\eta$ successively, yields
		 \begin{equation}
		 	\sup_{t\ge 0}\E\, \Big(| X_t-m_0|^2\Big)\le \bar A_\eta = \E\, \Big(| X_0-m_0|^2\Big) \vee \frac{c(1+\frac{\rho}{\eta}) |\sigma(m_0)|^2}{1-(1+\frac{\eta}{\rho})c[\sigma]^2_{\text{Lip}}} = \E\, \Big(| X_0-m_0|^2\Big) \vee c|\sigma(m_0)|^2\frac{\eta+\rho}{\eta(1-\rho-\eta)}.
		 \end{equation}
		\noindent A straightforward computation shows that 
		\(\eta \mapsto \bar{A}_{\eta}\)
		attains its minimum on \( (0,\,1 - \rho) \)
		at \(\eta = \sqrt{\rho} - \rho\).  This minimum is given by \(\frac{c}{(1-\sqrt{\rho})^2}|\sigma(m_0)|^2\) which completes the proof of the stated result.
		
		\smallskip
		\noindent$(b)$ Let $p\ge 2$. Set \( \rho_p := c\, (C_p^{BDG})^2\, [\sigma]^2_{\mathrm{Lip}} \in (0,\, 1)\). Owing to the triangle inequality and applying the   {\em BDG} inequality to the (a priori) local martingale $M_u = \int_0^u f_{\lambda}(t-s)\varsigma(s)\sigma(X_{s})dW_s$, $0\le s\le t$, (see~\cite[Proposition~4.3]{RevuzYor}) follow by the generalized Minkowski inequality, we get:
		\begin{align*}
			\Big\| |X_t-m_0|\Big\|_p &\le \Big\|  |X_0-m_0|\Big\|_p \Big|1-\frac{(f_\lambda * \mu)(t)}{\lambda m_0}\Big|  + \frac{C_p^{BDG}}{\lambda}\Big\|\left(f_{\lambda}^2 \stackrel{dt}{*} \varsigma^2(\cdot)|\sigma(X_{\cdot })|^2\right)_t \Big\|_{\frac p2}^{\frac 12}\\
			&\le  \Big\|  |X_0-m_0|\Big\|_p \Big|1-\frac{(f_\lambda * \mu)(t)}{\lambda m_0}\Big| +  \frac{C_p^{BDG}}{\lambda}\Big( \int_0^t f^2_{\lambda}(t-s)\varsigma^2(s) \big\||\sigma(X_{s})|^2\big\|_{\frac p2}\Big)^{\frac 12}.
		\end{align*}
		Owing to the elementary inequality $(a+b)^2 \le (1+\frac1\epsilon) a^2+(1+\epsilon )b^2\; \forall \epsilon\!\in (0,1/\rho_p-1)$, it follows that:
		{\small
			\begin{align*}
				\Big\| |X_t-m_0|\Big\|_p^2 &\le \Big\|  |X_0-m_0|\Big\|_p^2 \Big|1-\frac{(f_\lambda * \mu)(t)}{\lambda m_0}\Big|^2(1+1/\epsilon)
				+ \frac{(C^{BDG}_p)^2}{\lambda^2}(1+\epsilon) \int_0^t f^2_{\lambda}(t-s)\varsigma^2(s) \||\sigma(X_{s})|^2\|_{\frac p2}ds
			\end{align*}
		}
		Likewise, set \( \widetilde{\rho}_p := c\, (C_p^{BDG})^2\, [\sigma]^2_{\mathrm{Lip}}(1 + \varepsilon)= \rho_p(1 + \varepsilon) \in (0,\, 1)\) and let \(\epsilon\) in Remark \ref{rem:on_sigma} be equal to \( \frac{\widetilde{\rho}_p}{\eta}\) where \(\eta \in (0,\, 1 - \widetilde{\rho}_p)\) is  a free parameter such that \( \widetilde{\rho}_p + \eta \in (0,\, 1) \).
		From equation ~\eqref{eq:on_sigma},
		The real constants \( \kappa_i,\, i = 0, 2 \) depending on \( \eta \) are given by
		\(k_0 = k_0 (\eta):=(1+\frac{\widetilde{\rho}_p}{\eta}) |\sigma(m_0)|^2\) and \(k_2= k_2(\eta):=(1+\frac{\eta}{\widetilde{\rho}_p})[\sigma]^2_{\text{Lip}}\) so that \(c\, (C_p^{BDG})^2\,(1 + \varepsilon)\, \kappa_2 = \widetilde{\rho}_p + \eta < 1.\)
		As $\frac p2\ge 1$, according to the remark \ref{rem:on_sigma}, \(\Big\||\sigma(X_{s})|^2\Big\|_{\frac p2} \le \kappa_0 +\kappa_2\Big\| |X_{ s}-m_0| \Big\|^2_p\)
		which entails, combined with the identity  $f^2_{\lambda}*\varsigma^2 = c \lambda^2 (1-\big(1-\frac{(f_\lambda * \mu)}{\lambda m_0}\big)^2)$,   that, for every $t\ge 0$, 
		{\small
			\begin{align}\nonumber 
				&\Big\| |X_t-m_0|\Big\|_p^2 \le  \Big\|  |X_0-m_0|\Big\|_p^2 \Big|1-\frac{(f_\lambda * \mu)(t)}{\lambda m_0}\Big|^2 (1+\frac1\epsilon) \\ \label{eq:variancebound2}
				&\hspace{1.6cm}+ (C^{BDG}_p)^2(1+\epsilon)\Big( \kappa_0c\big( 1-\big(1-\frac{(f_\lambda * \mu)(t)}{\lambda m_0}\big)^2 \big) + \frac{\kappa_2}{\lambda^2} \int_0^t f^2_{\lambda}(t-s) \varsigma^2(s)\Big\| |X_{ s}-m_0|\Big\|^2_p
				ds\Big).
			\end{align}
		}
		Now let $A> \bar A_{\eta,\epsilon}:=\frac{\kappa_0c(C^{BDG}_p)^2(1+\epsilon)}{1-\kappa_2c(C^{BDG}_p)^2(1+\epsilon)}\vee  \Big[ (1+1/\epsilon)\Big\|  |X_0-m_0|\Big\|_p^2\Big]$, $\delta>0$ and \(t_{\delta} = \inf\Big\{ t \geq 0: \Big\| |X_t-m_0|\Big\|_p^2 \ge A + \delta\Big\}.\)
		If $t_{\delta}<+\infty$, then, on the one hand, it follows from the continuity of $t\mapsto\Big\| |X_t-m_0|\Big\|_p^2$   that    $A+\delta =\Big\| |X_{ t_\delta}-m_0|\Big\|_p^2$ and, on the other hand, from  Equation~\eqref{eq:VolterraStabilizer} satisfied by $\varsigma$,
		that \(\forall t\leq t_\delta \)
		\(\int_0^t f^2_{\lambda}(t-s) \varsigma^2(s)\Big\| |X_{s}-m_0|\Big\|_p^2ds\le  (A+\delta)c \lambda^2(1-\big(1-\frac{(f_\lambda * \mu)(t)}{\lambda m_0}\big)^2).\quad \)
		Moreover, since \( A > \left\| X_0 - m_0\right\|_p^2 (1 + \frac{1}{\epsilon}) \), we deduce from~\eqref{eq:variancebound2} the inequalities:
		\begin{align*}
			A+\delta &= \Big\| |X_{ t_\delta}-m_0|\Big\|_p^2< A \big(1-\frac{(f_\lambda * \mu)(t_\delta)}{\lambda m_0}\big)^2+(C^{BDG}_p)^2(1+\epsilon)\big(\kappa_0 c +\kappa_2c (A+\delta)\big)\big(1-\big(1-\frac{(f_\lambda * \mu)(t_\delta)}{\lambda m_0}\big)^2\big)\\
			&< A \big(1-\frac{(f_\lambda * \mu)(t_\delta)}{\lambda m_0}\big)^2 +  A(1-\big(1-\frac{(f_\lambda * \mu)(t_\delta)}{\lambda m_0}\big)^2) +   (C^{BDG}_p)^2(1+\epsilon)c\kappa_2 \delta\big((1-\big(1-\frac{(f_\lambda * \mu)(t_\delta)}{\lambda m_0}\big)^2\big)\\
			&< A + \delta \big((1-\big(1-\frac{(f_\lambda * \mu)(t_\delta)}{\lambda m_0}\big)^2\big) < A+ \delta.
		\end{align*}
		Here, the second inequality uses the bound  \( (C^{\mathrm{BDG}}_p)^2 (1 + \epsilon)\, c\, (\kappa_0 + \kappa_2 A) < A \)
		which holds by the very definition of \( A \), while the penultimate inequality follows from the assumption that $ (C^{BDG}_p)^2(1+\epsilon)c\kappa_2<1$. This  yields a contradiction. Consequently, $t_{\delta}= +\infty$ which implies that $\Big\|  |X_{ t}-m_0|\Big\|_p^2  \le A+\delta$ for every $t\ge 0$. Letting $\delta \to 0$ and 
		$A\to \bar A_{\eta,\epsilon}$ successively, yields 
		\begin{equation}
				\sup_{t\ge 0}\Big\|  |X_t-m_0|\Big\|_p\leq \bar A^{\frac1 2}_{\eta,\epsilon}= (1+\frac1\epsilon)\Big\|  |X_0-m_0|\Big\|_p \vee\Big( \frac{\widetilde{\rho}_p}{[\sigma]^2_{\text{Lip}}}|\sigma(m_0)|^2\frac{\eta+\widetilde{\rho}_p}{\eta(1-\widetilde{\rho}_p-\eta)} \Big)^{\frac1 2}<+\infty.
		\end{equation}
		\noindent A straightforward computation shows that the mapping 
		\(\eta \mapsto \bar{A}_{\eta,\epsilon}\) attains its minimum on the interval 
		\((0,\,1 - \widetilde{\rho}_p)\) at \(\eta = \sqrt{\widetilde{\rho}_p} - \widetilde{\rho}_p\), this minimum being 
		\(\frac{\widetilde{\rho}_p}{[\sigma]^2_{\text{Lip}}(1 - \sqrt{\widetilde{\rho}_p})^2} 
		\big|\sigma(m_0)\big|^2 = \frac{c(C^{BDG}_p)^2(1+\epsilon)}{(1 -[\sigma]_{\text{Lip}} \sqrt{c(C^{BDG}_p)^2(1+\epsilon)})^2} 
		\big|\sigma(m_0)\big|^2\), 
		which completes the proof.
		The stated results follows by setting \(C^{BDG}_p = 2\sqrt{p}\) owing to Lemma \ref{lm:best_bdg_constant}.\hfill $\square$
		
	    	\bigskip
	    \noindent {\bf Proof of Proposition \ref{prop:alphaFractKernel1_}.}
	    	\noindent Taking the Laplace transform of \eqref{eq:Deriv_stabil} and arguing as in Section~\ref{sec:sigma2rough}, we obtain, $L_{e^{2\rho\cdot}g'_{\alpha,\rho,\lambda}}(t) \stackrel{+\infty}{\sim} 2\lambda \,c \frac{\Gamma(\alpha)^2}{\Gamma(2\alpha-1)} t^{-(1-\alpha)}.$ By Karamata's Tauberian theorem, it follows that $g'_{\alpha,\rho,\lambda}(t) \stackrel{t\to0}{\sim} \frac{2\lambda c \Gamma(\alpha)^2}{\Gamma(2\alpha-1) \Gamma(2-\alpha)}  e^{-2 \rho t}t^{-\alpha}$.
	    Consequently, \(\lim_{t\to 0^+}g'_{\alpha,\rho,\lambda}(t)=+\infty\). One also verifies that \(\lim_{t\to 0^+}g''_{\alpha,\rho,\lambda}(t)=-\infty\) and following the argument of \cite[Proof of Proposition 5.3(c)]{Pages2024}, one then proves by contradiction that $g^{''}_{\alpha,\rho,\lambda} \le 0$ on $(0,+\infty)$, i.e. \(g_{\alpha,\rho,\lambda}\) is concave.
	    Moreover,  by Tauberian Final Value Theorem if $\lim_{t \to +\infty} g'_{\alpha,\rho,\lambda}(t)$ exists, then
	    $$\lim_{t \to +\infty} g'_{\alpha,\rho,\lambda}(t) = \lim_{z \to 0} z  L_{g'_{\alpha,\rho,\lambda}}(z) = \lim_{z \to 0} \left(z^2  L_{g_{\alpha,\rho,\lambda}}(z) -zg_{\alpha,\rho,\lambda}(0) \right)= \lim_{z \to 0} \left(z^2  L_{g_{\alpha,\rho,\lambda}}(0) -zg_{\alpha,\rho,\lambda}(0) \right)= 0$$ since \(g_{\alpha,\rho,\lambda} \) is integrable and thus have a finite Laplace transform.
	    Finally, $\lim_{t\to0^+}g'_{\alpha,\rho,\lambda}(t) = +\infty $ , \(\lim_{t \to +\infty} g'_{\alpha,\rho,\lambda}(t) = 0\) and \(g'_{\alpha,\rho,\lambda}\) is non-increasing on \((0,+\infty)\) (\(g''_{\alpha,\rho,\lambda}\leq 0\)), it follows that \(g'_{\alpha,\rho,\lambda}(t) \geq 0 \quad \forall t \in (0,+\infty)\).
	    Hence \(g_{\alpha,\rho,\lambda}\) is concave, non-decreasing and non-negative on \((0,+\infty)\), so that the function \( \sqrt{g_{\alpha,\rho,\lambda}} \)  is well-defined on \( (0, +\infty) \).     \hfill $\square$
	    
	    \bigskip
	    \noindent {\bf Proof of Proposition \ref{prop:main2}.}
	       From Equation~\eqref{eq:DerivSolventGammaKernel}, we have for all \( t > 0 \), \( f_{\alpha,\rho,\lambda}(t) := e^{-\rho t} f_{\alpha,\lambda}(t) \).  Thus, Proposition~\ref{prop:main2}~$(b)$ follows directly from \cite[Proposition~5.1]{Pages2024} for the \(\alpha-\)fractional kernel, upon multiplication by the exponential damping factor \(t\to e^{-\rho t}\), which is smooth, completely monotone and bounded by one. Therefore, \(f_{\alpha,\rho,\lambda}\) is CM and infinitely differentiable as the 
	       product of such functions. Moreover, the ${\cal L}^{2\beta}$-integrability of \(f_{\alpha,\rho,\lambda}\) $\forall \beta \in \big( 0, \frac{1}{2(1-\alpha)}\big)$
	       follows immediately from that of \(f_{\alpha,\lambda}\).
	       As for the ${\cal L}^2(\R_+)$-$\vartheta$-H\"older continuity of $f_{\alpha,\rho,\lambda}$ in~$(c)$, let $\delta >0$,
	       using that $0\le 1-e^{-v} \le (1-e^{-v})^{\vartheta} \le {v}^{\vartheta}$, for every $v\ge 0$, and $\vartheta\!\in (0,1]$, one may deduce
	       {\small
	       \begin{align*}
	       	\Big( \int_0^\infty \left| f_{\alpha,\rho,\lambda}(t+\delta) - f_{\alpha,\rho,\lambda}(t) \right|^2 dt \Big)^{\frac12}
	       	&\leq e^{-\rho \delta}\Big( \int_0^\infty \left| f_{\alpha, \lambda}(t+\delta) - f_{\alpha, \lambda}(t) \right|^2 dt \Big)^{1/2} + (\rho \delta )^\vartheta \Big( \int_0^\infty \left| f_{\alpha, \lambda}(t) \right|^2 dt \Big)^{\frac12} \\
	       	&\leq  e^{-\rho \delta} C_{\vartheta, \lambda}\delta^{\vartheta} + C_{f_\lambda}\delta^{\vartheta} := C_{\vartheta,\rho, \lambda}\delta^{\vartheta}.
	       \end{align*}
	       }
	      where in the last inequality, we used that \(f_{\alpha, \lambda} \in \mathcal{L}^2(\R_+)\) together with \cite[Proposition 5.1]{Pages2024}. 
	      For the claim~$(a)$, the identity \(R'_{\alpha, \rho, \lambda}(t) = - f_{\alpha, \rho, \lambda}(t)\), together with the complete monotonicity of \(f_{\alpha, \rho, \lambda}\), implies that\( R_{\alpha, \rho, \lambda} \) is completely monotone, hence  \( \mathcal{C}^\infty \) on  \( (0, +\infty) \), and decreases to the given limit computed in
	      the begining of the section~\ref{subsec:Gammaalphafrac}.
	    This completes the proof and we are done.	 \hfill $\square$

\end{document}